\documentclass[11pt,a4paper]{article}
\usepackage[utf8]{inputenc}
\usepackage[english]{babel}
\usepackage[T1]{fontenc}
\usepackage{amsmath, amsfonts, amssymb, amsthm}
\usepackage{lmodern}
\usepackage{fullpage}
\usepackage{microtype}
\usepackage{enumitem}
\usepackage{mathtools}
\usepackage[colorlinks=true, linkcolor=blue]{hyperref}
\usepackage{empheq}

\numberwithin{equation}{section}

\newtheorem{theorem}[equation]{Theorem}
\newtheorem*{theorem_nnum}{Theorem}
\newtheorem{proposition}[equation]{Proposition}
\newtheorem{lemma}[equation]{Lemma}
\newtheorem{corollary}[equation]{Corollary}

\theoremstyle{definition}
\newtheorem{definition}[equation]{Definition}

\theoremstyle{remark}
\newtheorem{remark}[equation]{Remark}
\newtheorem{example}[equation]{Example}

\renewcommand{\H}{R}
\newcommand{\tuple}{\boldsymbol}
\newcommand{\I}{\mathcal{I}}
\newcommand{\J}{\mathcal{J}}

\newcommand{\sizeglob}[1]{n_{#1}}
\newcommand{\sizepart}[2]{n_{#1}(#2)}

\newcommand{\ideal}{\mathfrak{I}}
\newcommand{\polrep}{\varphi}
\newcommand{\prep}[1]{{#1}'}
\newcommand{\idealB}{\mathfrak{J}}

\newcommand{\multgen}{}

\newcommand{\multx}{\multgen{x}}
\newcommand{\multone}{\multgen{\mathbf{1}}}
\newcommand{\actpol}[2]{\prescript{#1}{}{#2}}
\newcommand{\idealD}{\mathfrak{K}}
\newcommand{\prof}{\mathrm{Prof}}

\newcommand{\orbitsS}[2][\mathfrak{S}]{I^{#2} / #1_{#2}}

\allowdisplaybreaks[3]

\title{Morita equivalences for cyclotomic Hecke algebras of type~B and~D
\\
\medskip Équivalences de Morita pour les algèbres de Hecke cyclotomiques de type~B et~D
}
\author{Loïc \textsc{Poulain d'Andecy}\footnote{The first author is supported by \emph{Agence Nationale de la Recherche} through the JCJC project ANR-18-CE40-0001.}$\ $\thanks{Laboratoire de Math\'ematiques de Reims UMR 9008, Universit\'e de Reims Champagne-Ardenne,
Moulin de la Housse BP 1039, 51100 Reims, France} \and Salim \textsc{Rostam}\thanks{Univ Rennes, CNRS, IRMAR - UMR 6625, F-35000 Rennes, France}}

\date{}

\begin{document}

\maketitle

% j ai rajoute une phrase dans l'abstract
\abstract{We give a Morita equivalence theorem for so-called cyclotomic quotients of affine Hecke algebras of type~B and D, in the spirit of a classical result of Dipper--Mathas in type A for Ariki--Koike algebras. Consequently, the representation theory of affine Hecke algebras of type B and D reduces to the study of their cyclotomic quotients with eigenvalues in a single orbit under multiplication by $q^2$ and inversion. The main step in the proof consists in a decomposition theorem for generalisations of quiver Hecke algebras that appeared recently in the study of affine Hecke algebras of type B and~D. This theorem reduces the general situation of a disconnected quiver with involution to a simpler setting. To be able to treat types B and D at the same time we unify the different definitions of quiver Hecke algebra for type B that exist in the literature.}
{
\renewcommand{\abstractname}{Résumé}
\abstract{Nous énonçons un théorème d'équivalence de Morita pour les quotients cyclotomiques des algèbres de Hecke affines de type~B et~D, suivant un résultat classique de Dipper--Mathas en type A pour les algèbres d'Ariki--Koike. Ainsi, la théorie des représentations des algèbres de Hecke affines de type~B et~D se réduit à l'étude de leurs quotients cyclotomiques où les valeurs propres sont dans une unique orbite pour la multiplication par $q^2$ et l'inversion. La preuve consiste notamment en un théorème de décomposition pour des généralisations d'algèbres de Hecke carquois introduites récemment dans l'étude des algèbres de Hecke affines de type~B et~D, ramenant la situation générale d'un carquois non connexe avec involution à un cadre plus simple. Pour traiter simultanément les deux types, nous unifions les différentes définitions d'algèbres de Hecke carquois pour le type~B déjà existantes.}
}
\section*{Introduction}

Cyclotomic quotients of the affine Hecke algebra of type A, also known as Ariki--Koike algebras, have been extensively studied since their introduction by Broué--Malle~\cite{broue-malle} and Ariki--Koike~\cite{ariki-koike}. Given a field $K$, a subset $I \subseteq K^\times$, an element $q \in K^\times$ and a finitely-supported family $\Lambda = (\Lambda_i)_{i \in I}$ of non-negative integers, the Ariki--Koike algebra $H^\Lambda(\mathfrak{S}_n)$ is defined by the generators $g_0, \dots, g_{n-1}$ and the relations
\begin{align*}
g_i g_j &= g_j g_i, &&\text{for all } i, j \in \{0, \dots, n-1\}, \lvert i - j \rvert > 1,
\\
g_i g_{i+1} g_i &= g_{i+1} g_i g_{i+1}, &&\text{for all } i \in \{1, \dots, n-2\},
\\
g_0 g_1 g_0 g_1 &=  g_1 g_0 g_1 g_0,
\\
(g_i - q)(g_i + q^{-1}) &= 0, &&\text{for all } i \in \{1, \dots, n-1\},
\\
\prod_{i \in I} (g_0 - i)^{\Lambda_i} &= 0.
\end{align*}
We note that Ariki--Koike algebras are quotients, by the last relation, of affine Hecke algebras of type A and that the study of their representations (for all choices of $I$ and $\Lambda$) is equivalent to the study of finite-dimensional representations of affine Hecke algebras of type A.

By an important theorem of Dipper--Mathas~\cite{dipper-mathas}, we know that it suffices to study Ariki--Koike algebras when the set $I$ is \emph{$q^2$-connected}, that is,  in a single $q^2$-orbit (and even, up to a scalar renormalisation of the generator $g_0$, when $I \subseteq \langle q^2\rangle$). More precisely, if $I = \amalg_{j = 1}^d I^{(j)}$ is the decomposition of $I$ into $q^2$-connected sets then we have a Morita equivalence
\begin{equation}
\label{equation:intro_morita}
H^\Lambda(\mathfrak{S}_n) \overset{\text{Morita}}{\simeq}
\bigoplus_{\substack{n_1, \dots, n_d \geq 0 \\ n_1 + \dots + n_d = n}} \bigotimes_{j = 1}^d H^{\Lambda^{(j)}}(\mathfrak{S}_{n_j}),
\tag{$\clubsuit$}
\end{equation}
where $\Lambda^{(j)}$ is the restriction of $\Lambda$ to $I^{(j)}$.
(Note that the assumption in~\cite{dipper-mathas} is slightly stronger than the one above, but in practice it is this condition of $q^2$-connected sets that is used.)
Hence, this Morita equivalence allows to use  results that are only known when the set $I$ is $q^2$-connected, in particular, the celebrated Ariki's categorification theorem~\cite{ariki} that computes the decomposition numbers of Ariki--Koike algebras in terms of the canonical basis of a certain highest weight module over an affine quantum group.

\medskip
Another way to obtain this Morita equivalence was given by the second author~\cite[\textsection 3.4]{rostam2}, using the theory of quiver Hecke algebras. This is a family of graded algebras that was introduced a few years ago independently by Khovanov--Lauda~\cite{khovanov-lauda_diagrammaticI,khovanov-lauda_diagrammaticII} and Rouquier~\cite{rouquier_2}, in the context of categorification of quantum groups. If $\Gamma$ is a quiver, we denote by $\H_n(\Gamma)$ the associated quiver Hecke algebra (see~\textsection\ref{subsection:definition_klr}). For a certain quiver $\Gamma$ depending only on the order of $q^2$, Brundan--Kleshchev~\cite{BK} and independently Rouquier~\cite{rouquier_2} proved that a certain ``cyclotomic'' quotient of $\H_n(\Gamma)$ is isomorphic to an Ariki--Koike algebra. This result is now a basic tool in the study of Ariki--Koike algebras and their degenerations, including the symmetric group and the classical Hecke algebra of type A.  For instance, as  consequences first the Ariki--Koike algebra inherits the grading of the cyclotomic quiver Hecke algebra, and second depends on $q$ only through its order in $K^\times$.
Now if $\Gamma$ is of the form $\Gamma = \amalg_{j = 1}^d \Gamma^{(j)}$ where each $\Gamma^{(j)}$ is a full subquiver, it was shown in~\cite[\textsection 6]{rostam} that we have a decomposition
\begin{equation}
\label{equation:intro_dec_klr}
\H_n(\Gamma) \simeq \bigoplus_{\substack{n_1, \dots, n_d \geq 0 \\ n 1 + \dots + n_d = n}} \mathrm{Mat}_{\binom{n}{n_1, \dots, n_d}}\left(\bigotimes_{j = 1}^d \H_{n_j}(\Gamma^{(j)})\right).
\tag{$\spadesuit$}
\end{equation}
This isomorphism of algebras is compatible with cyclotomic quotients, and combining with the previous isomorphism of Brundan--Kleshchev and Rouquier allows to recover the Morita equivalence~\eqref{equation:intro_morita}. This Morita equivalence has been further generalised for the cyclotomic Hecke algebras of type $G(r,p,n)$ \cite{HM}. We indicate also the paper \cite{HZ} where the Dipper--Mathas result is studied and derived again from the point of view of affine Hecke algebras, and where the question of a similar result for other affine Hecke algebras is evoked.

\bigskip
The main point of this paper is to prove a similar decomposition theorem for some generalisations of quiver Hecke algebras and hence obtain an analogue of the Dipper--Mathas Morita equivalence for cyclotomic quotients of affine Hecke algebras of type B and D. Such generalisations of quiver Hecke algebras were introduced by Varagnolo and Vasserot~\cite{varagnolo-vasserot_canonical} (for type~B) and together with Shan~\cite{SVV} (for type D), in the course of their proofs of conjectures by Kashiwara--Enomoto \cite{EK} and Kashiwara--Miemietz \cite{KM}.  These algebras play for certain subcategories of representations of affine Hecke algebras of type B and D a similar role as quiver Hecke algebras for affine Hecke algebras of type A. Inspired by their results, the first author together with Walker~\cite{poulain_dandecy-walker_B,poulain_dandecy-walker_D} obtained an isomorphism theorem \textit{à la} Brundan--Kleshchev between cyclotomic quotients of affine Hecke algebras of type B and D and certain generalisations of cyclotomic quiver Hecke algebras.

The first step of this paper is to provide a definition of these generalisations of quiver Hecke algebras for the type B which encompasses all the slightly different versions previously defined. They are $\mathbb{Z}$-graded algebras and they depend upon a quiver with an involution and certain weight functions on the vertices. As for the type A case, that is, for usual quiver Hecke alegbras, the algebra that we define admits a PBW basis and this is a key ingredient to prove the decomposition theorem when the underlying quiver has several connected components.  The point of having defined a new algebra in Section~\ref{section:interpolating} is that we can now use the main results of~\cite{poulain_dandecy-walker_B,poulain_dandecy-walker_D} at the same time. We deduce our main theorem for type~B, Theorem~\ref{theorem:morita_B}, that we state now. 
Write $I \subseteq K^\times$ as $I = \amalg_{j = 1}^d I^{(j)}$ such that each $I^{(j)}$ is $q^2$-connected \emph{and} stable by scalar inversion. As in the type A case, for $\Lambda = (\Lambda_i)_{i \in I} \in \mathbb{N}^{(I)}$ we denote by $H^\Lambda(B_n)$ the quotient of the affine Hecke algebra of type B by the relation
\[
\prod_{i \in I} (X_1 - i)^{\Lambda_i} = 0
\]
(see~\textsection\ref{subsection:morita_B} for a precise definition).

\begin{theorem_nnum}
We have an (explicit) isomorphism
\[
H^\Lambda(B_n) \simeq
\bigoplus_{\substack{n_1, \dots, n_d \geq 0 \\ n_1 + \dots + n_d = n}} \mathrm{Mat}_{\binom{n}{n_1, \dots, n_d}}\left( \bigotimes_{j = 1}^d H^{\Lambda^{(j)}}(B_{n_j})\right),
\]
in particular, we have a Morita equivalence
\[
H^\Lambda(B_n) \overset{\text{Morita}}{\simeq}
\bigoplus_{\substack{n_1, \dots, n_d \geq 0 \\ n_1 + \dots + n_d = n}} \bigotimes_{j = 1}^d H^{\Lambda^{(j)}}(B_{n_j}).
\]
\end{theorem_nnum}
We also deduce that a similar result holds for the cyclotomic quotient $H^\Lambda(D_n)$ of the affine Hecke algebra of type D. Some technicalities typical to the type D situation result in a formulation of the final result a bit more complicated than for type B in the Theorem above, since it involves in addition a semi-direct product by powers of a cyclic group of order 2 (see Theorem \ref{theorem:morita_D}).

One motivation for considering cyclotomic quotients of affine Hecke algebras is that the study of (finite-dimensional) representations of the affine Hecke algebra is equivalent to the study of representations of all their cyclotomic quotients. As a consequence of our main results, we obtain that, for affine Hecke algebras of type B and D, this study reduces to considering the algebras $H^\Lambda(B_n)$ and $H^\Lambda(D_n)$ when the set $I$ is $q^2$-connected and stable by scalar inversion (see Corollaries~\ref{corollary:morita_B} and \ref{corollary:morita_D} for more details and a complete description of the finite number --- up to four --- of sets $I$ to be considered). This generalises the classical reduction for the affine Hecke algebras of type A (for which it is enough to consider $I=q^{2\mathbb{Z}}$) induced by the Dipper--Mathas result.
 
 % dernier paragraphe modifie et agrandi

% quelques modifs et raccourcissements ci-dessous
\paragraph{Organisation of the paper.} In Section~\ref{section:general_decomposition_theorem}, given an algebra $A$ and a set of idempotents satisfying certain properties  we prove a general decomposition theorem expressing $A$ in terms of a direct sum involving matrix algebras on idempotent truncations (Corollary~\ref{corollary:isom_A}).

Let $\Gamma$ be a (possible infinite) quiver with no $1$-loops, let $I$ be its vertex set and let $\alpha \subseteq \mathfrak{S}_n$ be  a finite union of $\mathfrak{S}_n$-orbits.
 In Section~\ref{section:toy_KLR} we recall the definition of the quiver Hecke algebra $\H_\alpha(\Gamma)$. We then review the proof, based on the general theorem from Section~\ref{section:general_decomposition_theorem}, of the decomposition isomorphism of~\cite{rostam} when $\Gamma$ has several connected components, generalising it to the case where $\Gamma$ is not necessarily finite (as it is assumed in~\cite{rostam}). In~\textsection\ref{subsection:cyclotomic_klr}, given a finitely-supported family $\Lambda$ of non-negative integers we define the cyclotomic quotient $\H_\alpha^\Lambda(\Gamma)$ of $\H_\alpha(\Gamma)$ and give the corresponding isomorphism when $\Gamma$ has several connected components.
 
Then we assume that $\Gamma$ is endowed with an involution $\theta$ and let $\beta \subseteq I^n$ be an orbit for the action of the Weyl group $B_n$ of type B and rank $n$.
 We begin Section~\ref{section:interpolating} by defining the algebra $V_\beta(\Gamma,\lambda,\gamma)$ depending in addition on $\lambda \in \mathbb{N}^I$ and $\gamma \in K^I$ satisfying certain conditions. This algebra generalises the constructions of~\cite{varagnolo-vasserot_canonical,poulain_dandecy-walker_B,poulain_dandecy-walker_D}, see Remarks~\ref{remark:coherence_VV}, \ref{remark:coherence_PAW} and \ref{remark:coherence_PAW_typeD} respectively. The algebra $V_\beta(\Gamma,\lambda,\gamma)$ is $\mathbb{Z}$-graded, and we prove in~\textsection\ref{subsection:basis_theorem} that it admits a PBW basis, using a polynomial realisation (the calculations  are postponed to Appendix~\ref{appendix_section:polynomial_realisation}). 
 
 Section~\ref{section:isomorphism B} is the heart of the paper. We prove a decomposition theorem, similar to~\eqref{equation:intro_dec_klr}, for  the algebra $V_\beta(\Gamma,\lambda,\gamma)$ when the quiver~$\Gamma$ is a disjoint union of $\theta$-stable full subquivers $\Gamma = \amalg_{j = 1}^d \Gamma^{(j)}$  (Theorem~\ref{theorem:isom_carquois_disjoints_B}). As in Section~\ref{section:toy_KLR}, we first use the results of Section~\ref{section:general_decomposition_theorem} and then prove  that some idempotent truncation of $V_\beta(\Gamma,\lambda,\gamma)$ can be expressed as a tensor product on smaller algebras involving the quivers $\Gamma^{(j)}$. Note here a technical difficulty comparing with the type A case: for $n_1 + \dots + n_d = n$, the group $\mathfrak{S}_{n_1} \times \dots\times\mathfrak{S}_{n_d}$ can be seen as a parabolic subgroup of~$\mathfrak{S}_n$ for its standard Coxeter structure, but it is no more the case for $B_{n_1} \times \dots \times B_{n_d} \subseteq B_n$, although this is still a subgroup. We prove in~\textsection\ref{subsection:cyclotomic_case_B} the cyclotomic analogue of the decomposition theorem (Corollary~\ref{corollary:isom_disjoint_quiver_cyclo_B}).
 
The shorter Section~\ref{section:QHA_D} is devoted to quiver Hecke algebras $W_\beta(\Gamma)$ for type D and their cyclotomic quotients $W^\Lambda_\beta(\Gamma)$. Using a result of~\cite{poulain_dandecy-walker_D} that expresses $W_\beta(\Gamma)$ as the subalgebra of fixed-points of a certain involutive automorphism of $V_\beta(\Gamma,0,0)$ (Proposition~\ref{prop-fixed-points}), we manage to give a decomposition isomorphism for $W_\beta(\Gamma)$ and its cyclotomic quotient when the quiver $\Gamma$ has several $\theta$-stable full subquivers (Theorem~\ref{theorem:isom_carquois_disjoints_D}). 

Finally, in Section~\ref{section:morita} we introduce the affine Hecke algebras $H(B_n)$ of type B and  $H(D_n)$ of type D, together with their cyclotomic quotients $H^\Lambda(B_n)$ and $H^\Lambda(D_n)$. We then use the analogues of Brundan--Kleshchev isomorphism theorem in types B and D from~\cite{poulain_dandecy-walker_B, poulain_dandecy-walker_D} to deduce from our disjoint quiver isomorphisms the announced Morita equivalences: Theorem~\ref{theorem:morita_B} for type B and Theorem~\ref{theorem:morita_D} for type D.

\paragraph{Acknowledgements} The authors would like to thank Ruari Walker for many interesting discussions initiating this work. The second author would like to thank Ruslan Maksimau for explaining a proof of Proposition~\ref{proposition:base_klr}. The authors are very grateful to an anonymous referee for many useful suggestions.%

\section{Decomposition in matrix algebras on idempotent truncations}
\label{section:general_decomposition_theorem}

The results in this section, or some versions of them, are probably known to specialists, but we could not find them in this precise form in the literature. So we state them in the form we need and provide complete proofs. The framework presented here encompasses several cases of proved isomorphism theorems such as in \cite{JaPo, rostam}.

Let $A$ be a unitary algebra over a ring $K$. Let $\I$ be a complete (finite) set of orthogonal idempotents, that is:
\begin{itemize}
\item for all $e \in \I$ we have $e^2 = e$;
\item for all $e, e' \in \I$, if $e \neq e'$ then $ee' = e'e = 0$;
\item we have $1 = \sum_{e \in \I} e$.
\end{itemize}

For any $e \in \I$, let $\phi_e, \psi_e \in A$ such that
\begin{subequations}
\label{subequations:phie_psie_e}
\begin{align}
\label{equation:phie_psie_e}
\phi_e \psi_e e &= e,
\\
\label{equation:e_phie_psie}
e \phi_e \psi_e &= e.
\end{align}
\end{subequations}

\begin{remark}
\label{remark:phie_psie_trivial_case}
Such elements necessarily exist, for instance $\phi_e = \psi_e = e$ for any $e \in \I$.  However, obviously this will not lead to interesting results.
\end{remark}

\begin{lemma}
\label{lemma:psie_e_phie_idemp}
For any $e \in \I$, the element $\psi_e e \phi_e$ is an idempotent.
\end{lemma}

\begin{proof}
Using~\eqref{equation:phie_psie_e}, we have
\begin{align*}
(\psi_e e \phi_e)^2
&=
\psi_e e (\phi_e \psi_e e) \phi_e
\\
&=
\psi_e e^2 \phi_e
\\
&=
\psi_e e \phi_e,
\end{align*}
as desired.
\end{proof}

Denote by $\J$ the image of the map $\begin{array}{|rcl}
\I & \longrightarrow & A
\\
e &\longmapsto& \psi_e e \phi_e\end{array}$ and write $\I_\epsilon$ for the fibre of any element $\epsilon \in \J$. We have
\[
\I_\epsilon = \{ e \in \I : \psi_e e \phi_e = \epsilon\},
\]
and
\[
%\label{equation:I_reunion_Jepsilon}
\bigsqcup_{\epsilon \in \J} \I_\epsilon = \I.
\]
By Lemma~\ref{lemma:psie_e_phie_idemp}, the set $\J$ consists of idempotents, however it is a priori not related to $\I$.

\begin{proposition}
\label{proposition:e_phie=phie_epsilon}
For any $\epsilon \in \J$ and any $e \in \I_\epsilon$ we have 
\begin{subequations}
\begin{align}
\label{equation:e_phie=phie_epsilon}
e \phi_e &= \phi_e \epsilon,
\\
\label{equation:epsilon_psie=psie_e}
\epsilon \psi_e &= \psi_e e.
\end{align}
\end{subequations}
\end{proposition}

\begin{proof}
We have
\begin{equation}
\label{equation:proof-psie_e_phie}
\psi_e e \phi_e = \epsilon,
\end{equation}
thus $(\phi_e \psi_e e) \phi_e = \phi_e \epsilon$. Using~\eqref{equation:phie_psie_e} we obtain the first equality. We also obtain $\psi_e (e \phi_e \psi_e) = \epsilon \psi_e$ from~\eqref{equation:proof-psie_e_phie} thus by~\eqref{equation:e_phie_psie} we obtain the second equality.
\end{proof}

\begin{proposition}
\label{proposition:psie_phie_epsilon}
For any $\epsilon \in \J$ and any $e \in \I_\epsilon$ we have
\begin{subequations}
\label{subequations:psie_phie_epsilon}
\begin{align}
\label{equation:psie_phie_epsilon}
\psi_e \phi_e \epsilon &= \epsilon,
\\
\label{equation:epsilon_psie_phie}
\epsilon \psi_e \phi_e &= \epsilon.
\end{align}
\end{subequations}
\end{proposition}

\begin{proof}
By~\eqref{equation:e_phie=phie_epsilon} we have $\phi_e \epsilon = e \phi_e$, thus
\[
\psi_e \phi_e \epsilon = \psi_e e \phi_e,
\]
and we conclude that~\eqref{equation:psie_phie_epsilon} holds since $\psi_e e \phi_e = \epsilon$ by definition of $\I_\epsilon$. Similarly, by~\eqref{equation:epsilon_psie=psie_e} we have
\[
\epsilon \psi_e \phi_e = \psi_e e \phi_e = \epsilon,
\]
thus~\eqref{equation:epsilon_psie_phie} holds.
\end{proof}

If $J$ is any finite set and $B$ any $K$-algebra, we denote by $\mathrm{Mat}_J(B)$ the $K$-algebra of matrices with rows and columns indexed by $J$ with entries in $B$.

\begin{definition}
For any $\epsilon \in \J$, we define the idempotent
\[
\hat\epsilon \coloneqq \sum_{e \in \I_\epsilon} e.
\]
\end{definition}

\begin{theorem}
\label{theorem:isomp_epsilon_A_epsilon}
Let $\epsilon \in \J$. We have the following isomorphism of $K$-algebras:
\[
\hat \epsilon A \hat \epsilon \simeq \mathrm{Mat}_{\I_\epsilon}(\epsilon A \epsilon).
\]
\end{theorem}

\begin{proof}
We first prove that for any $e', e \in \I_{\epsilon}$, the maps
\[
\begin{array}{c|rcl}
\theta_{e'e} :&
e'Ae &\longrightarrow& \epsilon A \epsilon M_{e'e}
\\
&
a & \longmapsto & \psi_{e'} a \phi_e M_{e'e}
\end{array},
\]
and
\[
\begin{array}{c|rcl}
\eta_{e'e} :&
\epsilon A \epsilon M_{e'e} &\longrightarrow & e'Ae
\\
&
aM_{e'e} & \longmapsto & \phi_{e'} a \psi_e,
\end{array}
\]
are well-defined and inverse isomorphism of $K$-modules. Here, we denoted by $M_{e'e} \in \mathrm{Mat}_{\I_\epsilon}(\epsilon A \epsilon)$ the matrix whose unique non-zero coefficient, which is $1$, is at row $e'$ and column~$e$. The maps $\theta_{e'e}$ and $\eta_{e'e}$ are well-defined by Proposition~\ref{proposition:e_phie=phie_epsilon}. Indeed, for any $a \in e' Ae$ then $a = e'ae$ and
\[
\psi_{e'} a \phi_e = (\psi_{e'} e') a (e \phi_e) = (\epsilon \psi_{e'}) a (\phi_e \epsilon) \in \epsilon A \epsilon,
\]
so $\theta_{e'e}$ is well-defined, and for any $a \in \epsilon A \epsilon$ then $a = \epsilon a \epsilon$ and
\[
\phi_{e'} a \psi_e = (\phi_{e'} \epsilon) a (\epsilon \psi_e) = (e'\phi_{e'}) a (\psi_e e) \in e'Ae,
\]
so $\eta_{e'e}$ is well-defined.
Now for any $a \in e' A e$ we have, using $a = e'ae$ and~\eqref{subequations:phie_psie_e},
\begin{align*}
\eta_{e'e}\bigl(\theta_{e'e}(a)\bigr)
&=
\eta_{e'e}(\psi_{e'} a \phi_e M_{e'e})
\\
&=
\phi_{e'}(\psi_{e'} a \phi_e)\psi_e
\\
&=
(\phi_{e'}\psi_{e'} e') a (e \phi_e \psi_e)
\\
&= e' a e
\\
&= a.
\end{align*}
Moreover, for any $a \in \epsilon A \epsilon$ we have, using $a = \epsilon a \epsilon$ and Proposition~\ref{proposition:e_phie=phie_epsilon},
\begin{align*}
\theta_{e'e}\bigl(\eta_{e'e}(aM_{e'e})\bigr)
&=
\theta_{e'e}(\phi_{e'} a \psi_e)
\\
&=
\psi_{e'} \phi_{e'} a \psi_e \phi_e M_{e'e}
\\
&=
(\psi_{e'} \phi_{e'} \epsilon) a (\epsilon\psi_e \phi_e) M_{e'e}
\\
&= \epsilon a \epsilon
\\
&= a.
\end{align*}

We now want to extend $\theta_{e'e}$ and $\eta_{e'e}$ to algebra isomorphisms. 
We have a direct sum decomposition
\begin{equation}
\label{equation:proof-direct_sum_decomposition}
\hat\epsilon A \hat\epsilon = \bigoplus_{e', e \in \I_\epsilon} e' A e.
\end{equation}
We define two maps
\begin{align*}
\theta_\epsilon &: \hat\epsilon A \hat\epsilon \to \mathrm{Mat}_{\I_\epsilon}(\epsilon A \epsilon),
\\
\eta_\epsilon &: \mathrm{Mat}_{\I_\epsilon}(\epsilon A \epsilon) \to \hat\epsilon A \hat\epsilon,
\end{align*}
by
\begin{align*}
\theta_\epsilon &\coloneqq \bigoplus_{e', e \in \I_\epsilon} \theta_{e'e},
\\
\eta_\epsilon &\coloneqq \bigoplus_{e',e \in \I_\epsilon} \eta_{e'e}.
\end{align*}
These two maps are inverse isomorphisms of $K$-modules. To prove that they are inverse isomorphism of $K$-algebras, it suffices to prove that $\theta_\epsilon$ is a morphism of $K$-algebras. Recalling the decomposition~\eqref{equation:proof-direct_sum_decomposition}, it suffices to prove that
\begin{equation}
\label{equation:proof-algebra_morphism}
\theta_\epsilon(a_1 a_2) = \theta_\epsilon(a_1) \theta_\epsilon(a_2),
\end{equation}
for any $a_i \in e'_i A e_i$ for any $e_i \in \I_\epsilon$. If $e_1 \neq e'_2$ then the left-hand side is zero, and so is the right-hand one since $M_{e'_1 e_1} M_{e'_2 e_2} = 0_{\mathrm{Mat}_{\I_\epsilon}(\epsilon A \epsilon)}$. Thus, we now assume that $e_1 = e'_2$. We have $a_1 = a_1 e_1$ and $a_1 a_2 = a_1 e_1 a_2 \in e'_1 A e_2$, thus using~\eqref{equation:e_phie_psie} we obtain
\begin{align*}
\theta_\epsilon(a_1 a_2)
&=
\theta_{e'_1 e_2}(a_1 a_2)
\\
&=
\psi_{e'_1} a_1 (e_1) a_2 \phi_{e_2} M_{e'_1 e_2}
\\
&=
\psi_{e'_1} a_1 (e_1 \phi_{e_1} \psi_{e_1}) a_2 \phi_{e_2} M_{e'_1 e_2}
\\
&=
(\psi_{e'_1} a_1 e_1 \phi_{e_1})(\psi_{e_1} a_2 \phi_{e_2}) M_{e'_1 e_2}
\\
&=
\left(\psi_{e'_1} a_1 \phi_{e_1} M_{e'_1 e_1}\right)\left(\psi_{e_1} a_2 \phi_{e_2} M_{e_1 e_2}\right)
\\
&=
\theta_{e'_1 e_1}(a_1) \theta_{e_1 e_2}(a_2)
\\
&=
\theta_\epsilon(a_1) \theta_\epsilon(a_2).
\end{align*}
This concludes the proof.
\end{proof}

\begin{corollary}
\label{corollary:isom_A}
Assume that for all $\epsilon, \epsilon' \in \J$ we have
\begin{equation}
\label{equation:assumption_strong_disjoint}
\epsilon \neq \epsilon' \implies \hat\epsilon A \hat\epsilon' = \{0\}.
\end{equation}
Then have the following isomorphism of $K$-algebras:
\[
A \simeq \bigoplus_{\epsilon \in \J} \mathrm{Mat}_{\I_\epsilon}(\epsilon A \epsilon).
\]
\end{corollary}

\begin{proof}
The assumption~\eqref{equation:assumption_strong_disjoint} implies that
\[
A \simeq \bigoplus_{\epsilon \in \J} \hat\epsilon A \hat\epsilon.
\]
We now use the result of Theorem~\ref{theorem:isomp_epsilon_A_epsilon}.
\end{proof}

\section{Application to quiver Hecke algebras}
\label{section:toy_KLR}

We here review and generalise the decomposition theorem from~\cite[\textsection 6]{rostam} to the case of a possibly infinite quiver. A careful analysis of the proofs in this section will be the starting point of several proofs later in the paper.  %

\subsection{Definition}
\label{subsection:definition_klr}

Let $\Gamma$ be a loop-free quiver, possibly infinite.  We write $I$ (respectively $A$) for the vertex (resp. arrow) set. We have a map $A \to I \times I$ given by $A \ni a \mapsto \bigl(o(a), t(a)\bigr) \in I \times I$. The loop-free condition says that for all $a \in A$ we have $o(a) \neq t(a)$.  For any $i, j \in I$, we write $\lvert i \to j\rvert$ for the (finite) number of $a \in A$ such that $o(a) = i$ and $t(a) = j$. We also define $i \cdot j \coloneqq \lvert i \to j \rvert + \lvert i \leftarrow j\rvert$. (We warn the reader that the usual quantity is $- i \cdot j$.) For any $i, j \in I$ we define
\[
d(i, j) \coloneqq \begin{cases}
i \cdot j, &\text{if } i \neq j,
\\
-2, &\text{otherwise}.
\end{cases}
\]
Let $u, v$ be two indeterminates over $K$. For any $i, j \in I$, we define the polynomial $Q_{ij}(u, v) \in K[u, v]$ by
\begin{equation}
\label{equation:def_Q}
Q_{ij}(u, v) \coloneqq \begin{cases}
(-1)^{\lvert i \to j \rvert} (u - v)^{i \cdot j}, &\text{if } i \neq j,
\\
0, &\text{otherwise,}
\end{cases}
\end{equation}
Note that
\begin{equation}
\label{equation:Q_ij(uv)=Q_ji(vu)=Q_ij(-v-u)}
Q_{ij}(u, v) = Q_{ji}(v, u) = Q_{ij}(-v, -u).
\end{equation}

Let $n \in \mathbb{N}$ and $\mathfrak{S}_n$ be the symmetric group on $n$ letters. We denote by $r_a$ the transposition $(a, a+1) \in \mathfrak{S}_n$ for any $a \in \{1, \dots, n-1\}$. We will consider the following two actions of $\mathfrak{S}_n$:
\begin{itemize}
\item the natural action on $\{1, \dots, n\}$, given by $r_a \cdot i \coloneqq r_a(i)$ for all $a \in \{1, \dots, n-1\}$ and $i \in \{1, \dots, n\}$ ;
\item the action on $I^n$ by place permutation, given by
\begin{equation}
\label{equation:action_ra}
r_a \cdot (\dots, i_a, i_{a+1},\dots) \coloneqq (\dots,i_{a+1}, i_a,\dots), 
\end{equation}
for any $\tuple i =(i_1, \dots, i_n) \in I^n$ and $a \in \{1, \dots, n-1\}$.
\end{itemize}
 Let $\alpha \subseteq I^n$ be a finite $\mathfrak{S}_n$-stable subset, that is, a finite union of $\mathfrak{S}_n$-orbits.

\begin{definition}[Khovanov--Lauda \cite{khovanov-lauda_diagrammaticI, khovanov-lauda_diagrammaticII}, Rouquier \cite{rouquier_2}]
The \emph{quiver Hecke algebra} associated with the quiver $\Gamma$ and the finite stable $\mathfrak{S}_n$-subset $\alpha \subseteq I^n$, denoted by $\H_\alpha(\Gamma)$, is the associative unitary $K$-algebra generated by elements
\[
\{y_a\}_{1 \leq a \leq n} \cup
\{\psi_b\}_{1 \leq b\leq n-1}  \cup \{e(\tuple{i})\}_{\tuple{i} \in \alpha},
\]
and relations, for any $\tuple i, \tuple j \in \alpha$ and $a, b \in \{1, \dots, n\}$,
\begin{align}
\label{relation:idempotents}
\sum_{\tuple{i} \in \alpha} e(\tuple{i}) &= 1,
&
e(\tuple{i}) e(\tuple{j}) &= \delta_{\tuple{i}\tuple{j}} e(\tuple{i}),
&
y_a y_b &= y_b y_a,
&
y_a e(\tuple{i}) &= e(\tuple{i})y_a,
\end{align}
and
\begin{align}
\psi_a e(\tuple{i}) &= e(r_a \cdot \tuple{i}) \psi_a,
\label{relation:psia_e(i)}
\\
\label{relation:psib_yj}
(\psi_a y_b - y_{r_a(b)} \psi_a)e(\tuple{i}) &= \begin{cases}
-e(\tuple{i}), & \text{if } b = a \text{ and } i_a = i_{a+1},
\\
e(\tuple{i}), &\text{if } b = a+1 \text{ and } i_a = i_{a+1},
\\
0, &\text{otherwise,}
\end{cases}
\end{align}
if $a \leq n - 1$,  and finally
\begin{align}
\psi_a \psi_b &= \psi_b \psi_a, \text{ if } \lvert a - b \rvert > 1,
\label{relation:psia_psib}
\\
\psi_a^2 e(\tuple{i}) &= Q_{i_a i_{a+1}}(y_a, y_{a+1})e(\tuple{i}),
\label{relation:psia2}
\\
(\psi_{b+1} \psi_b \psi_{b+1} - \psi_b \psi_{b+1} \psi_b )e(\tuple{i}) &=
\begin{dcases}
\frac{Q_{i_b i_{b+1}}(y_b,y_{b+1}) - Q_{i_b i_{b+1}}(y_{b+2}, y_{b+1})}{y_b - y_{b+2}}e(\tuple{i}), &\text{if } i_b = i_{b+2},
\\
0, &\text{otherwise},
\end{dcases}
\label{relation:psi_tresse3}
\end{align}
if $a \leq n - 1$ and $b \leq  n-2$.
\end{definition}

We may form the direct sum $\H_n(\Gamma) \coloneqq \bigoplus_\alpha \H_\alpha(\Gamma)$, where $\alpha$ runs over all the orbits of $I^n$ under the action of $\mathfrak{S}_n$.  If $\Gamma$ is finite, the direct sum is finite and $\H_n(\Gamma)$ is a unitary algebra, with unit $\sum_{\tuple i \in I^n} e(\tuple i)$.  Note that if $n = 0$ then $\H_\alpha(\Gamma) = \H_0(\Gamma) = K$.

\begin{proposition}[\cite{khovanov-lauda_diagrammaticI,khovanov-lauda_diagrammaticII,rouquier_2}]
\label{proposition:grading_klr}
The algebra $\H_\alpha(\Gamma)$ is endowed with the $\mathbb{Z}$-grading given by
\begin{align*}
\deg e(\tuple i) &=0,
\\
\deg y_a &= 2,
\\
\deg \psi_b e(\tuple i) &= d(i_b, i_{b+1}),
\end{align*}
for all $\tuple i \in \alpha$ and $a, b \in \{1, \dots, n\}$ with $b \leq n-1$.
\end{proposition}

For any $w \in \mathfrak{S}_n$, choose a reduced expression $w = r_{a_1} \cdots r_{a_k}$ and define $\psi_w \coloneqq \psi_{a_1} \cdots \psi_{a_k}$. Note that the element $\psi_w$ may depend on the chosen reduced expression.

\begin{proposition}[\cite{khovanov-lauda_diagrammaticI,khovanov-lauda_diagrammaticII,rouquier_2}]
\label{proposition:base_klr}
The algebra $\H_\alpha(\Gamma)$ is a free $K$-module, and
\[
\left\{y_1^{a_1} \cdots y_n^{a_n} \psi_w e(\tuple i) : a_i \in \mathbb{N}, w \in \mathfrak{S}_n, \tuple i \in \alpha\right\},
\]
is a $K$-basis.
\end{proposition}

\begin{remark}\label{rem-maps-orbits}
We recall that there is a one-to-one correspondence between $\mathfrak{S}_n$-orbits $\alpha\subset I^n$ and maps $\hat{\alpha}\ :\ I\to\mathbb{N}$ of \emph{weight} $n$, namely such that $\sum_{i\in I}\hat{\alpha}(i)=n$ (the number $\hat{\alpha}(i)$ counts the number of occurrence of $i$ in any element in the orbit $\alpha$).
\end{remark}

\subsection{Disjoint union of quivers}\label{subsection:disjoint-quiver}

Let $d \in \mathbb{N}^*$. Like in~\cite[\textsection 6.1.3]{rostam}, we assume that the quiver $\Gamma$ decomposes as a disjoint union of full subquivers 
$$\Gamma = \bigsqcup_{j =1}^d \Gamma^{(j)}\ ,$$
where there are no arrows between $\Gamma^{(j)}$ and $\Gamma^{(j')}$ if $j \neq j'$. We denote by $I = \amalg_{j =1}^d I^{(j)}$ the subsequent partition of the vertex set. Note that $Q_{ii'} = 1$ whenever $i\in I^{(j)}$ and $i' \in I^{(j')}$ with $j \neq j'$. 

Now we consider a special class of finite unions of $\mathfrak{S}_n$-orbits in $I^n$. We let $G$ be a finite group acting on $I$ and, for each $j \in \{1, \dots, d\}$, we assume that  $I^{(j)}$  is stable under the action of $G$. We denote
\[G_n=G^n \rtimes \mathfrak{S}_n\ ,\]
the semi-direct product where $\mathfrak{S}_n$ acts on place permutation on $G^n$.

The semidirect product $G_n$ acts naturally on $I^n$. For any $g = (g_1, \dots, g_n) \in G^n$ and $w \in \mathfrak{S}_n$ we have, for all $(i_1, \dots, i_n) \in I^n$,
\[
(g, w) \cdot (i_1, \dots, i_n) = \bigl(g_1 \cdot i_{w^{-1}(1)}, \dots, g_n \cdot i_{w^{-1}(n)}\bigr)\ .
\]

We fix $\alpha \subseteq I^n$ to be a $G_n$-orbit. Note that $\alpha$ is indeed a finite $\mathfrak{S}_n$-stable subset of $I^n$ as in~\textsection\ref{subsection:definition_klr}. 

\subsubsection{Decomposition of orbits}\label{subsubsection-orbits}
For any $\tuple i \in \alpha$ and $j \in \{1, \dots, d\}$, let $\tuple i^{(j)}$ be the tuple obtained from $\tuple i$ by removing the entries that are not in~$I^{(j)}$. We denote by $\sizepart j {\tuple i}$ the number of remaining entries, that is, the number of components of~$\tuple i^{(j)}$. It follows easily from the fact that each $I^{(j)}$ is stable under the action of $G$ that:
\begin{equation}
\label{equation:assumption_alpha}
\text{the tuple }\ (\sizepart 1 {\tuple i}, \dots, \sizepart d {\tuple i})\ \text{ is the same for each $\tuple i\in\alpha$.}
\end{equation}
Thus, we denote, for each $j \in \{1, \dots, d\}$, by $\sizepart j \alpha$ the unique value of $\sizepart j {\tuple i}$ for  $\tuple i\in\alpha$.  We may simply write $\sizeglob j$ instead of $\sizepart j \alpha$ when $\alpha$ is clear from the context. Note that $\sizeglob 1+\dots+\sizeglob d=n$.

We define
\[
\alpha^{(j)} \coloneqq \left\{ \tuple i^{(j)} : \tuple i \in \alpha\right\} \subseteq (I^{(j)})^{\sizeglob j}.
\]
The set $\alpha^{(j)}$ is a finite $\mathfrak{S}_{\sizeglob j}$-stable subset of $(I^{(j)})^{\sizeglob j}$. We will see in~\eqref{proposition:decomposition_orbits} that it is in fact a $G_{\sizeglob j}$-orbit.

In addition to \eqref{equation:assumption_alpha}, we will need the following property of $\alpha$.
\begin{proposition}\label{prop-assumption2}
Recall that $\alpha \subseteq I^n$ is a $G_n$-orbit.
We have:
\begin{equation}
\label{equation:assumption_alpha2}
\alpha^{(1)}\times \dots \times \alpha^{(d)}\subset \alpha\ .
\end{equation}
where we use implicitly the natural identification (by concatenation) of $I^{\sizeglob 1}\times \dots \times I^{\sizeglob d}$ with a subset of $I^n$.
\end{proposition}
\begin{proof}
Let us provide a proof which shows all the various elements explicitly. Since $\alpha$ is a $G_n$-orbit, it can be written of the form:
\[\alpha=\{\bigl(g_1 \cdot i_{w^{-1}(1)}, \dots, g_n \cdot i_{w^{-1}(n)}\bigr)\ |\ g_1,\dots,g_n\in G\,,\ w\in\mathfrak{S}_n\}\,,\]
for some element $(i_1,\dots,i_n)\in I^n$. By invariance under $\mathfrak{S}_n$, we can choose $(i_1,\dots,i_n)$ in an ordered form as follows:
\[\bigl(i^1_1,\dots,i^1_{\sizeglob 1},\dots\dots,i^d_1,\dots,i^d_{\sizeglob d}\bigr),\]
where $i^j_k\in I^{(j)}$ for all $j \in \{1, \dots, d\}$ and $k \in \{1, \dots, \sizeglob j\}$.
Then it is clear that for each $j \in \{1,\dots,d\}$, we have simply
\[\alpha^{(j)}=\left\{\bigl(g_1 \cdot i^j_{w^{-1}(1)}, \dots, g_{\sizeglob j} \cdot i^j_{w^{-1}(\sizeglob j)}\bigr)\left\vert\, g_1,\dots,g_{\sizeglob j}\in G\,,\ w\in\mathfrak{S}_{\sizeglob j}\right.\right\}\,.\]
Property~\eqref{equation:assumption_alpha2} is now immediate to check. 
\end{proof}

From the proof of the preceding proposition, it is easy to see that the map 
\begin{equation}
\label{proposition:decomposition_orbits}
\bigl\{G_n\text{-orbits of } I^n \bigr\} \longrightarrow \bigsqcup_{\substack{\sizeglob 1, \dots, \sizeglob d \geq 0 \\ \sizeglob 1 + \dots + \sizeglob d = n}} \prod_{j = 1}^d\bigl\{ G_{\sizeglob j}\text{-orbits of } \bigl(I^{(j)}\bigr)^{\sizeglob j}\bigr\},
\end{equation}
given by $\alpha \mapsto \bigl(\alpha^{(1)},\dots,\alpha^{(d)}\bigr)$ is a bijection. The inverse map associates to $\bigl(\alpha^{(1)},\dots,\alpha^{(d)}\bigr)$ the smallest $G_n$-stable subset in $I^n$ containing $\alpha^{(1)}\times\dots\times\alpha^{(d)}$.

\begin{remark}
What we actually need for the results of this section is a subset $\alpha$ satisfying properties \eqref{equation:assumption_alpha} and \eqref{equation:assumption_alpha2}. However, since we will use in all the paper only $G_n$-orbits, we find it more convenient to start directly with $G_n$-orbits. In fact we will only use the groups $G=\{1\}$ and $G=\mathbb{Z}/2\mathbb{Z}$, but considering an arbitrary finite group $G$ does not lead to any complication.
\end{remark}

\begin{remark}\label{rem-orbits}
$\bullet$ Let $\Omega$ be the set of $G$-orbits of $I$. Generalising Remark~\ref{rem-maps-orbits}, it is easy to see that there is a one-to-one correspondence between $G_n$-orbits $\alpha \subseteq I^n$ and maps $\hat\alpha : \Omega \to \mathbb{N}$ such that $\sum_{\omega \in \Omega} \hat\alpha(\omega) = n$. If $\alpha \subseteq I^n$ is a $G_n$-orbit and $\omega \in \Omega$, then $\hat\alpha(\omega)$ counts the number of occurrence of the elements of $\omega$ in any element of $\alpha$.

$\bullet$ For each $j=1,\dots,d$, let $\Omega^{(j)}$ be the set of $G$-orbits of $I^{(j)}$. We have $\Omega=\amalg_{j=1}^d\Omega^{(j)}$. Then the bijection \eqref{proposition:decomposition_orbits} in terms of maps simply associates to $\hat\alpha : \Omega \to \mathbb{N}$ the restrictions $\hat\alpha\rvert_{\Omega^{(j)}} : \Omega^{(j)} \to \mathbb{N}$ to each $\Omega^{(j)}$.
\end{remark}

\begin{example}
Let us give an example of a subset $\alpha$ not satisfying property \eqref{equation:assumption_alpha2}. Let $n=2$ and $\alpha=\{(a,A), (A,a), (b,B), (B,b)\}$ where $a,b\in I^{(1)}$ and $A,B\in I^{(2)}$. Then $\alpha$ is a union of two $\mathfrak{S}_2$-orbits and it satisfies \eqref{equation:assumption_alpha}. It does not satisfy \eqref{equation:assumption_alpha2}. Indeed, we have $\alpha^{(1)}=\{a,b\}$ and $\alpha^{(2)}=\{A,B\}$ but, for example, $(a,B)\notin\alpha$.
\end{example}

\subsubsection{Decomposition along the connected components of the quiver}
\label{subsection:decomposition_KLR}

We keep $\alpha \subseteq I^n$ a $G_n$-orbit for some finite group $G$ acting on each set $I^{(j)}$. We may (and we will) simply write~$\sizeglob j$ instead of $\sizepart j \alpha$.

For each $i\in I$, we set $p(i) = j\in\{1,\dots,d\}$ if $i\in I^{(j)}$. Then for each $\tuple i = (i_1, \dots, i_n) \in I^n$, we define its  \emph{profile} by $p(\tuple i)=\bigl( p(i_1),\dots,p(i_n)\bigr)\in\{1,\dots,d\}^n$.
Let 
$$\prof^{\alpha} \coloneqq \{p(\tuple i)\,,\ \tuple i\in\alpha\} \subseteq\{1,\dots,d\}^n$$
be the set of all profiles of elements of $\alpha$. Note that \eqref{equation:assumption_alpha} ensures that $\prof^{\alpha}$ is also a single orbit, now for the action of $\mathfrak{S}_n$ on $\{1,\dots,d\}^n$ by place permutation.

A natural element to consider in this orbit $\prof^{\alpha}$ is
\[
\mathfrak{t}^\alpha \coloneqq (1, \dots, 1,2,\dots,2,\dots \dots, d, \dots, d),
\]
where each $j \in \{1,\dots,d\}$ appears exactly $\sizeglob j$ times. Then every element $\mathfrak{t}\in \prof^{\alpha}$ can be reordered to obtain the distinguished element $\mathfrak{t}^\alpha$.
More precisely, for any $\mathfrak{t}\in \prof^{\alpha}$, the set of elements $w\in\mathfrak{S}_n$ such that $w\cdot\mathfrak{t}=\mathfrak{t}^\alpha$ forms a right coset in $\mathfrak{S}_n$ for the subgroup $\mathfrak{S}_{\sizeglob 1}\times\dots\times \mathfrak{S}_{\sizeglob d}$ (the stabiliser of $\mathfrak{t}^\alpha$). There is a unique minimal length element in this coset (see \emph{e.g.}~\cite{geck-pfeiffer_characters}) and we denote it $\pi_{\mathfrak{t}}$. In particular,  the element $\pi_{\mathfrak{t}}$ is the unique minimal length element of $\mathfrak{S}_n$ such that $\pi_{\mathfrak{t}}\cdot\mathfrak{t}=\mathfrak{t}^\alpha$.

%Recall that $\alpha \subseteq I^n$ is a finite subset stable under the action of $\mathfrak{S}_n$ satisfying~\eqref{equation:assumption_alpha}.
For any $\mathfrak{t}\in \prof^{\alpha}$, we define the idempotent
\[e(\mathfrak{t})=\sum_{\substack{\tuple{i} \in \alpha \\ p(\tuple i) = \mathfrak{t}}} e(\tuple{i}) \in \H_\alpha(\Gamma),\]
and we set
\[
\I \coloneqq \bigl\{e(\mathfrak{t}) : \mathfrak{t} \in \prof^\alpha\bigr\}.
\]
It is a complete set of orthogonal idempotent and its cardinality is $\binom{n}{ \sizeglob 1, \dots, \sizeglob d}$. 
Then, for any $\mathfrak{t}\in \prof^{\alpha}$ we fix a reduced expression $\pi_\mathfrak{t} = r_{a_1} \cdots r_{a_k}$  and define
\begin{subequations}
\label{subequations:def_psit_phit}
\begin{align}
\psi_\mathfrak{t} &\coloneqq \psi_{a_1} \cdots \psi_{a_k} \in \H_\alpha(\Gamma),
\\
\phi_\mathfrak{t} &\coloneqq \psi_{a_k} \cdots \psi_{a_1} \in \H_\alpha(\Gamma).
\end{align}
\end{subequations}

In the following proposition, the grading on $\mathrm{Mat}_{\I}\bigl(e(\mathfrak{t}^\alpha) \H_\alpha(\Gamma)e(\mathfrak{t}^\alpha)\bigr)$ is trivially induced from the grading on $e(\mathfrak{t}^\alpha) \H_\alpha(\Gamma)e(\mathfrak{t}^\alpha)$ (an homogeneous element of degree $N$ is a matrix where all coefficients are homogeneous elements of degree $N$).

\begin{proposition}
\label{proposition:isomorphism_KLR_mat_idempotents}
We have an isomorphism of graded algebras:
\[
\H_\alpha(\Gamma) \simeq \mathrm{Mat}_{\binom{n}{n_1, \dots, n_d}}\bigl(e(\mathfrak{t}^\alpha) \H_\alpha(\Gamma)e(\mathfrak{t}^\alpha)\bigr).
\]
\end{proposition}

\begin{proof}
The proof follows the same steps as in \cite{rostam} and we only give a sketch and the precise references to \cite{rostam}.
First we have that the data $\{e(\mathfrak{t}),\psi_\mathfrak{t},\phi_\mathfrak{t}\}_{\mathfrak{t}\in\prof^{\alpha}}$ in $\H_\alpha(\Gamma)$ enters the general setting~\eqref{subequations:phie_psie_e} of Section~\ref{section:general_decomposition_theorem}, namely we have, for any $\mathfrak{t}\in\prof^{\alpha}$ (see~\cite[Proposition 6.18]{rostam}),
\begin{equation}
\label{equation:phit_psit_klr}
\phi_\mathfrak{t} \psi_\mathfrak{t} e(\mathfrak{t})
=
e(\mathfrak{t}) \phi_\mathfrak{t} \psi_\mathfrak{t} = e(\mathfrak{t}).
\end{equation}
The main point to prove~\eqref{equation:phit_psit_klr} is the following fact:
\begin{equation}
\label{equation:psi2_e(s)}
\psi_a^2 e(\mathfrak{t}) = e(\mathfrak{t}),
\end{equation}
for any $a \in \{1, \dots, n-1\}$ and $\mathfrak{t} \in \prof^\alpha$ such that $\mathfrak{t}_a \neq \mathfrak{t}_{a+1}$ (see~\cite[Lemma 6.15]{rostam}).
Similarly, we obtain, for any $\mathfrak{t}\in\prof^{\alpha}$,
\begin{equation}
\label{equation:psit_et_phit}
\psi_\mathfrak{t} e(\mathfrak{t}) \phi_\mathfrak{t} = \psi_\mathfrak{t} \phi_\mathfrak{t} e(\mathfrak{t}^\alpha) = e(\mathfrak{t}^\alpha).
\end{equation}
This last equality ensures that the set $\J$ in the notation of \textsection\ref{section:general_decomposition_theorem} is $\J=\{e(\mathfrak{t}^\alpha)\}$.
Since $\J$ is reduced to one element, we deduce that the assumption~\eqref{equation:assumption_strong_disjoint} is automatically satisfied, and we can use Corollary~\ref{corollary:isom_A} to obtain the proposition.
Finally, the fact that the isomorphism is homogeneous follows from $\deg \psi_\mathfrak{t} e(\mathfrak{t})= \deg \phi_\mathfrak{t} e(\mathfrak{t}) = 0$ for any $\mathfrak{t} \in \prof^\alpha$ (see~\cite[Remark 6.29]{rostam}).
\end{proof}

\begin{remark}
Similarly to~\eqref{equation:psi2_e(s)}, we have (see~\cite[Lemma 6.20]{rostam})
\begin{equation}
\label{equation:ya_psit}
y_a \phi_\mathfrak{t} e(\mathfrak{t}) = \phi_\mathfrak{t} y_{\pi_\mathfrak{t}(a)} e(\mathfrak{t}),
\end{equation}
for any $a \in \{1, \dots, n-1\}$  and $\mathfrak{t} \in \prof^\alpha$ such that $\mathfrak{t}_a \neq \mathfrak{t}_{a+1}$, and also (see~\cite[Lemma 6.15]{rostam})
\begin{equation}
\label{equation:psi_tresse_exact}
\psi_{a+1} \psi_a \psi_{a+1} e(\mathfrak{t}) = \psi_a \psi_{a+1} \psi_a e(\mathfrak{t}),
\end{equation}
for any $a \in \{1, \dots, n-2\}$  and $\mathfrak{t} \in \prof^\alpha$ such that $\mathfrak{t}_a \neq \mathfrak{t}_{a+2}$. In particular, \eqref{equation:psi_tresse_exact} implies that the quantities $\psi_\mathfrak{t} e(\mathfrak{t})$ and $e(\mathfrak{t})\phi_\mathfrak{t}$ do not depend on the chosen reduced expression for $\pi_\mathfrak{t}$.
\end{remark}

\subsubsection{Expression as a tensor product}

We now want to write the algebra $e(\mathfrak{t}^\alpha) \H_\alpha(\Gamma)e(\mathfrak{t}^\alpha)$ as a tensor product. Recall that $\alpha$ is a $G_n$-orbit and thus satisfies properties~\eqref{equation:assumption_alpha} and~\eqref{equation:assumption_alpha2}. We have already used the first property. The second will be explicitly used during the proof of the next result.

Note that, for any $j \in \{1, \dots, d\}$, the algebra $\H_{\alpha^{(j)}}(\Gamma^{(j)})$ is well-defined since $\alpha^{(j)}$ consists of $\sizeglob j$-tuples of vertices $I^{(j)}$ of $\Gamma^{(j)}$ and is stable under permutations (see \textsection\ref{subsubsection-orbits}). 
\begin{theorem}
\label{theorem:klr_isom_idempot_tensor}
We have an (explicit) isomorphism of graded algebras:
\[
e(\mathfrak{t}^\alpha)\H_{\alpha}(\Gamma)e(\mathfrak{t}^\alpha) \simeq \H_{\alpha^{(1)}}(\Gamma^{(1)}) \otimes \dots \otimes \H_{\alpha^{(d)}}(\Gamma^{(d)}).
\]
\end{theorem}

\begin{proof}
We construct an algebra homomorphism $f$ from the tensor product to $e(\mathfrak{t}^\alpha)\H_{\alpha}(\Gamma)e(\mathfrak{t}^\alpha)$ as follows. For any $\tuple i^{(j)} \in \alpha^{(j)} \subseteq (I^{(j)})^{\sizeglob j}$ with $j \in \{1, \dots, d\}$ we define
\[
f\bigl(e(\tuple i^{(1)}) \otimes \dots \otimes e(\tuple i^{(d)})\bigr) \coloneqq e(\tuple i^{(1)}, \dots, \tuple i^{(d)}).
\]
Note that $\bigl(\tuple i^{(1)},\dots,\tuple i^{(d)}\bigr) \in \alpha$ due to Proposition \ref{prop-assumption2}. Moreover, for any $j \in \{1, \dots, d\}$ we 
denote $y_a^{(j)}$ and $\psi_b^{(j)}$ the generators of $\H_{\alpha^{(j)}}(\Gamma^{(j)})$ in the tensor product and we
define
\begin{align*}
f(y^{(j)}_a) &\coloneqq e(\mathfrak{t}^\alpha)y_{\sizeglob 1 + \dots + \sizeglob{j-1} + a}e(\mathfrak{t}^\alpha),
\\
f(\psi^{(j)}_b) &\coloneqq e(\mathfrak{t}^\alpha)\psi_{\sizeglob 1 + \dots + \sizeglob{j-1} + b}e(\mathfrak{t}^\alpha),
\end{align*}
for all $a, b \in \{1, \dots, \sizeglob j\}$ with $b \leq \sizeglob j - 1$.  By~\cite[Lemma 6.24]{rostam}, the map $f$ is indeed a homomorphism. Using the basis of Proposition~\ref{proposition:base_klr}, we can prove that $f$ sends a basis onto a basis and thus is an isomorphism (see~\cite[Proposition 6.25]{rostam}). Finally, the isomorphism $f$ is clearly homogeneous.
\end{proof}

Combining Theorem~\ref{theorem:klr_isom_idempot_tensor} with Proposition~\ref{proposition:isomorphism_KLR_mat_idempotents}, we obtain the main result of this section.

\begin{corollary}
\label{corollary:isom_dec_klr_block}
We have an (explicit) isomorphism of graded algebras:
\[
\H_{\alpha}(\Gamma) \simeq \mathrm{Mat}_{\binom{n}{n_1, \dots, n_d}}\left(\bigotimes_{j  = 1}^d \H_{\alpha^{(j)}}(\Gamma^{(j)})\right).
\]
\end{corollary}

\begin{remark}
If $\alpha=\amalg_{i=1}^k \alpha_i$ the decomposition of $\alpha$ into $\mathfrak{S}_n$-orbits, then we have $\H_{\alpha}(\Gamma)=\oplus_{i=1}^k \H_{\alpha_i}(\Gamma)$. So of course, as far as the algebras $\H_{\alpha}(\Gamma)$ are concerned, taking $\alpha$ a single $\mathfrak{S}_n$-orbit would be enough. However, we really needed a more general setting since we will apply later the results above for orbits $\alpha\subset I^n$ of the Weyl group of type B.
\end{remark}

%Contrary to the results of~\cite[\textsection 6]{rostam}, the isomorphism of Corollary~\ref{corollary:isom_dec_klr_block} is valid even if the quiver $\Gamma$ is infinite.
We now show how to recover~\cite[Theorem 6.26]{rostam}, with the difference that the result we obtain here is also valid if the quiver $\Gamma$ is infinite.

\begin{corollary}
\label{corollary:isom_dec_klr_n}
We have an (explicit) isomorphism of graded algebras:
\[
\H_n(\Gamma) \simeq \bigoplus_{\substack{n_1, \dots, n_d \geq 0 \\ n_1 + \dots + n_d = n}} \mathrm{Mat}_{\binom{n}{n_1, \dots, n_d}} \left(\bigotimes_{j = 1}^d \H_{n_j}(\Gamma^{(j)})\right).
\]
\end{corollary}

\begin{proof}
We write $\orbitsS{n}$ to denote the $\mathfrak{S}_n$-orbits in $I^n$. We apply the isomorphism of Corollary~\ref{corollary:isom_dec_klr_block} in each term of right-hand side of the equality $\H_n(\Gamma) = \oplus_\alpha \H_\alpha(\Gamma)$, where $\alpha$ runs over $\orbitsS{n}$ (so we use the situation $G=\{1\}$ here). Recalling the $1{:}1$-correspondence in~\eqref{proposition:decomposition_orbits}, we obtain
\begin{align*}
\H_n(\Gamma)
&\simeq
\bigoplus_{\alpha \in \orbitsS{n}} \H_\alpha(\Gamma)
\\
&\simeq
\bigoplus_{\alpha \in \orbitsS{n}} \mathrm{Mat}_{\binom{n}{\sizepart 1 \alpha, \dots, \sizepart d \alpha}}\left(\bigotimes_{j = 1}^d \H_{\alpha^{(j)}}(\Gamma^{(j)})\right)
\\
&\simeq
\bigoplus_{\substack{n_1, \dots, n_d \geq 0 \\ n_1 + \dots + n_d = n}}
\;
\bigoplus_{\substack{\alpha \in\orbitsS{n} \\ \sizepart j \alpha = n_j}} \mathrm{Mat}_{\binom{n}{\sizepart 1\alpha, \dots, \sizepart d \alpha}}\left(\bigotimes_{j = 1}^d \H_{\alpha^{(j)}}(\Gamma^{(j)})\right)
\\
&\simeq
\bigoplus_{\substack{n_1, \dots, n_d \geq 0 \\ n_1 + \dots + n_d = n}}
 \mathrm{Mat}_{\binom{n}{n_1, \dots, n_d}}\left(\bigoplus_{\substack{\alpha \in\orbitsS{n} \\ \sizepart j \alpha = n_j}}
 \;
 \bigotimes_{j = 1}^d \H_{\alpha^{(j)}}(\Gamma^{(j)})\right)
 \\
&\simeq
\bigoplus_{\substack{n_1, \dots, n_d \geq 0 \\ n_1 + \dots + n_d = n}}
 \mathrm{Mat}_{\binom{n}{n_1, \dots, n_d}}\left(
 \bigoplus_{\alpha^{(1)} \in \orbitsS{n_1}}
 \dots
 \bigoplus_{\alpha^{(d)} \in \orbitsS{n_d}}
 \;
 \bigotimes_{j = 1}^d \H_{\alpha^{(j)}}(\Gamma^{(j)})\right)
  \\
&\simeq
\bigoplus_{\substack{n_1, \dots, n_d \geq 0 \\ n_1 + \dots + n_d = n}}
 \mathrm{Mat}_{\binom{n}{n_1, \dots, n_d}}\left(
 \bigotimes_{j = 1}^d \H_{n_j}(\Gamma^{(j)})\right),
\end{align*}
as desired.
\end{proof}

\subsubsection{Cyclotomic case}
\label{subsection:cyclotomic_klr}

We keep the above setting with the quiver $\Gamma$, its full subquivers $\Gamma^{(j)}$ and a $G_n$-orbit $\alpha$.
In addition, let $\Lambda = (\Lambda_i)_{i \in I}$ be a finitely-supported family of non-negative integers.

\begin{definition}[\protect\cite{rouquier_2,BK}]
The \emph{cyclotomic} quiver Hecke algebra $\H_\alpha^{\Lambda}(\Gamma)$ is the quotient of the quiver Hecke algebra $\H_\alpha(\Gamma)$ by the two-sided ideal $\ideal_\alpha^{\Lambda}$ generated by the relations
\begin{equation}
\label{equation:ideal_KLR}
y_1^{\Lambda_{i_1}} e(\tuple i) = 0,
\end{equation}
for all $\tuple i = (i_1, \dots, i_n) \in \alpha$.
\end{definition}

Since the above relations are homogeneous, the cyclotomic quiver Hecke algebras is graded, as in Proposition~\ref{proposition:grading_klr}. Note that if $\Lambda_i = 0$ for all $i$ then
\[
\H_\alpha^\Lambda(\Gamma) = \begin{cases}
\{0\}, &\text{if } n \geq 1,
\\
K, &\text{if } n = 0.
\end{cases}
\]

As in~\cite[\textsection 6.4.1]{rostam}, we want to state Corollaries~\ref{corollary:isom_dec_klr_block} and~\ref{corollary:isom_dec_klr_n} in the cyclotomic setting. First, for any $j \in \{1, \dots, d\}$ let $\Lambda^{(j)}$ be the restriction of $\Lambda$ to $I^{(j)}$.

\begin{theorem}
\label{theorem:isom_dec_klr_cyclo_block}
We have an (explicit) isomorphism of graded algebras:
\[
\H_\alpha^{\Lambda}(\Gamma) \simeq \mathrm{Mat}_{\binom{n}{n_1, \dots, n_d}} \left(\bigotimes_{j = 1}^d \H_{\alpha^{(j)}}^{\Lambda^{(j)}}(\Gamma^{(j)})\right)\ .
\]
\end{theorem}

\begin{proof} %modifs dans la preuve
The proof is similar to the one of~\cite[Theorem 6.30]{rostam}. We provide details since it will be used later in the paper.

Note that $\otimes_{j = 1}^d \H_{\alpha^{(j)}}^{\Lambda^{(j)}}(\Gamma^{(j)})$ is the quotient of $\otimes_{j = 1}^d \H_{\alpha^{(j)}}(\Gamma^{(j)})$  by the two-sided ideal
\[
\ideal_{\alpha, \otimes}^{\Lambda} \coloneqq \langle 1\otimes \dots \otimes 1\otimes \ideal_{\alpha^{(j)}}^{\Lambda^{(j)}}\otimes 1\otimes\dots\otimes 1\,,\ j=1,\dots,d\rangle
\]
generated by the ideals $\ideal_{\alpha^{(j)}}^{\Lambda^{(j)}}$ in position $j$ in the tensor product. 
We will identify the algebra $\otimes_{j = 1}^d \H_{\alpha^{(j)}}(\Gamma^{(j)})$ with the algebra $e(\mathfrak{t}^\alpha)\H_\alpha(\Gamma)e(\mathfrak{t}^\alpha)$ thanks to the explicit isomorphism given in the proof of Theorem \ref{theorem:klr_isom_idempot_tensor}. With this identification, the ideal $\ideal_{\alpha, \otimes}^{\Lambda}$ is generated by the elements
\[y_b^{\Lambda_{i_b}}e(\tuple i)\ ,\]
where $\tuple i \in \alpha$ is of profile $\mathfrak{t}^{\alpha}$, and $b$ is of the form $b=n_1+\dots+n_{j-1}+1$ for $j\in\{1,\dots,d\}$.

Now let $\theta$ be the isomorphism of Proposition~\ref{proposition:isomorphism_KLR_mat_idempotents} and $\eta$ its inverse. For convenience we denote during the proof $N:=\binom{n}{n_1, \dots, n_d}$. We will prove the following two inclusions:
\begin{subequations}
\label{subequations:theta_eta_ideals}
\begin{align}
\theta\left(\ideal_\alpha^{\Lambda}\right)
&\subseteq \mathrm{Mat}_{N} \left(\ideal_{\alpha, \otimes}^{\Lambda}\right),
\label{subequation:theta_ideal}
\\
\ideal_\alpha^{\Lambda}
&\supseteq \eta\left(\mathrm{Mat}_{N} \left(\ideal_{\alpha, \otimes}^{\Lambda}\right)\right).
\label{subequation:eta_ideal}
\end{align}
\end{subequations}

 Let $\mathfrak{t}', \mathfrak{t} \in \prof^\alpha$. First recall that for $h \in e(\mathfrak{t}') \H_\alpha(\Gamma)e(\mathfrak{t})$, we have
\[
\theta(h) = \psi_{\mathfrak{t}'} h \phi_\mathfrak{t} M_{\mathfrak{t}'\mathfrak{t}} \in \mathrm{Mat}_{N}\left(e(\mathfrak{t}^\alpha) \H_\alpha(\Gamma) e(\mathfrak{t}^\alpha)\right).
\]
while for $h \in e(\mathfrak{t}^\alpha) \H_\alpha(\Gamma)e(\mathfrak{t}^\alpha)$ we have
\[
\eta\left(h M_{\mathfrak{t}'\mathfrak{t}}\right) = \phi_{\mathfrak{t}'} h \psi_\mathfrak{t}\ ,
\]
where the elements $\phi_\mathfrak{t},\psi_\mathfrak{t}$ were introduced in \eqref{subequations:def_psit_phit}.

$\bullet$ Let $\tuple i \in \alpha$ of profile $\mathfrak{t}$. By~\eqref{equation:phit_psit_klr},~\eqref{relation:psia_e(i)} and~\eqref{equation:ya_psit} we have:
\[
y_1^{\Lambda_{i_1}} e(\tuple i)=y_1^{\Lambda_{i_1}} e(\tuple i)e(\mathfrak{t})=y_1^{\Lambda_{i_1}} e(\tuple i)\phi_\mathfrak{t}\psi_\mathfrak{t}e(\mathfrak{t})
=
\phi_\mathfrak{t} y_{\pi_\mathfrak{t}(1)}^{\Lambda_{i_1}} e(\pi_\mathfrak{t} \cdot \tuple i) \psi_\mathfrak{t}e(\mathfrak{t}).
\]
Thus, to prove~\eqref{subequation:theta_ideal} it suffices to show that
\[
\theta\left(y_{\pi_\mathfrak{t}(1)}^{\Lambda_{i_1}} e(\pi_{\mathfrak{t}}\cdot \tuple i)\right) \in \mathrm{Mat}_{N} \left(\ideal_{\alpha, \otimes}^{\Lambda}\right)\ .
\]
By definition of $\pi_{\mathfrak{t}}$, we have that $\tuple i'\coloneqq \pi_\mathfrak{t} \cdot \tuple i$ has profile $\mathfrak{t}^\alpha$ and therefore $y_{\pi_\mathfrak{t}(1)}^{\Lambda_{i_1}} e(\tuple i')\in e(\mathfrak{t}^{\alpha}) \H_\alpha(\Gamma)e(\mathfrak{t}^{\alpha})$. Let $b \coloneqq \pi_{\mathfrak{t}}(1)$ so that we have $i_1=i'_b$, and moreover, by~\cite[Proposition 6.7]{rostam}, the element $b$ is of the form $\sizeglob 1+ \dots + \sizeglob{j-1} + 1$. We conclude that
\[
\theta\left(y_{\pi_\mathfrak{t}(1)}^{\Lambda_{i_1}} e(\pi_{\mathfrak{t}}\cdot \tuple i)\right)=y_b^{\Lambda_{i'_b}} e(\tuple i') M_{\mathfrak{t}^\alpha \mathfrak{t}^\alpha}\in \mathrm{Mat}_{N} \left(\ideal_{\alpha, \otimes}^{\Lambda}\right)\ .
\]

$\bullet$
Let $\tuple i\in\alpha$ with profile $\mathfrak{t}^\alpha$ and let $b=n_1+\dots+n_{j-1}+1$ with $j\in\{1,\dots,d\}$ such that $n_j \neq 0$. Let us prove that
\[
\eta\left(y_b^{\Lambda_{i_b}} e(\tuple i) M_{\mathfrak{t}'\mathfrak{t}}\right) \in \ideal_\alpha^{\Lambda}.
\]
Since $M_{\mathfrak{t}'\mathfrak{t}}=M_{\mathfrak{t}'\mathfrak{t''}}M_{\mathfrak{t}''\mathfrak{t}}$ for any $\mathfrak{t''}$ it is enough to prove it for a single value of $\mathfrak{t}'$. So without loss of generality, since $n_j \neq 0$ we can assume that $\mathfrak{t}'$ starts with $j$ so that $\pi_{\mathfrak{t'}}(1) = b$. We conclude that
\[\eta\left(y_b^{\Lambda_{i_b}} e(\tuple i) M_{\mathfrak{t}'\mathfrak{t}}\right)=
\phi_{\mathfrak{t}'} y_b^{\Lambda_{i_b}} e(\tuple i) \psi_\mathfrak{t}
= y_1^{\Lambda_{i_b}} e(\pi_{\mathfrak{t'}}^{-1}\cdot\tuple i) \phi_{\mathfrak{t}'} \psi_\mathfrak{t} \in \ideal_\alpha^{\Lambda},
\]
since, if we denote $\tuple i'=\pi_{\mathfrak{t'}}^{-1}\cdot\tuple i$ then we have $i'_1=i_b$.

\medskip
This concludes the proof of~\eqref{subequations:theta_eta_ideals} showing that we have
\[
\theta\left(\ideal_\alpha^{\Lambda}\right)
= \mathrm{Mat}_{N} \left(\ideal_{\alpha, \otimes}^{\Lambda}\right)\ .
\]
Thus we can deduce the isomorphism of Theorem~\ref{theorem:isom_dec_klr_cyclo_block} from Corollary~\ref{corollary:isom_dec_klr_block}.
\end{proof}

\begin{remark}\label{rem-cancellation_beta}%petite remarque sur les annulations
\begin{itemize}
\item We saw that if $\Lambda^{(j)}\equiv 0$ on $\Gamma^{(j)}$ for some $j$ then, if moreover $\alpha^{(j)}\neq\emptyset$ (that is, if $n_j(\alpha)\neq0$), we have $\H_{\alpha^{(j)}}^{\Lambda^{(j)}}(\Gamma^{(j)})=\{0\}$  from the defining relations. So in turn, Theorem~\ref{theorem:isom_dec_klr_cyclo_block} implies that $\H_\alpha^{\Lambda}(\Gamma)=\{0\}$.

\item The conclusion of the preceding item can in fact be seen more directly. Indeed the cyclotomic relations in $\H_\alpha^{\Lambda}(\Gamma)$ imply that $e(\tuple i)=0$ for all $\tuple i \in \alpha$ with $i_1\in\Gamma^{(j)}$. So we have that the idempotent $e(\mathfrak{t})$ is 0 for any profile $\mathfrak{t}$ starting with $j$ (and at least one profile like this exists in $\text{Prof}^\alpha$ when $n_j(\alpha)\neq0$). %Now recall from our construction in the first sections (see Formulas \eqref{equation:phit_psit_klr} and \eqref{equation:psit_et_phit}) that we have elements $\psi_\mathfrak{t},\phi_\mathfrak{t}$ for any profile $\mathfrak{t}$ such that:
Since:
\[\psi_\mathfrak{t} e(\mathfrak{t}) \phi_\mathfrak{t} = e(\mathfrak{t}^\alpha)\ \ \ \ \ \text{and}\ \ \ \ \ \phi_\mathfrak{t} e(\mathfrak{t}^{\alpha}) \psi_\mathfrak{t} = e(\mathfrak{t})\ ,\]
it follows immediately that if $n_j(\alpha)\neq 0$ then all idempotents $e(\mathfrak{t})$ are 0 and in turn all idempotents $e(\tuple i)$, $\tuple i\in\alpha$, are 0, which shows that $\H_\alpha^{\Lambda}(\Gamma)=\{0\}$.
\end{itemize}
\end{remark}

As in Corollary~\ref{corollary:isom_dec_klr_n}, we deduce the following corollary.

\begin{corollary}
\label{corollary:isom_dec_klr_cyclo_n}
We have an (explicit) isomorphism of graded algebras:
\[
\H_n^{\Lambda}(\Gamma) \simeq \bigoplus_{\substack{n_1, \dots, n_d \geq 0 \\ n_1 + \dots + n_d = n}} \mathrm{Mat}_{\binom{n}{n_1, \dots, n_d}} \left(\bigotimes_{j = 1}^d \H_{n_j}^{\Lambda^{(j)}}(\Gamma^{(j)})\right).
\]
\end{corollary}

\begin{remark}
\label{rem-cancellation_n}
It follows from Remark~\ref{rem-cancellation_beta} that we can assume that $\Lambda$ is supported on all components of $\Gamma$, that is $\Lambda^{(j)}\not\equiv 0$ for all $j\in\{1,\dots,d\}$. In other words, we can replace from the beginning $\Gamma$ by $\tilde{\Gamma}$ where we removed the components $\Gamma^{(j)}$ such that $\Lambda^{(j)}\equiv 0$. In particular, we have $\H_n^{\Lambda}(\Gamma)=\H_n^{\Lambda_{\vert\tilde I}}(\tilde{\Gamma})$, where $\tilde I$ denotes the vertex set of $\tilde\Gamma$. We could have done that but it turned out to be not really necessary to state Theorem~\ref{theorem:isom_dec_klr_cyclo_block} or Corollary~\ref{corollary:isom_dec_klr_cyclo_n}. For example, in Corollary~\ref{corollary:isom_dec_klr_cyclo_n}, if $\Lambda^{(j)}\equiv 0$ for some $j$ then all the summands with $n_j\neq 0$ are $\{0\}$ and can thus be removed from the direct sum.
\end{remark}

\section{Interpolating quiver Hecke algebras for type B}
\label{section:interpolating}

The aim of this section is to unite the definitions of quiver Hecke algebras for type B that are introduced in~\cite{varagnolo-vasserot_canonical} by Varagnolo and Vasserot and in~\cite{poulain_dandecy-walker_B,poulain_dandecy-walker_D} by the first author and Walker.

\subsection{Definition}
\label{subsection:definition_QHA_typeB}

Let $\Gamma$ be a quiver as in~\textsection\ref{subsection:definition_klr}. We also adopt the notation of this subection.
Let $\theta$ be an involution of $\Gamma$, that is, the map $\theta$ is an involution on both sets $I$ and $A$ and satisfies
\begin{equation}
\label{equation:theta_involution}
\theta(o(a)) = t(\theta(a)),
\end{equation}
for all $a \in A$.
Note the following consequence: for any $i, j \in I$ we have $\lvert i \to j \rvert = \lvert \theta(j) \to \theta(i)\rvert$ and thus
\begin{equation}
\label{equation:theta_dot_invariant}
i \cdot j = \theta(i) \cdot \theta(j).
\end{equation}
It follows from the definition~\eqref{equation:def_Q} of the polynomials $Q_{ij}$  and from~\eqref{equation:theta_involution} again that
\begin{equation}
\label{equation:Qij=Qthetaithetaj}
Q_{ij}(u, v) = Q_{\theta(j)\theta(i)}(u, v),
\end{equation} 
for any $i, j \in I$.

Let $B_n$ be the group of signed permutations of $\{\pm 1, \dots, \pm n\}$, that is, the group of permutations $\pi$ of $\{\pm 1, \dots, \pm n\}$ satisfying $\pi(-i) = -\pi(i)$ for all $i \in \{1, \dots, n\}$. We have a natural isomorphism $B_n \simeq (\mathbb{Z}/2\mathbb{Z})^n \rtimes \mathfrak{S}_n$. In particular we are in the setting of~\textsection\ref{subsection:disjoint-quiver} with $G = \mathbb{Z}/2\mathbb{Z}$, which acts on $I$ via the canonical surjection $G \twoheadrightarrow \langle\theta\rangle$.
We have a natural inclusion $\mathfrak{S}_n \subseteq B_n$, where $r_a$ is identified with $(a, a+1)(-a, -a-1)$ for all $a \in \{1, \dots, n-1\}$. We see $B_n$ as a Weyl group of type B by adding the generator $r_0 \coloneqq (-1, 1)$. The action of $B_n$ on $I^n$ is given by~\eqref{equation:action_ra} and
\[
r_0 \cdot (i_1, \dots, i_n) \coloneqq (\theta(i_1), i_2, \dots, i_n),
\]
for any $\tuple{i} = (i_1, \dots, i_n) \in I^n$. Let $\beta \subseteq I^n$ be a $B_n$-orbit. In particular, the set $\beta$ is a finite $\mathfrak{S}_n$-stable subset of $I^n$.

\begin{remark}
The result of %Lemma~\ref{lemma:correspondance_Gorbits_maps} 
Remark~\ref{rem-orbits} can here be written as follows. There is a one-to-one correspondence between $B_n$-orbits $\beta\subset I^n$ and maps $\hat{\beta}\ :\ I\to\mathbb{N}$ such that $\hat\beta = \hat\beta \circ \theta$ and $\frac{1}{2} \sum_{\substack{i \in I \\ \theta(i) \neq i}} \hat\beta(i) + \sum_{\substack{i \in I \\\theta(i) = i}} \hat\beta(i) = n$ (the number $\hat{\beta}(i)$ counts the number of occurrence of both $i$ and $\theta(i)$  in any element in the orbit $\beta$). See also \cite[Remark 2.5]{poulain_dandecy-walker_B}.
\end{remark}

Let $\lambda \in \mathbb{N}^I$ and  $\gamma \in K^I$. Define
\[
d(i) \coloneqq \begin{cases}
\lambda_i + \lambda_{\theta(i)}, &\text{if } \gamma_i = 0,
\\
-2, &\text{otherwise}.
\end{cases}
\]
For any $i \in I$, we make the following assumptions:
\begin{subequations}
\label{subequations:conditions_gamma}
\begin{align}
\label{equation:condition_gamma}
\theta(i) \neq i &\implies \gamma_i = 0,
\\
\label{equation:condition_strong}
\gamma_i = 0 &\implies [\theta(i) \neq i \text{ or } d(i) = 0].
\end{align}
\end{subequations}

% Let $\lambda \comp_I n$ be a finitely-supported sequence of integers $\lambda \in \mathbb{N}^{(I)}$. For any $i \in I$ we define
%\[
%d^{\lambda, \theta}(i) \coloneqq \begin{cases}
%\lambda_i + \lambda_{\theta(i)}, &\text{if } \theta(i) \neq i,
%\\
%\lambda_i, &\text{if } \theta(i) = i \text{ and } \gamma_i = 0,
%\\
%-2, &\text{otherwise}.
%\end{cases}
%\]
Note that $\gamma$ is $\theta$-invariant, that is, we have
\begin{equation}
\label{equation:gamma_thetainv}
\gamma_{\theta(i)} = \gamma_i, \text{ for all } i \in I.
\end{equation}

\begin{remark}
\label{remark:condition_strong}
$\bullet$ Condition~\eqref{equation:condition_strong} may seem strong; without it we encounter in~\textsection\ref{subsection:gammai1=0=gammai2} useless complications for our means (see also Remark~\ref{remark:not_strong_assumption}). 

$\bullet$ Similarly, one could consider a more general definition than the one below. As for example in \cite[\textsection 3.2]{rouquier_2}, we could remove any reference to a quiver and start only with a family of polynomials associated to the set $I$ with involution $\theta$ (namely, $Q_{ij}[u,v]$ and a polynomial replacing $(-1)^{\lambda_{\theta(i_1)}} y_1^{d(i_1)}$ in the definition below). Then one should look for conditions ensuring the existence of a polynomial representation. We do not pursue in this direction to avoid adding another layer of technicalities.
\end{remark}

\begin{definition}
The algebra $V_\beta(\Gamma, \lambda, \gamma)$ is the unitary associative $K$-algebra generated by elements
\[
\{y_a\}_{1 \leq a \leq n} \cup 
\{\psi_b\}_{0 \leq b \leq n-1}  \cup \{e(\tuple{i})\}_{\tuple{i} \in \beta},
\]
with the relations \eqref{relation:idempotents}--\eqref{relation:psi_tresse3} of Section~\ref{section:toy_KLR} involving all the generators but $\psi_0$,  together with
\begin{align}
\psi_0 e(\tuple{i}) &= e(r_0 \cdot \tuple{i}) \psi_0,
\label{relation:psi0_e(i)}
\\
\psi_0 \psi_b &= \psi_b \psi_0, \qquad\text{for all } b \in \{2, \dots, n-1\},
\label{relation:psi0_psib}
\\
\label{relation:psi0_y1}
(\psi_0 y_1 + y_1 \psi_0)e(\tuple{i}) &= 2\gamma_{i_1} e(\tuple{i}),
\\
\label{relation:psi0_yj}
\psi_0 y_a &= y_a \psi_0, \qquad\text{for all } a \in \{2, \dots, n\},
\\
\label{relation:psi0square}
\psi_0^2 e(\tuple{i}) &= \begin{cases}
(-1)^{\lambda_{\theta(i_1)}} y_1^{d(i_1)} e(\tuple{i}),
&\text{if } \gamma_{i_1} = 0,
\\
0, &\text{otherwise},
\end{cases}
\\
\label{relation:braid4}
\left((\psi_0 \psi_1)^2 - (\psi_1 \psi_0)^2\right)e(\tuple{i}) &=\\
\nonumber & \hspace{-1.2cm}\begin{cases}
(-1)^{\lambda_{\theta(i_1)}} \frac{(-y_1)^{d(i_1)} - y_2^{d(i_1)}}{y_1 + y_2} \psi_1  e(\tuple{i}),
&
\text{if } \gamma_{i_1} = 0 \text{ and }\theta(i_1) = i_2,
\\[1em]
\gamma_{i_2}  \frac{Q_{i_2 i_1}(y_1, -y_2) - Q_{i_2 i_1}(y_1, y_2)}{y_1 y_2} \left(y_1 \psi_0 - \gamma_{i_1}\right) e(\tuple{i}) &\text{otherwise,}
\end{cases}
\end{align}
for all $\tuple i \in \beta$.
%nouvelle equation l'ancienne est la dessous
%\begin{align}
%\label{relation:braid4}
%\left((\psi_0 \psi_1)^2 - (\psi_1 \psi_0)^2\right)e(\tuple{i}) &=
%\begin{cases}
%(-1)^{\lambda_{\theta(i_1)}} \frac{(-y_1)^{d(i_1)} - y_2^{d(i_1)}}{y_1 + y_2} \psi_1  e(\tuple{i}),
%&
%\text{if } \gamma_{i_1} = 0 \text{ and }\theta(i_1) = i_2,
%\\
%\gamma_{i_2} \frac{Q_{\theta(i_1) i_2}(y_1, -y_2) - Q_{\theta(i_1) i_2}(y_1, y_2)}{y_2}\psi_0 %e(\tuple{i}),
%&
%\text{if } \gamma_{i_1} = 0 \neq \gamma_{i_2},
%\\
%\gamma_{i_2}  \frac{Q_{i_1 i_2}(y_1, -y_2) - Q_{i_1 i_2}(y_1, y_2)}{y_1 y_2} \left(y_1 \psi_0 - \gamma_{i_1}\right) e(\tuple{i}),
%&
%\text{if } \gamma_{i_1} \neq 0 \neq \gamma_{i_2} \text{ and } i_1 \neq i_2,
%\\
%0,
%&\text{otherwise},
%\end{cases}
%\end{align}

\end{definition}

% modif des commentaires ci-dessous 
%Note that, in relation~\eqref{relation:braid4}, when $\gamma_{i_1} = 0$ and $\theta(i_1) = i_2$ then by~\eqref{equation:gamma_thetainv} we have $\gamma_{i_2} = 0$, thus the enumerated cases are indeed disjoint. Moreover, when $\gamma_{i_1} \neq 0 \neq \gamma_{i_2}$ and $i_1 \neq i_2$ then by~\eqref{equation:condition_gamma} we have $\theta(i_1) = i_1$ and $\theta(i_2) = i_2$ thus we can use~\eqref{equation:Q_ij(uv)=Q_ji(vu)=Q_ij(-v-u)} and~\eqref{equation:Qij=Qthetaithetaj} so that
It is clear that the fraction in the first line of the right hand side in \eqref{relation:braid4} is a polynomial in $y_1,y_2$.
Then we note that the second line in the right hand side of \eqref{relation:braid4} is 0 when $\gamma_{i_2} = 0$ or when $i_1=i_2$ (recalling~\eqref{equation:def_Q}), and is a polynomial in $y_1, y_2$ when $\gamma_{i_1} = 0$. So for the second line, if $\gamma_{i_1}\neq 0\neq\gamma_{i_2}$ and $i_1\neq i_2$ then by~\eqref{equation:condition_gamma} we have $\theta(i_1) = i_1$ and $\theta(i_2) = i_2$, and thus we can use~\eqref{equation:Q_ij(uv)=Q_ji(vu)=Q_ij(-v-u)} and~\eqref{equation:Qij=Qthetaithetaj} so that
\begin{align*}
\frac{Q_{i_1 i_2}(u, -v) - Q_{i_1 i_2}(u, v)}{uv}
&=
\frac{Q_{i_1 i_2}(v, -u) - Q_{\theta(i_2) \theta(i_1)}(u, v)}{uv}
\\
&=
\frac{Q_{i_1 i_2}(v, -u) - Q_{i_2 i_1}(u, v)}{uv}
\\
&=
\frac{Q_{i_1 i_2}(v, -u) - Q_{i_1 i_2}(v, u)}{uv},
\end{align*}
is a polynomial. 

Finally, note that when $n = 0$ then $V_\beta(\Gamma,\lambda,\gamma) = K$.

\begin{remark}
\label{remark:quotient}
Since $\beta$ is a finite $\mathfrak{S}_n$-stable subset of $I^n$, we can also consider the algebra $\H_\beta(\Gamma)$ as defined in~\textsection\ref{subsection:definition_klr}. The subalgebra of $V_\beta(\Gamma,\lambda,\gamma)$ generated by all the generators but $\psi_0$ is an obvious quotient of $\H_\beta(\Gamma)$ (see also Corollary~\ref{corollary:klr_subalgebra_vv}).
\end{remark}

\begin{remark}
\label{remark:coherence_VV}
If $\theta$ has no fixed point in $I$ then $V_\beta(\Gamma, \lambda,\gamma)$ is exactly the algebra defined in~\cite{varagnolo-vasserot_canonical}. In this case, by~\eqref{equation:condition_gamma} we necessarily have $\gamma_i = 0$ for any $i$ and~\eqref{equation:condition_strong} is automatically satisfied. In particular, in \eqref{relation:braid4} the second line is always zero in this situation.
\end{remark}

\begin{remark}
\label{remark:coherence_PAW}
Assume that $K$ is field of characteristic different from $2$ and let $p, q \in K^\times$ with $q^2 \neq 1 \neq p^2$. Let $\theta : K^\times \to K^\times$ be the scalar inversion. For any $x \in K^\times$, we define the set $I_x \coloneqq \{ x^\epsilon q^{2l} : \epsilon \in \{\pm 1\}, l \in \mathbb{Z}\}$. Let $x_1, \dots, x_k \in K^\times$ such that the sets $I_{x_a}$ are pairwise disjoint. Let $\Gamma$ be the quiver with vertices $I \coloneqq \amalg_{a = 1}^k I_{x_a}$ and arrows between $v$ and $q^2 v$ for all $v \in I$. Finally let $\lambda$ be the indicator function of $P \coloneqq \{\pm p\} \cap I$ and define $\gamma_i \coloneqq 1$ if $\theta(i) = i$ and $\gamma_i \coloneqq 0$ otherwise (thus~\eqref{subequations:conditions_gamma} is satisfied). Then $V_\beta(\Gamma, \lambda, \gamma)$ is exactly the algebra $V_{\tuple x}^{I^\beta}$ defined in~\cite{poulain_dandecy-walker_B}. This is, together with the next remark, the situation relevant for the applications to affine Hecke algebras, see Section \ref{section:morita}.
\end{remark}

\begin{remark}
\label{remark:coherence_PAW_typeD}
The algebra of~\cite[\textsection 3.1]{poulain_dandecy-walker_D} is obtained  with the same choice of $\Gamma,\theta$ as in the preceding remark, together with $\gamma_i\coloneqq 0$ and $\lambda_i \coloneqq 0$ for all $i$. In particular, Condition~\eqref{equation:condition_strong} is satisfied since $d(i) = 0$ for all $i \in I$. We will come back to this particular situation in Section \ref{section:QHA_D}.
\end{remark}

 The algebra $V_\beta(\Gamma, \lambda,\gamma)$ is endowed with the $\mathbb{Z}$-grading given by
\begin{subequations}
\label{gradingB}
\begin{align}
\deg e(\tuple{i}) &= 0,
\\
\deg y_a &= 2,
\\
\deg \psi_0 e(\tuple{i}) &= 
d(i_1),
\\
\deg \psi_b e(\tuple{i}) &= d(i_b, i_{i+1}).
\end{align}
\end{subequations}

The homogeneity of the defining relations that do not involve $\psi_0$ is as in Section~\ref{section:toy_KLR}, the other ones being a simple calculation. For~\eqref{relation:psi0_y1} note that if $\gamma_{i_1} = 0$ there is nothing to check, and if $\gamma_{i_1} \neq 0$ then by definition we have $d(i_1) = -2$  thus $\deg \psi_0 y_1 e(\tuple i) = \deg y_1 \psi_0 e(\tuple i) = 0$. To check the last relation, let us write  $i_1 i_2$ instead of $\tuple{i}$ and even $a$ instead of $i_a$ and $\bar{a}$ instead of $\theta(i_a)$. We have
\begin{align}
\left((\psi_0 \psi_1)^2 - (\psi_1 \psi_0)^2\right)e(12)
&=
\psi_0 \psi_1 \psi_0 \psi_1 e(12) - \psi_1 \psi_0 \psi_1 \psi_0 e(12)
\notag
\\
\label{equation:example_short}
&=
\psi_0 e(1\bar{2}) \psi_1 e(\bar{2}1) \psi_0 e(21) \psi_1 e(12) - \psi_1 e(\bar{2}\bar{1})\psi_0 e(2\bar{1}) \psi_1 e(\bar{1}2) \psi_0 e(12).
\end{align}
We have:
\begin{align*}
\deg \psi_0 e(1\bar{2}) = \deg \psi_0 e(12) &= d(1),
\\
\deg \psi_0 e(21) = \deg\psi_0 e(2\bar{1}) &= d(2).
\end{align*}
Moreover, by~\eqref{equation:theta_dot_invariant} we have
\begin{gather*}
\deg \psi_1 e(\bar{2}1) = d(\bar{2}, 1) = d(1, \bar{2}) = d(\bar{1}, 2) =  \deg\psi_1 e(\bar{1}2),
\\
\deg \psi_1 e(12)  = d(1,2) = d(2,1) = d(\bar{2} ,\bar{1}) = \deg\psi_1 e(\bar{2}\bar{1}).
\end{gather*}
Thus, the quantity $\left((\psi_0 \psi_1)^2 - (\psi_1 \psi_0)^2\right)e(\tuple{i})$ is homogeneous of degree
\[
d(i_1) + d(i_2) + d\bigl(i_1, i_2\bigr) + d\bigl(i_1, \theta(i_2)\bigr).
\]
A quick calculation now shows that the last relation is homogeneous (note that in the first case we have $\gamma_{i_2} = 0$ by~\eqref{equation:gamma_thetainv}).

\subsection{Basis theorem}
\label{subsection:basis_theorem}

We now want to give an analogue of the basis theorem Proposition~\ref{proposition:base_klr} for quiver Hecke algebras. As in~\cite{khovanov-lauda_diagrammaticI,khovanov-lauda_diagrammaticII,rouquier_2}, we will construct a polynomial realisation of $V_\beta(\Gamma,\lambda,\gamma)$.  Let $(P_{ij}(u, v))_{i, j \in I}$ be a family of polynomials satisfying 
\begin{subequations}
\label{subequations:properties_P}
\begin{align}
\label{equation:P(uv)=P(-v-u)}
P_{ij}(u, v) &= P_{ij}(-v, -u),
\\
\label{equation:Pij=Pthetaithetaj}
P_{ij}(u, v) &= P_{\theta(j)\theta(i)}(u, v),
\end{align}
\end{subequations}
and such that
\begin{equation}
\label{equation:PP=Q}
P_{ij}(u, v) P_{ji}(v, u) = Q_{ij}(u, v).
\end{equation}
Note that %$P_{ij}(u, v) \coloneqq \delta_{ij}(u - v)^{\lvert j \to i\rvert}$% c'etait faux comme ca
$P_{ij}(u, v) \coloneqq (u - v)^{\lvert j \to i\rvert}$ if $i\neq j $ and $P_{ij}(u, v) \coloneqq 0$ if $i=j$ is an example of such a family, by~\eqref{equation:theta_involution}. 
Now let $(\alpha_i(y))_{i \in I}$ be a family of polynomials such that
\begin{align}
\label{equation:equation_alpha}
\alpha_{\theta(i)}(y) \alpha_i(-y)
&= (-1)^{\lambda_{\theta(i)}} y^{d(i)},
& \text{if } \gamma_i = 0,
\\
\label{equation:alpha_0}
\alpha_i(y) &= 0, &\text{otherwise}.
\end{align}
Note that if $\gamma_i = 0$ we can just set $\alpha_i(y) \coloneqq y^{\lambda_{\theta(i)}}$. We now consider the sum of polynomials algebras $K[x, \beta] \coloneqq\oplus_{\tuple{i} \in \beta} K[x_1, \dots, x_n]\mathbf{1}_{\tuple{i}}$, where $\mathbf{1}_{\tuple{i}}$ denotes the unit of the summand corresponding to $\tuple i$, so that
\begin{align*}
f \mathbf{1}_{\tuple i} &= \mathbf{1}_{\tuple i} f, &&\text{for all } f \in K[x_1, \dots, x_n] \text{ and }\tuple i \in \beta,
\\
\mathbf{1}_{\tuple i} \mathbf{1}_{\tuple j} &= \delta_{\tuple{ij}} \mathbf{1}_{\tuple i}, &&\text{for all }  \tuple i, \tuple j \in \beta.
\end{align*}
The Weyl group $B_n$ acts on $K[x_1, \dots,x_n]$ by $\prescript{w}{}{f}(x_1, \dots, x_n) \coloneqq f\bigl(w^{-1} \cdot (x_1, \dots, x_n)\bigr)$ for any $w \in B_n$ and $f \in K[x_1, \dots, x_n]$,
%on n'avait pas donne l'action
 where the action of the generator $r_0$ on $(x_1,\dots,x_n)$ is by multiplying $x_1$ by $-1$, and the action of the generator $r_a$, $a=1,\dots,n-1$, on $(x_1,\dots,x_n)$ is by exchanging $x_a$ and $x_{a+1}$. 
The action of $B_n$ on $K[x_1, \dots, x_n]$ extends by linearity to $K[x, \beta]$ by setting $w \star f \mathbf{1}_{\tuple i} \coloneqq \actpol{w}{f} \mathbf{1}_{w \cdot \tuple i}$ for any $\tuple i \in \beta$.

We now consider the linear action of $V_\beta(\Gamma, \lambda,\gamma)$ on  $K[x, \beta]$ given on the generators by
\begin{align*}
e(\tuple{j})\cdot f\mathbf{1}_{\tuple{i}} &\coloneqq \delta_{\tuple{ij}}f\mathbf{1}_{\tuple{i}} = \delta_{\tuple{ij}} \mathbf{1}_{\tuple i}f,
\\
y_a \cdot f\mathbf{1}_{\tuple{i}} &\coloneqq   x_a f\mathbf{1}_{\tuple{i}} = x_a \mathbf{1}_{\tuple i}f,
\\
\psi_b \cdot f\mathbf{1}_{\tuple{i}} &\coloneqq \delta_{i_b, i_{b+1}} \frac{\actpol{r_b}{f} - f}{x_b - x_{b+1}} \mathbf{1}_{\tuple{i}}
+ 
P_{i_b, i_{b+1}}(x_{b+1}, x_b) \actpol{r_b}{f}\mathbf{1}_{r_b \cdot \tuple{i}}
\\
&=	\Bigl(\delta_{i_b, i_{b+1}} (x_b - x_{b+1})^{-1}(r_b - 1)
+ 
P_{i_b, i_{b+1}}(x_{b+1}, x_b) r_b\Bigr) \star f\mathbf{1}_{\tuple i},
\\
\psi_0 \cdot f\mathbf{1}_{\tuple{i}} &\coloneqq \left( \gamma_{i_1}\frac{f - \actpol{r_0}{f}}{x_1} + \alpha_{i_1}(x_1)\actpol{r_0}{f}\right) \mathbf{1}_{r_0 \cdot \tuple{i}}
\\
&=
	\left( \gamma_{i_1}x_1^{-1}(1 - r_0) + \alpha_{i_1}(x_1)r_0 \right) \star f \mathbf{1}_{\tuple i},
\end{align*}
for any  $\tuple i, \tuple j \in \beta$ and $f \in K[x_1, \dots, x_n]$.

\begin{lemma}
\label{lemma:pbw_action_well_defined}
The previous action is well-defined.
\end{lemma}

The proof of Lemma~\ref{lemma:pbw_action_well_defined} is given in Appendix~\ref{appendix_section:polynomial_realisation}.
For each $w \in B_n$ we now fix a reduced expression $w = r_{a_1} \cdots r_{a_k}$ and define $\psi_w \coloneqq \psi_{a_1} \cdots \psi_{a_k} \in V_\beta(\Gamma, \lambda, \gamma)$. Note that the element $\psi_w$ may depend on the chosen reduced expression.

\begin{theorem}
\label{theorem:pbw}
The algebra $V_\beta(\Gamma, \lambda, \gamma)$ is a free $K$-module, and
\[
\left\{y_1^{a_1} \cdots y_n^{a_n} \psi_w e(\tuple{i}) : a_i \in \mathbb{N},  w \in B_n, \tuple{i} \in \beta\right\},
\]
is a $K$-basis.
\end{theorem}

\begin{proof}
As in~\cite{khovanov-lauda_diagrammaticI,khovanov-lauda_diagrammaticII,rouquier_2}, successively applying the defining relations of $V_\beta(\Gamma,\lambda,\gamma)$ we can see that the above family is a spanning set, hence it remains to prove that it is linearly independent. 
For any $b \in \{0, \dots, n-1\}$, $\tuple i \in \beta$ and $f \in K[x_1, \dots, x_n]$ we can write
\[
\psi_b \cdot f \mathbf{1}_{\tuple i}= \left(A_{\tuple i}^{r_b} r_b + A_{\tuple i}^{1,r_b} \right)\star f \mathbf{1}_{\tuple i},
\]
where $A_{\tuple{i}}^{r_b},A_{\tuple i}^{1,r_b}\in K(x_1, \dots, x_n)$ with $A_{\tuple i}^{r_b}$ non-zero (recall that $P_{ij} \neq 0$ if $i \neq j$). If $<$ is the Bruhat order on $B_n$, we deduce that  for each $w \in B_n$ we can write
\[
\psi_w \cdot f \mathbf{1}_{\tuple{i}}= \left(A_{\tuple{i}}^w  w  + \sum_{w' < w} A_{\tuple{i}}^{w', w} w' \right) \star f \mathbf{1}_{\tuple i},
\]
where $A_{\tuple i}^w,A_{\tuple i}^{w',w} \in K(x_1, \dots, x_n)$ with $A_{\tuple i}^w$  non-zero. Thus,
\[
y_1^{a_1}\cdots y_n^{a_n} \psi_w \cdot f \mathbf{1}_{\tuple{i}} = \left(A_{\tuple{i}}^w x_1^{a_1} \cdots x_n^{a_n} w + \sum_{w' < w} A_{\tuple{i}}^{w', w} x_1^{a_1} \cdots x_n^{a_n} w'\right)\star f \mathbf{1}_{\tuple i},
\]
for any $a_1, \dots, a_n \in \mathbb{N}$.
%%%%%%
% Note that we can write
%\[
%\polrep(\psi_b e(\tuple i)) = A_{\tuple i}^{r_b} r_b + A_{\tuple i}^{1, r_b},
%\]
%for all $b \in \{0, \dots, n-1\}$, where $A_{\tuple i}^{r_b},A_{\tuple i}^{1, r_b} \in K(x_1, \dots, x_n)$ with $A_{\tuple i}^{r_b} \neq 0$. We deduce that  for each $w \in B_n$ we can write
%\[
%\polrep(\psi_w e(\tuple{i})) = A_{\tuple{i}}^w r_w + \sum_{w' < w} A_{\tuple{i}}^{w',w} r_{w'},
%\]
%where $A_{\tuple{i}}^w,A_{\tuple{i}}^{w', w} \in K(x_1, \dots, x_n)$ with $A_{\tuple{i}}^w \neq 0$ and $<$ is the Bruhat order on $B_n$. Thus,
%\[
%\polrep(y_1^{a_1}\cdots y_n^{a_n} \psi_w e(\tuple{i})) = A_{\tuple{i}}^w x_1^{a_1} \cdots x_n^{a_n} r_w + \sum_{w' < w} A_{\tuple{i}}^{w', w} x_1^{a_1} \cdots x_n^{a_n} r_{w'}.
%\]
We now use the following basic Lemma~\ref{lemma:dedekind} from field theory and notice that the elements of $B_n$ induce distinct field homomorphisms of $K(x_1, \dots, x_n)$.

\begin{lemma}[Dedekind]
\label{lemma:dedekind}
If $u_1, \dots, u_n : F \to G$ are distinct field homomorphisms then they form a linearly independent family over $G$.
\end{lemma}
So we can use reverse induction in the Bruhat order to show that the images of the basis elements are linearly independent in $\mathrm{End}_K(K[x, \beta])$ and thus conclude the proof.
\end{proof}

As a corollary, we obtain the sequel of Remark~\ref{remark:quotient}.

\begin{corollary}
\label{corollary:klr_subalgebra_vv} %correction de l'enonce
The  subalgebra of $V_\beta(\Gamma,\lambda,\gamma)$ generated by all generators but $\psi_0$  is isomorphic to  $\H_\beta(\Gamma)$.
\end{corollary}

\section{Disjoint quiver isomorphism}
\label{section:isomorphism B}

Let $\Gamma$ be a quiver with an involution $\theta$ and $\lambda \in \mathbb{N}^I, \gamma \in K^I$ as in~\textsection\ref{subsection:definition_QHA_typeB}. Let $d$ be a positive integer and write $\Gamma = \amalg_{j = 1}^d \Gamma^{(j)}$ such that
\begin{itemize}
\item each $\Gamma^{(j)}$ is a full subquiver of $\Gamma$;
\item each $\Gamma^{(j)}$ is stable under $\theta$.
\end{itemize}
We write $I = \amalg_{j = 1}^d I^{(j)}$ the corresponding partition of the vertex set of $\Gamma$. 
Recall that $B_n \simeq G^n \rtimes \mathfrak{S}_n$ with $G = \mathbb{Z}/2\mathbb{Z}$ acting on $I$ via $G \twoheadrightarrow \langle \theta\rangle$. In particular, each $I^{(j)}$ for $j \in \{1, \dots, d\}$ is stable under the action of $G$ so that we are in the setting of~\textsection\ref{subsection:disjoint-quiver}. 

%tout petit peu plus de details
Let $\beta$ be a $B_n$-orbit in $I^n$. As explained in~\textsection\ref{subsection:disjoint-quiver}, both properties \eqref{equation:assumption_alpha} and \eqref{equation:assumption_alpha2} are satisfied. In particular, for any $j \in \{1, \dots, d\}$ we have an integer $\sizepart j \beta = \sizeglob j$ and we have a $B_{\sizeglob j}$-orbit $\beta^{(j)} \subseteq (I^{(j)})^{\sizeglob j}$.

 For any $j \in \{1, \dots, d\}$, we define $\lambda^{(j)} \in \mathbb{N}^{I^{(j)}}$ (respectively $\gamma^{(j)} \in K^{I^{(j)}}$) to be the restriction of $\lambda$ (resp.  $\gamma$) to $I^{(j)}$.

\begin{theorem}
\label{theorem:isom_carquois_disjoints_B}
We have an (explicit) isomorphism of graded algebras
\[
V_\beta(\Gamma, \lambda, \gamma) \simeq \mathrm{Mat}_{\binom{n}{\sizeglob 1, \dots, \sizeglob d}}\left(\bigotimes_{j = 1}^d V_{\beta^{(j)}}\left(\Gamma^{(j)}, \lambda^{(j)}, \gamma^{(j)}\right )\right).
\]
\end{theorem}

As in \textsection\ref{section:toy_KLR}, will first apply the result of \textsection\ref{section:general_decomposition_theorem} and then prove an isomorphism with a tensor product. Parts~\ref{subsection:costless_isomorphism_B} and~\ref{subsection:embedding_B} are devoted to the proof of Theorem~\ref{theorem:isom_carquois_disjoints_B}, which is a direct consequence of~\eqref{equation:costless_B} and Proposition~\ref{proposition:isom_tensor_product}.

\subsection{Fixing the profile}
\label{subsection:costless_isomorphism_B}

As defined in~\textsection\ref{subsection:decomposition_KLR}, to each $\tuple i \in \beta$ we associate its profile $p(\tuple i) \in \{1, \dots, d\}^n$, and we write $\prof^\beta \subseteq \{1, \dots, d\}^n$ to denote the set of all profiles of elements of $\beta$. Any element of $\prof^\beta$ can be reordered so that we obtain
\[
\mathfrak{t}^\beta = (1, \dots, 1, \dots, d, \dots, d),
\]
where each $j \in \{1, \dots, d\}$ appears exactly $\sizeglob j$ times.
To any $\mathfrak{t} \in \prof^\beta$, we define the idempotent
\[
e(\mathfrak{t}) \coloneqq \sum_{\substack{\tuple{i} \in \beta \\ p(\tuple i) = \mathfrak{t}}} e(\tuple{i})  \in V_\beta(\Gamma, \lambda, \gamma),
\]
and we define
\[
\I \coloneqq \{e(\mathfrak{t}) : \mathfrak{t} \in \prof^\beta\}.
\]
It is a complete set of orthogonal idempotents and its cardinality is exactly $\binom{n}{\sizeglob 1, \dots, \sizeglob d}$. Since any reduced expression in $\mathfrak{S}_n$ in the generators $r_1, \dots, r_{n-1}$ is also reduced in $B_n$ for these same generators, the definitions~\eqref{subequations:def_psit_phit} make sense in $V_\beta(\Gamma,\lambda,\gamma)$ for any $\mathfrak{t} \in \prof^\beta$.  Moreover, since the defining relations of $\H_\beta(\Gamma)$ are also satisfied in $V_\beta(\Gamma,\lambda,\gamma)$, we deduce that equations~\eqref{equation:phit_psit_klr} and~\eqref{equation:psit_et_phit} are still satisfied in $V_\beta(\Gamma,\lambda,\gamma)$ thus as in~\textsection\ref{subsection:decomposition_KLR} we conclude that
\begin{equation}
\label{equation:costless_B}
V_\beta(\Gamma,\lambda,\gamma) \simeq \mathrm{Mat}_{\binom{n}{\sizeglob 1, \dots, \sizeglob d}}\left(e(\mathfrak{t}^\beta) V_\beta(\Gamma,\lambda,\gamma) e(\mathfrak{t}^\beta)\right).
\end{equation}

\subsection{Embedding the tensor product}
\label{subsection:embedding_B}

The aim of this section is to prove the following proposition.

\begin{proposition}
\label{proposition:isom_tensor_product}
We have an (explicit) isomorphism of graded algebras
\[
e(\mathfrak{t}^\beta)V_\beta(\Gamma, \lambda, \gamma)e(\mathfrak{t}^\beta) \simeq \bigotimes_{j=1 }^d V_{\beta^{(j)}}\left(\Gamma^{(j)}, \lambda^{(j)}, \gamma^{(j)}\right).
\]
\end{proposition}

\subsubsection{Images of the generators}

Set $n=n_1+\dots+n_d$. % where $n_j\coloneqq|\beta^{(j)}|$. 
We start by defining a map from the set of generators of the algebra $\bigotimes_{j=1}^d V_{\beta^{(j)}}\left(\Gamma^{(j)}, \lambda^{(j)}, \gamma^{(j)}\right)$ to  $e(\mathfrak{t}^\beta)V_\beta(\Gamma, \lambda, \gamma)e(\mathfrak{t}^\beta)$.

Let $j\in \{1,\dots,d\}$. We denote $\psi_0^{(j)},\dots,\psi_{n_j-1}^{(j)}$, $y^{(j)}_1,\dots,y_{n_j}^{(j)}$, $e(\tuple{i}^j)$ with $\tuple{i}^j\in\beta^{(j)}$, the generators of $V_{\beta^{(j)}}\left(\Gamma^{(j)}, \lambda^{(j)}, \gamma^{(j)}\right)$. Then we consider the map
\begin{align}
e(\tuple{i}^1) \otimes \dots \otimes e(\tuple{i}^d) &\mapsto e(\tuple{i}^1, \dots, \tuple{i}^d)\,,\label{map1}\\
\psi^{(j)}_0 &\mapsto e(\mathfrak{t}^\beta)\psi_{n_1 + \dots + n_{j-1}} \dots \psi_1\psi_0\psi_1 \dots \psi_{n_1 + \dots + n_{j-1}}e(\mathfrak{t}^\beta)\,,\label{map2}\\
\psi^{(j)}_a &\mapsto e(\mathfrak{t}^\beta)\psi_{n_1 + \dots + n_{j-1} + a}e(\mathfrak{t}^\beta)\,,\ \ \ \ \ \ a=1,\dots,n_j-1\,,\label{map3}\\
y_b^{(j)} &\mapsto  e(\mathfrak{t}^\beta)y_{n_1+\dots+n_{j-1} + b}e(\mathfrak{t}^\beta)\,,\ \ \ \ \ \ b=1,\dots,n_j\,,\label{map4}
\end{align}
where each $\tuple{i}^j\in\beta^{(j)}$ and $(\tuple{i}^1, \dots, \tuple{i}^d)$ is simply the concatenation. Note that  $(\tuple{i}^1, \dots, \tuple{i}^d) \in \beta$ since $\beta$ is a $B_n$-orbit, using Proposition \ref{prop-assumption2}. Moreover, the profile of $(\tuple{i}^1, \dots, \tuple{i}^d)$ is $\mathfrak{t}^\beta$ and thus $e(\tuple{i}^1, \dots, \tuple{i}^d)e(\mathfrak{t}^\beta)=e(\mathfrak{t}^\beta)e(\tuple{i}^1, \dots, \tuple{i}^d)=e(\tuple{i}^1, \dots, \tuple{i}^d)$. By convention, $n_1+\dots+n_{j-1}=0$ if $j=1$ (and $\psi_0^{(1)}\mapsto\psi_0$).
Note also that the Formula \eqref{map1} extended by linearity gives the image of an idempotent $e(\tuple i^j)\in V_{\beta^{(j)}}\left(\Gamma^{(j)}, \lambda^{(j)}, \gamma^{(j)}\right)$:
\begin{equation}
e(\tuple{i}^j)\mapsto \sum_{\substack{j'=1\\ j'\neq j}}^d\ \ \sum_{\tuple i^{j'}\in\beta^{(j')}} e(\tuple i^1,\dots,\tuple i^d)\,.\label{map1'}
\end{equation}
Equivalently, the image of $e(\tuple{i}^j)$ is the sum of the idempotents $e(\tuple i)$ where the sum is taken over $\tuple i\in\beta$ such that the profile of $\tuple i$ is $\mathfrak{t}^\beta$ and moreover $(i_{n_1+\dots+n_{j-1}+1},\dots,i_{n_1+\dots+n_j})=\tuple i^j$.

\vskip .2cm
We will prove that the map given in \eqref{map1}--\eqref{map4} extends to an homomorphism of graded algebras denoted $\rho$ and that $\rho$ is bijective.

\subsubsection{Grading}

 We check that the map given in \eqref{map1}--\eqref{map4} preserves the grading given in \eqref{gradingB}. For the images of the idempotents and of the generators $y_b^{(j)}$, there is nothing to check. 

Let $\tuple i^j\in\beta^{(j)}$ and $\tuple i\in\beta$ such that $(i_{n_1+\dots+n_{j-1}+1},\dots,i_{n_1+\dots+n_j})=\tuple i^j$.
Let $a\in\{1,\dots,n_j-1\}$. On the one hand, we have $\text{deg}\,\psi^{(j)}_ae(\tuple i^j)=d(i^j_a,i^j_{a+1})$. On the other hand, we have 
$$\text{deg}\,\psi_{n_1+\dots+n_{j-1}+a}e(\tuple i)=d(i_{n_1+\dots+n_{j-1}+a},i_{n_1+\dots+n_{j-1}+a+1})=d(i^j_a,i^j_{a+1})\ .$$

Finally, on the one hand, we have $\text{deg}\,\psi^{(j)}_0e(\tuple i^j)=d(i^j_1)$. On the other hand, we claim that we have
$$\text{deg}\,\psi_{k}\dots\psi_0\dots\psi_ke(\tuple j)=d(j_{k+1})\ ,$$
for any $k\geq 0$ and any $\tuple j\in\beta$ such that $j_{k+1}$ is not in the same component as $j_1,\dots,j_k$ for the decomposition of the quiver $\Gamma = \amalg_{j=1}^d \Gamma^{(j)}$. Taking $k=n_1+\dots+n_{j-1}$ and $\tuple j=\tuple i$, this concludes the verification.

To prove the claim, we use induction on $k$. For $k=0$, this is the definition of the degree of $\psi_0e(\tuple j)$. For $k>0$, we have $\text{deg}\,\psi_ke(\tuple j)=j_k\cdot j_{k+1}=\lvert j_k \to j_{k+1} \rvert + \lvert j_k \leftarrow j_{k+1}\rvert=0$ by assumption on $\tuple j$. Similarly, $\text{deg}\,\psi_ke(\tuple j')=0$ where $\tuple j'=r_{k-1}\dots r_0\dots r_{k-1} r_k(\tuple j)$, since $(j'_k,j'_{k+1})=(\theta(j_{k+1}),j_k)$. It remains to use the induction hypothesis, namely that $\text{deg}\,\psi_{k-1}\dots\psi_0\dots\psi_{k-1}e(r_k(\tuple j))=d(j_{k+1})$, valid since $r_k(\tuple j)$ has $j_{k+1}$ in position $k$.

\subsubsection{Bijectivity}

We assume for a moment that the map given in \eqref{map1}--\eqref{map4} extends to an algebra homomorphism. We denote this map by $\rho$ and we prove here that $\rho$ is bijective. 

For any $j \in \{1,\dots,d\}$, we write $B^{(j)} \coloneqq B_{n_j}$ and we rename its generators to $r_0^{(j)}, \dots, r_{n_j - 1}^{(j)}$. We recall the following fact.
\begin{lemma}
\label{lemma:parabolic_Bn}
We have an injective group homomorphism
\begin{align*}
B^{(1)} \times \dots \times B^{(d)} \to B_n\\
(w_1,\dots,w_d)\mapsto \overline{w}_1\dots\overline{w}_d
\end{align*}
given on the generators by, for $j\in\{1,\dots,d\}$,
\begin{align*}
r^{(j)}_0 &\mapsto r_{n_1 + \dots + n_{j-1}} \dots r_1r_0r_1 \dots r_{n_1 + \dots + n_{j-1}}\,,\\
r^{(j)}_a &\mapsto r_{n_1 + \dots + n_{j-1} + a}\,,\ \ \ \ \ \ a=1,\dots,n_j-1\,.
\end{align*}
By convention, $n_1+\dots+n_{j-1}=0$ if $j=1$ (and $r_0^{(1)}\mapsto r_0$). Moreover, any $d$-tuple of reduced expressions is sent onto a reduced expression in $B_n$.
\end{lemma} 

\begin{proof}
Recall that $B_n = \langle r_0, \dots, r_{n-1} \rangle$ is the group of signed permutations of $\{\pm 1, \dots, \pm n\}$, with $r_0 = (1,-1)$ and $r_a = (a,a+1)(-a,-a-1)$ for $a=1,\dots,n-1$. Let $t_1\coloneqq r_0$ and $t_{a+1}\coloneqq r_at_ar_a$ for $a=1,\dots,n-1$. The element $t_a$ corresponds to the transposition $(-a,a)$.

For any $i \in \{1, \dots, n\}$ and $a \in \{1,\dots,i\}$, we set by convention $r_a \dots r_{i-1} = 1$ if $a = i$.
It is easy to see (for example \cite[Figure 9]{OP}) that:
\[B_n=\bigsqcup_{a=1}^nr_a\dots r_{n-1}B_{n-1}\sqcup\bigsqcup_{a=1}^{n}t_ar_a\dots r_{n-1}B_{n-1}\ .\]
So, if we define, for $i\in\{1,\dots,n\}$,
\[R^{(i)}:=\{t_a^{\epsilon}r_a\dots r_{i-1}\ |\ a\in\{1,\dots,i\}\,,\ \epsilon\in\{0,1\}\}\ ;\]
then we have that
\begin{equation}\label{setWn}\{u_n\dots u_1\ |\ u_i\in R^{(i)}\ \}\ ,
\end{equation}
forms a complete set of pairwise distinct elements of $B_n$. Moreover this set consists of reduced expressions in terms of the generators $r_0,r_1,\dots,r_{n-1}$, since the polynomial $\sum_k a_k t^k$, where $a_k$ records the number of elements in (\ref{setWn}) written as a product of $k$ generators, is easily found to be $\prod_{i=1}^n\frac{1-t^{2i}}{1-t}$ which is the Poincar\'e polynomial $\sum_{w\in B_n}t^{\ell(w)}$ of the Coxeter group of type $B_n$ (see, for instance, \cite[Theorem 7.1.5]{BB}).

Now, to prove the lemma, we note that the subgroup permuting only the numbers $\pm1,\dots,\pm n_1$ is isomorphic to $B^{(1)}$, the subgroup permuting only the numbers $\pm (n_1+1),\dots,\pm (n_1+n_2)$ is isomorphic to $B^{(2)}$ and so on. These subgroups commute and therefore we have an embedding of $B^{(1)} \times \dots \times B^{(d)}$ inside $B_n$ (though not as a parabolic subgroup). It is straightforward to see that this corresponds to the embedding described at the level of the generators in the lemma.

For the statement about the reduced expressions, let us first recall that the length function of the Coxeter group $B_n$ can be expressed in terms on inversions as follows (see for example \cite[\S 8.1]{BB}):
\[%\begin{equation}\label{lWn}
\ell(\pi)=\sharp\left\{1\leq i<j\leq n\ \vert\ \pi(i)>\pi(j)\right\}+\sharp\left\{1\leq i\leq j\leq n\ |\ \pi(-i)>\pi(j)\right\}\ .
\]%\end{equation}
Using the notations of the lemma, we obtain that $\ell(\overline{w}_1\dots\overline{w}_d)=\ell(\overline{w}_1)+\dots+\ell(\overline{w}_d)$, since $\overline{w}_1$ permutes only the numbers $\pm1,\dots,\pm n_1$, $\overline{w}_2$ permutes only the numbers $\pm (n_1+1),\dots,\pm (n_1+n_2)$, and so on. So it remains to show that a reduced expression in $B^{(j)}$, $j=1,\dots,d$, is sent to a reduced expression in $B_n$.

Let $j\in\{1,\dots,d\}$. We claim that it is enough to show our assertion for a single reduced expression for each element of $B^{(j)}$. Indeed the number of occurrences of $r_0$ in different reduced expressions of a same element remains constant (due to the homogeneity in $r_0$ of the braid relations of $B_n$), and therefore, the number of generators in the images of these different reduced expressions is also constant. So if one of these images is reduced, they are all reduced.

Finally, to conclude the proof of the lemma, we observe that the set of reduced expressions of the form \eqref{setWn} in $B^{(j)}$ is sent to expressions of the same form in $B_n$, which are therefore reduced as well.
\end{proof}

To prove that $\rho$ is bijective, we use first that we know a basis of $\bigotimes_{j =1}^d V_{\beta^{(j)}}\bigl(\Gamma^{(j)}, \lambda^{(j)}, \gamma^{(j)}\bigr)$ by Theorem~\ref{theorem:pbw}. A basis element is of the form
\begin{equation}\label{basisprod}
\bigotimes_{j=1}^d (y_1^{(j)})^{a^{(j)}_1} \dots (y_{n_j}^{(j)})^{a^{(j)}_{n_j}} \psi_{w_j}^{(j)} e(\tuple{i}^j),
\end{equation}
where $a^{(j)}_{1},\dots a^{(j)}_{n_j}\in\mathbb{N}$, $\tuple i^j \in \beta^{(j)}$ and $w_j\in B^{(j)}$. Note that we have fixed a reduced expression for each element $w_j\in B^{(j)}$ for each $j=1,\dots,d$, in order to define $\psi_{w_j}^{(j)}$.

On the other hand, we also know a basis of $e(\mathfrak{t}^\beta)V_\beta(\Gamma, \lambda, \gamma)e(\mathfrak{t}^\beta)$ again by Theorem~\ref{theorem:pbw}. Indeed note that $e(\tuple i)e(\mathfrak{t}^\beta)=e(\tuple i)$ if the profile of $\tuple i$ is $\mathfrak{t}^\beta$ and $e(\tuple i)e(\mathfrak{t}^\beta)=0$ otherwise. Moreover, $\psi_w e(\tuple i)=e(w\cdot \tuple i)\psi_w$. So it is straightforward to conclude that a basis element of $e(\mathfrak{t}^\beta)V_\beta(\Gamma, \lambda, \gamma)e(\mathfrak{t}^\beta)$ is of the form
\begin{equation}\label{basisidem}
y_1^{a_1} \dots y_n^{a_n} \psi_{w} e(\tuple i)\,,
\end{equation}
where $a_{1},\dots, a_{n}\in\mathbb{N}$, $\tuple i\in \beta$ with profile $\mathfrak{t}^\beta$ and $w$ is in the subgroup of $B_n$ isomorphic to $B^{(1)} \times \dots \times B^{(d)}$ from Lemma \ref{lemma:parabolic_Bn} (the stabiliser of $\mathfrak{t}^\beta$). We must fix reduced expressions for such $w$ in order to define $\psi_w$. We fix them as the images of the reduced expressions of elements $B^{(1)} \times \dots \times B^{(d)}$ chosen in the preceding paragraph. That we can do so is the last statement in Lemma \ref{lemma:parabolic_Bn}.

Finally, the image of a basis element \eqref{basisprod} under the homomorphism $\rho$ is
\begin{equation}
\label{equation:element_tensor_inside_full_algebra}
y_1^{b_1} \dots y_n^{b_n} \psi_{\overline{w}_1} \cdots \psi_{\overline{w}_d} e(\tuple i^1,\dots,\tuple i^d),
\end{equation}
where $b_{n_1 + \dots + n_{j-1} + k}=a^{(j)}_{k}$ and the notation $\overline{w}_j$ comes from Lemma~\ref{lemma:parabolic_Bn}. The concatenation $(\tuple i^1,\dots,\tuple i^d)$ has the profile $\mathfrak{t}^\beta$ since each $\tuple i^j\in\beta^{(j)}$, and due to our choices of reduced expressions, we have $\psi_{\overline{w}_1} \cdots \psi_{\overline{w}_d}=\psi_{\overline{w}_1\cdots \overline{w}_d}$. So we conclude that the element~\eqref{equation:element_tensor_inside_full_algebra} is of the form \eqref{basisidem}. Further, it is immediate that we can obtain in this way all the basis elements of $e(\mathfrak{t}^\beta)V_\beta(\Gamma, \lambda, \gamma)e(\mathfrak{t}^\beta)$. We conclude that the homomorphism $\rho$ sends a basis onto a basis and thus is bijective.

\subsubsection{Homomorphism property}
To finish the proof of Proposition~\ref{proposition:isom_tensor_product}, it remains to check that the map defined in \eqref{map1}--\eqref{map4} extends to an algebra homomorphism. It is possible but quite lengthy to check explicitly that all defining relations are preserved. Instead we are going to use the polynomial representation introduced in \textsection\ref{subsection:basis_theorem}. We keep in use the notations introduced in \textsection\ref{subsection:basis_theorem}.

From the proof of Theorem \ref{theorem:pbw}, we see that the action of the algebra $V_\beta(\Gamma, \lambda, \gamma)$ on $K[x,\beta]$ is faithful, or in other words, we have an embedding of $V_\beta(\Gamma, \lambda, \gamma)$ in $\mathrm{End}_K(K[x, \beta])$. Therefore, if we denote $\phi\bigl(e(\mathfrak{t}^\beta)\bigr)$ the image of $e(\mathfrak{t}^\beta)$ by this embedding, we obtain an embedding of the algebra $e(\mathfrak{t}^\beta)V_\beta(\Gamma, \lambda, \gamma)e(\mathfrak{t}^\beta)$ in $\mathrm{End}_K\bigl(\phi\bigl(e(\mathfrak{t}^\beta)\bigr)K[x, \beta]\bigr)$. We have immediately:
\begin{equation}\label{space1}
\phi\bigl(e(\mathfrak{t}^\beta)\bigr)K[x, \beta]=\bigoplus_{\substack{\tuple i\in\beta\\ p(\tuple i)=\mathfrak{t}^\beta}}K[x_1,\dots,x_n]\mathbf{1}_{\tuple{i}}\ .
\end{equation}

On the other hand, we also have an embedding of the algebra $\bigotimes_{j=1 }^d V_{\beta^{(j)}}\left(\Gamma^{(j)}, \lambda^{(j)}, \gamma^{(j)}\right)$ in $\mathrm{End}_K(\bigotimes_{j=1}^dK[x, \beta^{(j)}])$, and we have the natural identification:
\begin{equation}\label{space2}
\bigotimes_{j=1}^dK[x, \beta^{(j)}]=\bigotimes_{j=1}^d\bigoplus_{\tuple i^j\in\beta^{(j)}}K[x_1^{(j)},\dots,x_{\sizeglob j}^{(j)}]\mathbf{1}_{\tuple{i}^j}\cong \bigoplus_{\substack{\tuple i\in\beta\\ p(\tuple i)=\mathfrak{t}^\beta}}K[x_1,\dots,x_n]\mathbf{1}_{\tuple{i}}\ .
\end{equation}
The identification simply maps $f_1\mathbf{1}_{\tuple{i}^1}\otimes \dots \otimes f_d\mathbf{1}_{\tuple{i}^d}$ to $f_1\dots f_d\mathbf{1}_{(\tuple{i}^1,\dots,\tuple i^d)}$.

Through the identifications we just made, both algebras related by the map in \eqref{map1}--\eqref{map4} are seen as algebras of endomorphisms of the same space, in \eqref{space1} and \eqref{space2}. So in order to check the homomorphism property, it is enough to check that both sides of Formulas \eqref{map1}--\eqref{map4} are in fact the same elements in the endomorphism algebra.

This verification is immediate for \eqref{map1}--\eqref{map2} and \eqref{map4}. For the image of $\psi_0^{(j)}$, we proceed as follows. First, it is convenient to choose a polynomial representation as in \textsection\ref{subsection:basis_theorem} for which $P_{ij}(u, v) \coloneqq (u - v)^{\lvert j \to i\rvert}$ if $i\neq j $ and $P_{ij}(u, v) \coloneqq 0$ if $i=j$.

Let $\tuple i\in\beta$ such that $p(\tuple i)=\mathfrak{t}^\beta$. It means that $\tuple i=(\tuple i^1,\dots,\tuple i^d)$ where $\tuple i^j\in\beta^{(j)}$. Fix $j\in\{1,\dots,d\}$ and set for brevity $k=\sizeglob 1+\dots+\sizeglob{j-1}$. Through the identifications explained above, the action of $\psi_0^{(j)}$ is given by:
\[f\mathbf{1}_{\tuple{i}}\mapsto \left( \gamma_{i_{k+1}}\frac{f - \actpol{r^{(j)}_0}{f}}{x_{k+1}} + \alpha_{i_{k+1}}(x_{k+1})\actpol{r^{(j)}_0}{f}\right) \mathbf{1}_{r^{(j)}_0 \cdot \tuple{i}}\ ,\]
where we recall that $r_0^{(j)}=r_k\dots r_1r_0r_1\dots r_k$ acts on $\tuple i$ simply by replacing $i_{k+1}$ by $\theta(i_{k+1})$.

On the other hand, we need to calculate the action of $\psi_{k} \dots \psi_1\psi_0\psi_1 \dots \psi_{k}$. We note that, with our choice of $P_{ij}(u,v)$, we have that $P_{ij}(u,v)=1$ if one index is among $\{i_1,\dots,i_k\}$ and the other is $i_{k+1}$ or $\theta(i_{k+1})$. Indeed, $i_{k+1}$ is not in the same connected component of the quiver than $i_1,\dots,i_k$ since $p(\tuple i)=\mathfrak{t}^\beta$. This is also true for $\theta(i_{k+1})$ since $\theta$ leaves stable the set $I^{(j)}$.

Then the calculation is made in three steps, corresponding respectively to the action of $\psi_1 \dots \psi_{k}$, the action of $\psi_0$ and the action of $\psi_k \dots \psi_{1}$:
\begin{multline*}
f\mathbf{1}_{\tuple{i}}  \mapsto \actpol{r_1\dots r_k}{f}\mathbf{1}_{r_1\dots r_k\cdot\tuple{i}}
\\
\mapsto \left( \gamma_{i_{k+1}}\frac{\actpol{r_1\dots r_k}{f} - \actpol{r_0r_1\dots r_k}{f}}{x_{1}} + \alpha_{i_{k+1}}(x_{1})\actpol{r_0r_1\dots r_k}{f}\right)\mathbf{1}_{r_0r_1\dots r_k\cdot\tuple{i}}
\\
\mapsto \left( \gamma_{i_{k+1}}\frac{f - \actpol{r^{(j)}_0}{f}}{x_{k+1}} + \alpha_{i_{k+1}}(x_{k+1})\actpol{r^{(j)}_0}{f}\right) \mathbf{1}_{r^{(j)}_0 \cdot \tuple{i}}\ .
\end{multline*}
This concludes the verification of the homomorphism property and the proof of Proposition~\ref{proposition:isom_tensor_product}.

\subsection{Cyclotomic quotients}
\label{subsection:cyclotomic_case_B}

As in~\textsection\ref{subsection:cyclotomic_klr}, let $\Lambda = (\Lambda_i)_{i \in I}$ be a finitely-supported family of non-negative integers. In the same way as~\cite{varagnolo-vasserot_canonical,poulain_dandecy-walker_B,poulain_dandecy-walker_D}, we define the  cyclotomic quotient of the algebra $V_\beta(\Gamma,\lambda,\gamma)$.

\begin{definition}
\label{def-cycloB}
We define the algebra $V_\beta^\Lambda(\Gamma,\lambda,\gamma)$ as the quotient of $V_\beta(\Gamma,\lambda,\gamma)$ by the two-sided ideal $\idealB_\beta^{\Lambda}$ generated by the relations
\[
y_1^{\Lambda_{i_1}}e(\tuple i) = 0\,,\ \ \ \ \ \ \text{for all $\tuple i = (i_1, \dots, i_n) \in \beta$.}
\]
\end{definition}

The above relations are homogeneous so that $V_\beta^\Lambda(\Gamma,\lambda,\gamma)$ is graded. Note that if $\Lambda_i = 0$ for all $i$ then
\[
V_\beta^\Lambda(\Gamma,\lambda,\gamma) =
\begin{cases}
\{0\}, &\text{if } n \geq 1,
\\
K, &\text{if } n = 0.
\end{cases}
\]

As in~\textsection\ref{subsection:cyclotomic_klr}, for any $j \in \{1, \dots, d\}$ let $\Lambda^{(j)}$  be the restriction of $\Lambda$ to the vertex set $I^{(j)}$ of~$\Gamma^{(j)}$.

\begin{corollary}
\label{corollary:isom_disjoint_quiver_cyclo_B}
We have an (explicit) isomorphism of graded algebras:
\[
V_\beta^{\Lambda}(\Gamma, \lambda, \gamma) \simeq \mathrm{Mat}_{ \binom{n}{\sizeglob 1, \dots, \sizeglob d}}\left(\bigotimes_{j = 1}^d V_{\beta^{(j)}}^{\Lambda^{(j)}}\left(\Gamma^{(j)}, \lambda^{(j)}, \gamma^{(j)}\right )\right).
\]
\end{corollary}

\begin{proof} % petites modifications de la preuve
Recall that the algebra $\H_\beta(\Gamma)$ is isomorphic to a subalgebra of $V_\beta(\Gamma,\lambda,\gamma)$ (see Corollary \ref{corollary:klr_subalgebra_vv}). Moreover, if $\vartheta$ denotes the isomorphism of Theorem~\ref{theorem:isom_carquois_disjoints_B}, its restriction to $\H_\beta(\Gamma)$ is by construction the isomorphism of Corollary \ref{corollary:isom_dec_klr_block}. Therefore it is immediate that the calculations made in the proof of Theorem \ref{theorem:isom_dec_klr_cyclo_block} can be repeated verbatim here. They show that, if we denote $\idealB_{\beta,\otimes}^{\Lambda}$ the ideal of $\otimes_{j = 1}^d V_{\beta^{(j)}}(\Gamma^{(j)}, \lambda^{(j)}, \gamma^{(j)})$ such that the quotient is $\otimes_{j = 1}^d V_{\beta^{(j)}}^{\Lambda^{(j)}}(\Gamma^{(j)}, \lambda^{(j)}, \gamma^{(j)})$ (see the proof of Theorem \ref{theorem:isom_dec_klr_cyclo_block}), we have
\[
\vartheta(\idealB_\beta^{\Lambda}) = \mathrm{Mat}_{ \binom{n}{\sizeglob 1, \dots, \sizeglob d}}\left(\idealB_{\beta,\otimes}^{\Lambda}\right)\ .
\]
This concludes the proof.
\end{proof}

%definir V_n^{\Lambda}
We define $V_n^{\Lambda}(\Gamma,\lambda,\gamma):=\bigoplus_{\beta}V_\beta^\Lambda(\Gamma,\lambda,\gamma)$ where the direct sum is over the $B_n$-orbits $\beta$ in~$I^n$. As in Corollary~\ref{corollary:isom_dec_klr_n}, using the bijection~\eqref{proposition:decomposition_orbits}  we deduce the following corollary. Note that we now use  \eqref{proposition:decomposition_orbits} with $G=\mathbb{Z}/2\mathbb{Z}$.

\begin{corollary}
\label{corollary:isom_disjoint_quiver_cyclo_B_n}
We have an (explicit) isomorphism of graded algebras:
\[
V_n^{\Lambda}(\Gamma,\lambda,\gamma) \simeq \bigoplus_{\substack{n_1, \dots, n_d \geq 0 \\ n_1 + \dots + n_d = n}} \mathrm{Mat}_{\binom{n}{n_1, \dots, n_d}} \left(\bigotimes_{j = 1}^d V_{n_j}^{\Lambda^{(j)}}(\Gamma^{(j)},\lambda^{(j)},\gamma^{(j)})\right).
\]
\end{corollary}

\begin{remark}
\label{remark:cancellation_n_B}
As in Remark~\ref{rem-cancellation_n}, we deduce that we can assume that $\Lambda$ is supported on all components of $\Gamma$.
\end{remark}

\section{Quiver Hecke algebras for type D}
\label{section:QHA_D}

To fit with the setting of~\cite{poulain_dandecy-walker_D}, we now assume that $K$ is a field with $char(K)\neq 2$.

Let $\Gamma$ be a quiver with an involution $\theta$ as in~\textsection\ref{subsection:definition_QHA_typeB} and let $\beta$ be a $B_n$-orbit in $I^n$.  As before, let $\Lambda = (\Lambda_i)_{i \in I}$ be a finitely-supported family of non-negative integers.

In this section, as in Remark~\ref{remark:coherence_PAW_typeD} we consider the situation $\lambda_i=\gamma_i=0$ for all $i\in I$, and we denote simply $V_{\beta}(\Gamma)=V_{\beta}(\Gamma,0,0)$ the resulting algebra, defined in Section \ref{subsection:definition_QHA_typeB} (note that Conditions~\eqref{subequations:conditions_gamma} are satisfied with this choice of $\lambda$ and $\gamma$).
The defining relations \eqref{relation:psi0_e(i)}--\eqref{relation:braid4} (those involving the generator $\psi_0$) become simply:
\begin{align}
\psi_0 e(\tuple{i}) &= e(r_0 \cdot \tuple{i}) \psi_0,
\label{relation:psi0_e(i)2}
\\
\psi_0 \psi_b &= \psi_b \psi_0, \qquad\text{for all } b \in \{2, \dots, n-1\},
\label{relation:psi0_psib2}
\\
\label{relation:psi0_y12}
\psi_0 y_1 & = -y_1 \psi_0,
\\
\label{relation:psi0_yj2}
\psi_0 y_a &= y_a \psi_0, \qquad\text{for all } a \in \{2, \dots, n\},
\\
\label{relation:psi0square2}
\psi_0^2 & =1\,,
\\
\label{relation:braid42}
(\psi_0 \psi_1)^2  & = (\psi_1 \psi_0)^2\ .
\end{align}
So we see immediately that we have an homogeneous involutive algebra automorphism $\iota$ of $V_{\beta}(\Gamma)$ given on the generators by:
\begin{equation}\label{iota}
\iota(\psi_0)=-\psi_0\ \ \ \ \ \text{and}\ \ \ \ \iota(X)=X\ \ \ \text{for $X\in\{\psi_1,\dots,\psi_{n-1},y_1,\dots,y_n\}\cup\{e(\tuple i)\}_{\tuple i\in\beta}$\ .}
\end{equation}
Note that $\iota$ is the identity map if $n = 0$.
We denote by $V_{\beta}(\Gamma)^{\iota}$ the fixed-point subalgebra of $V_{\beta}(\Gamma)$, that is, $V_{\beta}(\Gamma)^{\iota}=\{x\in V_{\beta}(\Gamma)\ |\ \iota(x)=x\}$. The subalgebra $V_{\beta}(\Gamma)^{\iota}$ is a graded subalgebra of $V_{\beta}(\Gamma)$ since $\iota$ is homogeneous.

\paragraph{Cyclotomic quotients.} We recall that $V_\beta^\Lambda(\Gamma)$ is the quotient of $V_{\beta}(\Gamma)$ by the two-sided ideal $\idealB_\beta^\Lambda$ generated by
\[y_1^{\Lambda_ {i_1}}e(\tuple i) = 0\,,\ \ \ \ \ \ \text{for all $\tuple i\in\beta$.}\]
These relations are homogeneous so that the algebra $V_\beta^\Lambda(\Gamma)$ inherits the grading of $V_{\beta}(\Gamma)$.
The same formulas as in \eqref{iota} define an homogeneous involutive algebra automorphism of $V_{\beta}^\Lambda(\Gamma)$, and we make the slight abuse of notation of keeping the name $\iota$ for this automorphism. The fixed-point subalgebra is denoted $V_{\beta}^\Lambda(\Gamma)^{\iota}$.

\subsection{Definition and main property of \texorpdfstring{$W_{\delta}(\Gamma)$}{Wdelta(Gamma)}} %Des modifs et des ajouts jusqu'a la fin de la section
\label{subsection:def_W}

We recall some definitions and the results we need from \cite{poulain_dandecy-walker_D}.

\medskip
If $n\geq 2$, we identify the Weyl group $D_n$ of type D as the subgroup of $B_n$ generated by $s_0\coloneqq r_0r_1r_0$, $s_1\coloneqq r_1$, $\dots$, $s_{n-1}\coloneqq r_{n-1}$. The convention we need here is that $D_n = \{1\}$ if $n \in \{0,1\}$. The group $D_n$ then acts on~$I^n$ by, if $n \geq 2$,
\begin{align*}
s_0\cdot (i_1,i_2,\dots,i_n)&=(\theta(i_2),\theta(i_1),i_3,\dots,i_n)\,,\\
s_a\cdot (\dots,i_a,i_{a+1},\dots)&=(\dots,i_{a+1},i_{a},\dots)\ \ \ \ \ a=1,\dots,n-1.
\end{align*}
Let $\delta$ be a finite subset of $I^n$ stable by $D_n$, that is a finite union of $D_n$-orbits.
\begin{definition}
Let $n\geq 2$. The algebra $W_\delta(\Gamma)$ is the unitary associative $K$-algebra generated by elements
\[
\{y_a\}_{1 \leq a \leq n} \cup
\{\psi_b\}_{1 \leq b \leq n-1}  \cup  \{\Psi_0\} \cup \{e(\tuple{i})\}_{\tuple{i} \in \delta},
\]
with the relations \eqref{relation:idempotents}--\eqref{relation:psi_tresse3} of Section~\ref{section:toy_KLR} involving all the generators but $\Psi_0$,  together with
\begin{align}
\Psi_0 e(\tuple{i}) &= e(s_0 \cdot \tuple{i}) \Psi_0,
\label{relation:psi0_e(i)D}
\\
\Psi_0 \psi_b &= \psi_b \Psi_0, \qquad\text{for all } b \in \{1, \dots, n-1\}\text{ with $b\neq 2$\,},
\label{relation:psi0_psibD}
\\
\label{relation:psi0_y1D}
(\Psi_0 y_a + y_{r_1(a)} \Psi_0)e(\tuple{i}) &= \left\{\begin{array}{ll}
e(\tuple{i}) & \text{if $\theta(i_1)=i_2$,}\\
0 & \text{otherwise,}
\end{array}\right.\ \ \ \ \text{for $a\in\{1,2\}$,}
\\
\label{relation:psi0_yjD}
\Psi_0 y_a &= y_a \Psi_0, \qquad\text{for all } a \in \{3, \dots, n\},
\\
\label{relation:psi0squareD}
\Psi_0^2 e(\tuple{i}) &= Q_{\theta(i_1),i_2}(-y_1,y_2)e(\tuple i),
\\
\label{relation:braid3D}
(\Psi_0 \psi_2\Psi_0 - \psi_2\Psi_0\psi_2)e(\tuple{i}) &= \begin{dcases}
\frac{Q_{\theta(i_1),i_2}(-y_1,y_2) - Q_{\theta(i_1),i_2}(y_3,y_2)}{y_1+y_3}e(\tuple i),
&\text{if } \theta(i_1) = i_3,
\\
0, &\text{otherwise},
\end{dcases}
\end{align}
for all $\tuple i \in \delta$. 
\end{definition}

By convention, we set $W_\delta(\Gamma) = \H_\delta(\Gamma)$ if $n \in \{0, 1\}$. Explicitly, $W_\delta(\Gamma)=K$ if $n=0$ and $W_\delta(\Gamma)=\sum_{\tuple i\in\delta}K[y_1]e(\tuple i)$ if $n=1$. This choice for $n\in \{0, 1\}$ is important for the statements of the results in the next subsection.

\begin{remark}
\label{remark:coherence_PAW_D_W}
With the choices of $\Gamma$, $\theta$ and the notations of Remark~\ref{remark:coherence_PAW}, the algebra $W_\delta(\Gamma)$ is exactly the algebra $W_{\tuple x}^\delta$ defined in~\cite{poulain_dandecy-walker_D}.
\end{remark}

 The algebra $W_\delta(\Gamma)$ is $\mathbb{Z}$-graded with
\begin{align*}
\deg e(\tuple{i}) &= 0,
\\
\deg y_a &= 2,
\\
\deg \Psi_0 e(\tuple{i}) &= d(\theta(i_1),i_2),
\\
\deg \psi_b e(\tuple{i}) &= d(i_b, i_{b+1}).
\end{align*}

\begin{definition}
The cyclotomic quotient $W^\Lambda_\delta(\Gamma)$ is the quotient of the algebra $W_\delta(\Gamma)$ by the %two-sided ideal $\idealD_\delta^\Lambda$ generated by the
relations
\[y_1^{\Lambda_ {i_1}}e(\tuple i) = 0\,,\ \ \ \ \ \ \text{for all $\tuple i\in\delta$.}\]
\end{definition}
The algebra $W^\Lambda_\delta(\Gamma)$ inherits the grading from $W_\delta(\Gamma)$ since the additional relations are homogeneous. If $\Lambda_i = 0$ for all $i$ then
\[
W_\delta^\Lambda(\Gamma) = 
\begin{cases}
\{0\}, &\text{if } n \geq 1,
\\
K, &\text{if } n = 0.
\end{cases}
\]

\paragraph{Fixed-point isomorphism.} Let $\beta$ be a $B_n$-orbit in $I^n$. Note that $\beta$ is a finite union of $D_n$-orbits, so that both algebras $V_{\beta}(\Gamma)$ and $W_{\beta}(\Gamma)$ are defined.

We recall the following results from \cite{poulain_dandecy-walker_D}. Note that they were proved for a particular choice of $\Gamma$ and $\theta$ (the one relevant for the next section). However, the proof does not depend on this choice and can be repeated verbatim in our general setting.
\begin{proposition}[\cite{poulain_dandecy-walker_D}]
\label{prop-fixed-points}
\leavevmode % to fix problems when using \ref with the label
\begin{enumerate}[label=(\roman*)]
\item
\label{item:fixed_point_VV}
The algebra $W_{\beta}(\Gamma)$ is isomorphic to the subalgebra $V_{\beta}(\Gamma)^{\iota}$ of $V_\beta(\Gamma)$. 

\item
\label{item:fixed_points_VV_cyclo}
Assume that $\Lambda$ satisfies $\Lambda_{\theta(i)}=\Lambda_i$ for all $i\in I$. The cyclotomic quotient $W_{\beta}^\Lambda(\Gamma)$ is isomorphic to $V_{\beta}^{\Lambda}(\Gamma)^{\iota}$.
\end{enumerate}

%\medskip
In both cases, an isomorphism is given by $\Psi_0\mapsto \psi_0\psi_1\psi_0$ and $X\mapsto X$ for all the generators~$X$ but $\Psi_0$.
\end{proposition}

\begin{remark}
Note that it is assumed in~\cite{poulain_dandecy-walker_D} that $n \geq 2$. With our conventions, the statements are also true for $n \in \{0, 1\}$, in which cases the verification is straightforward.
\end{remark}

\begin{remark}% ajout d'une remarque
\label{remark:isom_Lambdatilde}
Recall the defining relations~\eqref{relation:psi0_e(i)2}, \eqref{relation:psi0_y12} and~\eqref{relation:psi0square2}  of $V^{\Lambda}_{\beta}(\Gamma)$. Conjugating the cyclotomic relations of $V^\Lambda_\beta(\Gamma)$ by $\psi_0$, we obtain that $y_1^{\Lambda_{i_1}}e(r_0\cdot \tuple i)=0$ for any $\tuple i \in\beta$. From this remark, it is easy to see that we have in fact $V^{\Lambda}_{\beta}(\Gamma)=V^{\widetilde{\Lambda}}_{\beta}(\Gamma)$, where $\widetilde{\Lambda}$ is now given by $\widetilde{\Lambda}_i=\text{min}\{\Lambda_i,\Lambda_{\theta(i)}\}$.
This phenomenon does not necessarily occur also in $W^{\Lambda}_{\beta}(\Gamma)$ (where $\psi_0$ is not present) and this explains the assumptions on $\Lambda$ in Proposition~\ref{prop-fixed-points}\ref{item:fixed_points_VV_cyclo}.
\end{remark}

We note that the isomorphisms given in the preceding proposition are isomorphisms of graded algebras. Indeed, in $V_{\beta}(\Gamma)$ we have $\deg \psi_0=0$ and so it is straightforward to check that the given map is homogeneous.

From Proposition~\ref{prop-fixed-points}\ref{item:fixed_point_VV} and Corollary \ref{corollary:klr_subalgebra_vv}, one obtains immediately the following statement.
\begin{corollary}
\label{corollary:klr_subalgebra_W}
The  subalgebra of $W_\beta(\Gamma)$ generated by all generators but $\Psi_0$  is isomorphic to $\H_\beta(\Gamma)$.
\end{corollary}

\paragraph{Semi-direct product.} In this paragraph, assume that $n\geq 1$. Since $\iota$ is involutive, the vector space $V_\beta(\Gamma)$ decomposes as
\[V_\beta(\Gamma)=V_\beta(\Gamma)^{\iota}\oplus V_\beta(\Gamma)^-\ ,\]
where $V_\beta(\Gamma)^-$ is the eigenspace of $\iota$ for the eigenvalue $-1$. Moreover, the generator $\psi_0$ is invertible (in fact, $\psi_0^2=1$) and satisfies $\iota(\psi_0)=-\psi_0$. So the multiplication by $\psi_0$ provides an isomorphism of vector spaces between $V_\beta(\Gamma)^{\iota}$ and $V_\beta(\Gamma)^-$, so that $V_\beta(\Gamma)^-$ can be written $V_\beta(\Gamma)^{\iota}\psi_0$. Working out the multiplication in $V_\beta(\Gamma)$
\[(x+y\psi_0)(x'+y'\psi_0)=xx'+y\psi_0y'\psi_0+(y\psi_0x'\psi_0+xy')\psi_0\ ,\]
one obtains as a standard consequence that $V_\beta(\Gamma)$ is isomorphic to the semi-direct product $V_\beta(\Gamma)^{\iota}\rtimes C_2$, where the action of the cyclic group $C_2$ of order 2 on $V_\beta(\Gamma)^{\iota}$ is by conjugation by $\psi_0$. Recall that as a vector space $V_\beta(\Gamma)^{\iota}\rtimes C_2$ is the tensor product $V_\beta(\Gamma)^\iota \otimes K[C_2]$, and the multiplication is given by
\[(x \otimes \psi_0^{\epsilon})(x'\otimes\psi_0^{\epsilon'})=(x\psi_0^{\epsilon}x'\psi_0^{\epsilon})\otimes\psi_0^{\epsilon+\epsilon'}\ .\]

Then we formulate the preceding standard facts taking into account Proposition \ref{prop-fixed-points}. First we give explicitly the automorphism of $W_{\beta}(\Gamma)$ induced by conjugation by $\psi_0$ in $V_{\beta}(\Gamma)$. We denote this automorphism of order 2 by $\pi$. It is given on the generators by:
\begin{equation}\label{auto-pi}
\pi\ :\ \ \ \Psi_0\mapsto \psi_1\,,\ \ \ \psi_1\mapsto\Psi_0\,,\ \ \ \ y_1\mapsto -y_1\,,\ \ \ \ \ e(\tuple i)\mapsto e(r_0\cdot \tuple i)\ ,
\end{equation}
and the identity on all the other generators. As a consequence of Proposition \ref{prop-fixed-points} together with the preceding discussion, we conclude that
\[V_{\beta}(\Gamma)\simeq W_{\beta}(\Gamma)\rtimes \langle\pi\rangle\ ,\]
and similarly, for $\Lambda$ as in Proposition~\ref{prop-fixed-points}\ref{item:fixed_points_VV_cyclo},
\begin{equation}
\label{equation:VLambda_sdirect_product}
V^{\Lambda}_{\beta}(\Gamma)\simeq W^{\Lambda}_{\beta}(\Gamma)\rtimes \langle\pi\rangle\ ,
\end{equation}
where we still denote by $\pi$ the automorphism of order 2 of $W^{\Lambda}_{\beta}(\Gamma)$ given by the same formulas~\eqref{auto-pi}. This is indeed an automorphism since $\Lambda$ satisfies the assumption of Proposition~\ref{prop-fixed-points}\ref{item:fixed_points_VV_cyclo}.

With these descriptions as semi-direct products, the involution $\iota$ on $V_{\beta}(\Gamma)$ (and on $V^{\Lambda}_{\beta}(\Gamma)$) is simply given by:
\begin{equation}\label{iota-semi}
\iota\bigl(x\otimes\pi^{\epsilon}\bigr)=(-1)^{\epsilon} x\otimes \pi^{\epsilon}\ ,
\end{equation}
where $\epsilon\in\{0,1\}$ and $x\in W_{\beta}(\Gamma)$ (or $x\in W^{\Lambda}_{\beta}(\Gamma)$).

\subsection{Disjoint quiver isomorphism for \texorpdfstring{$W_{\delta}(\Gamma)$}{Wdelta(Gamma)}}
\label{subsection:disjoint_quiver_isom_W}

Now let $d$ be a positive integer and assume that the quiver $\Gamma $ admits a decomposition $\Gamma = \amalg_{j = 1}^d \Gamma^{(j)}$ as in~\textsection\ref{section:isomorphism B}. Let $\beta$ be a $B_n$-orbit in $I^n$.  As in~\textsection\ref{section:isomorphism B}, for any $j \in \{1,\dots,d\}$, we have an integer $\sizepart j \beta = \sizeglob j$ and a $B_{\sizeglob j}$-orbit $\beta^{(j)}$ in $(I^{(j)})^{\sizeglob j}$. 

If $n_j(\beta)=0$ for some $j\in\{1,\dots,d\}$ then consider $\tilde{\Gamma}$ the quiver where we removed the component $\Gamma^{(j)}$. It is immediate from the definitions that $W_\beta(\Gamma)$ is the same algebra as $W_\beta(\tilde{\Gamma})$. So we lose no generality by assuming that $n_j(\beta)\neq 0$ for all $j\in\{1,\dots,d\}$.

\paragraph{Fixed points of tensor products.} Since $n_j(\beta) \geq 1$ for all $j \in \{1, \dots, d\}$, by the preceding section we have $V_{\beta^{(j)}}(\Gamma^{(j)}) \simeq W_{\beta^{(j)}}(\Gamma^{(j)}) \rtimes C_2$ for all $j$. Hence,
\[\bigotimes_{j=1}^dV_{\beta^{(j)}}(\Gamma^{(j)})\simeq \Bigl(\bigotimes_{j=1}^dW_{\beta^{(j)}}(\Gamma^{(j)})\Bigr)\rtimes C_2^d\ ,\]
where $C_2^d$ acts on the tensor product by the automorphism $\pi$ from \eqref{auto-pi} on each factor.

We would like to describe the fixed points of $\bigotimes_{j=1}^dV_{\beta^{(j)}}(\Gamma^{(j)})$ for the involutive automorphism $\iota^{\otimes}$ given by the tensor product of $\iota$ for each factor. From Formula \eqref{iota-semi}, it is immediate to see that
\begin{equation}\label{fixed-points-tensor}
\Bigl(\bigotimes_{j=1}^dV_{\beta^{(j)}}(\Gamma^{(j)})\Bigr)^{\iota^{\otimes}}\simeq \Bigl(\bigotimes_{j=1}^dW_{\beta^{(j)}}(\Gamma^{(j)})\Bigr)\rtimes C_2^{d-1}\ ,
\end{equation}
where $C_2^{d-1}$ is seen as  the subgroup of ``even'' elements of $C_2^d$, namely
\begin{equation}\label{subgroup}
C_2^{d-1}=\{(\pi^{\epsilon_1},\dots,\pi^{\epsilon_d})\in C_2^d\ \ \text{such that } \epsilon_1+\dots+\epsilon_d=0 \pmod 2\,\}\ .
\end{equation}
%The action of this subgroup on $\bigotimes_{j=1}^dW_{\beta^{(j)}}(\Gamma^{(j)})$ is by restriction of the action of $C_2^d$. %Note that we used the assumption $n_j(\beta)\neq 0$ in order to have a description of each $V_{\beta^{(j)}}(\Gamma^{(j)})$ as a semi-direct product.

\paragraph{Disjoint quiver isomorphism.}
We can now formulate the  main result of this section. Recall that $n_j(\beta)\neq 0$ for all $j \in \{1,\dots,d\}$.
\begin{theorem}
\label{theorem:isom_carquois_disjoints_D}
 We have (explicit) isomorphisms of graded algebras:
\begin{equation}
\label{equation:isom_carquois_disjoint_D}
W_\beta(\Gamma) \simeq \mathrm{Mat}_{\binom{n}{\sizeglob 1, \dots, \sizeglob d}}\left(\Bigl(\bigotimes_{j=1}^dW_{\beta^{(j)}}(\Gamma^{(j)})\Bigr)\rtimes C_2^{d-1}\right)\ ,
\end{equation}
and, assuming $d>1$,
\begin{equation}
\label{equation:isom_carquois_disjoint_D_cyclo}
W^\Lambda_\beta(\Gamma) \simeq \mathrm{Mat}_{\binom{n}{\sizeglob 1, \dots, \sizeglob d}}\left(\Bigl(\bigotimes_{j=1}^dW^{\widetilde{\Lambda}^{(j)}}_{\beta^{(j)}}(\Gamma^{(j)})\Bigr)\rtimes C_2^{d-1}\right)\,,
\end{equation}
where  $\widetilde{\Lambda} = (\widetilde{\Lambda}_i)_{i \in I}$ is defined by $\widetilde{\Lambda}_i \coloneqq \min\{\Lambda_i, \Lambda_{\theta(i)}\}$. 
\end{theorem}

Note that in both formulas above, the group $C_2^{d-1}$ is as given in \eqref{subgroup}. Moreover, the semi-direct product in Formula \eqref{equation:isom_carquois_disjoint_D_cyclo} is well-defined since each $\widetilde{\Lambda}^{(j)}$ satisfies the condition $\widetilde{\Lambda}^{(j)}_i=\widetilde{\Lambda}^{(j)}_{\theta(i)}$ of
Proposition~\ref{prop-fixed-points}\ref{item:fixed_points_VV_cyclo} (see~\eqref{equation:VLambda_sdirect_product}).

\begin{remark} The reader may have noticed that the assumptions $d>1$ and $n_j(\beta)\neq 0$ (which do not reduce the generality as explained above) were not present in the preceding section for the type B in Theorem \ref{theorem:isom_carquois_disjoints_B} and Corollary \ref{corollary:isom_disjoint_quiver_cyclo_B}. Indeed those statements are more uniform in the sense that they are also valid as they are, even if some $n_j(\beta)$ are 0 or if $d=1$. In particular, for $d = 1$ we do \emph{not} necessarily have $W_\beta^\Lambda(\Gamma) = W_\beta^{\widetilde{\Lambda}}(\Gamma)$ (cf. Remark~\ref{remark:isom_Lambdatilde}).
\end{remark}

\begin{proof}
$\bullet$ Recall from Theorem \ref{theorem:isom_carquois_disjoints_B} that we have an isomorphism between $V_\beta(\Gamma)$ and the algebra $\mathrm{Mat}_{\binom{n}{\sizeglob 1, \dots, \sizeglob d}}\left(\bigotimes_{j=1}^d V_{\beta^{(j)}}(\Gamma^{(j)})\right)$. This isomorphism was obtained with the following two steps:
\[V_\beta(\Gamma)\simeq \mathrm{Mat}_{\binom{n}{\sizeglob 1, \dots, \sizeglob d}}\left(e(\mathfrak{t}^{\beta})V_\beta(\Gamma)e(\mathfrak{t}^{\beta})\right)\ \ \ \ \text{and}\ \ \ \ \bigotimes_{j=1}^d V_{\beta^{(j)}}(\Gamma^{(j)})\simeq e(\mathfrak{t}^{\beta})V_\beta(\Gamma)e(\mathfrak{t}^{\beta})\ .\]
For the first isomorphism, see~\textsection\ref{subsection:costless_isomorphism_B}, the construction of the idempotent $e(\mathfrak{t}^{\beta})$ does not involve~$\psi_0$, and neither does the construction of the matrix units (that is, the construction of the elements $\psi_{\mathfrak{t}}$ and $\phi_{\mathfrak{t}}$ given by Formulas \eqref{subequations:def_psit_phit}). So we deduce immediately how the automorphism $\iota$ of $V_\beta(\Gamma)$ behaves with respect to this isomorphism, namely we have that
\[V_\beta(\Gamma)^{\iota}\simeq \mathrm{Mat}_{\binom{n}{\sizeglob 1, \dots, \sizeglob d}}\left(e(\mathfrak{t}^{\beta})V_\beta(\Gamma)^{\iota}e(\mathfrak{t}^{\beta})\right)\ .\]
According to Formula \eqref{fixed-points-tensor} (that we can use since $n_j(\beta) \neq 0$), to prove~\eqref{equation:isom_carquois_disjoint_D} it remains only to show that 
\[e(\mathfrak{t}^{\beta})V_\beta(\Gamma)^{\iota}e(\mathfrak{t}^{\beta})\simeq \Bigl(\bigotimes_{j=1}^d V_{\beta^{(j)}}(\Gamma^{(j)})\Bigr)^{\iota^{\otimes}}\ .\]
So if we denote $\rho$ the isomorphic map from $\bigotimes_{j=1}^d V_{\beta^{(j)}}(\Gamma^{(j)})$ to $e(\mathfrak{t}^{\beta})V_\beta(\Gamma)e(\mathfrak{t}^{\beta})$, it remains to check that
\[\rho\circ\iota^{\otimes}=\iota\circ\rho\ .\]
This is immediately verified from Formulas \eqref{map1}--\eqref{map4} giving the map $\rho$ in the proof of Proposition~\ref{proposition:isom_tensor_product}. Moreover, the isomorphism~\eqref{equation:isom_carquois_disjoint_D} is graded since it is the restriction of a graded isomorphism (to a graded subalgebra). 

\medskip
$\bullet$ To prove~\eqref{equation:isom_carquois_disjoint_D_cyclo}, we start exactly as in the proof of Corollary~\ref{corollary:isom_disjoint_quiver_cyclo_B}, namely we repeat the calculations in the proof of Theorem \ref{theorem:isom_dec_klr_cyclo_block}. We can do so since $\H_\beta(\Gamma)$ is a subalgebra of $W_\beta(\Gamma)$ by Corollary \ref{corollary:klr_subalgebra_W}. 

Let $\vartheta$ denote the isomorphism in \eqref{equation:isom_carquois_disjoint_D} and let $\idealD_{\beta}^{\Lambda}$ denote the ideal of $W_\beta(\Gamma)$ giving the cyclotomic quotient $W^{\Lambda}_\beta(\Gamma)$. The proofs of Corollary~\ref{corollary:isom_disjoint_quiver_cyclo_B} and Theorem \ref{theorem:isom_dec_klr_cyclo_block} show that 
$$\vartheta(\idealD_{\beta}^{\Lambda})=\mathrm{Mat}_{\binom{n}{\sizeglob 1, \dots, \sizeglob d}}(\idealD_{\beta,\otimes}^{\Lambda})\,,$$
where $\idealD_{\beta,\otimes}^{\Lambda}$ is the ideal of $\Bigl(\bigotimes_{j=1}^dW_{\beta^{(j)}}(\Gamma^{(j)})\Bigr)\rtimes C_2^{d-1}$ generated by the elements
\begin{equation}\label{el-quot}
y_b^{\Lambda_{i_b}}e(\tuple i)\ ,
\end{equation}
where $\tuple i$ is of profile $\mathfrak{t}^{\beta}$, and $b$ is of the form $b=n_1+\dots+n_{j-1}+1$ for $j\in\{1,\dots,d\}$. Note that, as in the proof of Theorem \ref{theorem:isom_dec_klr_cyclo_block} we slightly abuse notations: if $\tuple i=(\tuple i^1,\dots,\tuple i^d)$ with $\tuple i^k\in \beta^{(k)}$, we identify $y_b^{\Lambda_{i_b}}e(\tuple i) \in W_\beta(\Gamma)$ with the element of $\bigotimes_{j=1}^dW_{\beta^{(j)}}(\Gamma^{(j)})$ which is $e(\tuple i^k)$ in the $k$-th factor with $k\neq j$ and $(y_1^{(j)})^{\Lambda_{(\tuple i^j)_1}}e(\tuple i^j)$ in the $j$-th factor (where $y_1^{(j)}$ denotes the generator $y_1$ of $W_{\beta^{(j)}}(\Gamma^{(j)})$).

Contrary to the type A and B, we need to show something more here to prove~\eqref{equation:isom_carquois_disjoint_D_cyclo}. In particular, we cannot consider the semi-direct product $\Bigl(\bigotimes_{j=1}^dW_{\beta^{(j)}}^{\Lambda^{(j)}}(\Gamma^{(j)})\Bigr)\rtimes C_2^{d-1}$  since the elements $\Lambda^{(j)}$ do not necessarily satisfy the stability condition of Proposition~\ref{prop-fixed-points}\ref{item:fixed_points_VV_cyclo}. Thus, let $\idealD_{\beta,\otimes}^{\widetilde{\Lambda}}$ be the ideal of $\Bigl(\bigotimes_{j=1}^dW_{\beta^{(j)}}(\Gamma^{(j)})\Bigr)\rtimes C_2^{d-1}$ generated by the elements
\begin{equation}\label{el-quot2} y_b^{\widetilde{\Lambda}_{i_b}}e(\tuple i)\ ,
\end{equation}
where $\tuple i \in \beta$ is of profile $\mathfrak{t}^{\beta}$, and $b$ is of the form $b=n_1+\dots+n_{j-1}+1$ for $j\in\{1,\dots,d\}$, and where $\widetilde{\Lambda}$ is defined in Theorem~\ref{theorem:isom_carquois_disjoints_D}. We will show that
$$\idealD_{\beta,\otimes}^{\Lambda}=\idealD_{\beta,\otimes}^{\widetilde{\Lambda}}\ .$$

First, since $\widetilde{\Lambda}_i\leq \Lambda_i$ for all $i\in I$, we have $\idealD_{\beta,\otimes}^{\Lambda}\subset \idealD_{\beta,\otimes}^{\widetilde{\Lambda}}$.
For the reverse inclusion, take an element $y_b^{\widetilde\Lambda_{i_b}} e(\tuple i)$ as in \eqref{el-quot}. If $\widetilde{\Lambda}_{i_b} = \Lambda_{i_b}$ then $y_b^{\widetilde\Lambda_{i_b}} e(\tuple i) \in \idealD_{\beta, \otimes}^\Lambda$, thus we assume that $\widetilde{\Lambda}_{i_b} = \Lambda_{\theta(i_b)}$. Let $\xi\in C_2^{d-1}$ such that the component of $\xi$ in position $j$ is $\pi$. Such an element exists since we assumed that $d>1$. Then, using Formulas \eqref{auto-pi} for the action of $\pi$ on $W_{\beta^{(j)}}(\Gamma^{(j)})$, we have, where $\tuple i' \in \beta$ of profile $\mathfrak{t}^{\beta}$ is such that $i'_b=\theta(i_b)$,
\[
\xi \cdot \left( y_b^{\widetilde{\Lambda}_{i_b}}e(\tuple i)\right)=(-y_b)^{\widetilde\Lambda_{i_b}}e(\tuple i')
= (-y_b)^{\Lambda_{\theta(i_b)}} e(\tuple i') = (-y_b)^{\Lambda_{i'_b}} e(\tuple i')
\, .
\]
 Since the action of $\xi$ is invertible, we thus deduce that $y_b^{\widetilde{\Lambda}_{i_b}} e(\tuple i) \in \idealD_{\beta,\otimes}^\Lambda$. Finally, we showed that all elements in \eqref{el-quot2} are in $\idealD_{\beta,\otimes}^{\Lambda}$, and thus $\idealD_{\beta,\otimes}^{\widetilde{\Lambda}}\subset \idealD_{\beta,\otimes}^{\Lambda}$. This concludes the proof.
\end{proof}

We define $W_n(\Gamma)\coloneqq \oplus_\delta W_\delta(\Gamma)$, where $\delta$ runs over all the orbits of $I^n$ under the action of~$D_n$, and similarly $W^\Lambda_n(\Gamma)= \oplus_\delta W^{\Lambda}_\delta(\Gamma)$. In the type D situation, the statements below are less clean that those of Corollary \ref{corollary:isom_dec_klr_n} or Corollary \ref{corollary:isom_disjoint_quiver_cyclo_B_n}. Nevertheless, it still explicitly reduces the study of $W_n(\Gamma)$ and $W^{\Lambda}_n(\Gamma)$ to the situation of a quiver with a single component.

For $(n_1,\dots,n_d)\in(\mathbb{Z}_{\geq0})^d$, we denote $l(n_1,\dots,n_d)$ the number of its non-zero components. Assume that $n\geq 1$ to avoid a trivial situation.
\begin{corollary}
\label{corollary:disjoint_quiver_isom_D}
We have (explicit) isomorphisms of graded algebras:
\begin{align*}
W_n(\Gamma) &\simeq \bigoplus_{\substack{n_1, \dots, n_d \geq 0 \\ n_1 + \dots + n_d = n}} \mathrm{Mat}_{\binom{n}{n_1, \dots, n_d}} \left(\Bigl(\bigotimes_{\substack{j=1\\n_j\neq 0}}^d W_{n_j}(\Gamma^{(j)})\Bigr)\rtimes C_2^{l(n_1,\dots,n_d)-1}\right)\ ,
\\
W^{\Lambda}_n(\Gamma) &\simeq \bigoplus_{\substack{n_1, \dots, n_d \geq 0 \\ n_1 + \dots + n_d = n}} \mathrm{Mat}_{\binom{n}{n_1, \dots, n_d}} \Bigl(\mathcal{W}(n_1,\dots,n_d)\Bigr)\ ,
\end{align*}
where:
\begin{itemize}
\item If $l(n_1,\dots,n_d)=1$ then $\mathcal{W}(n_1,\dots,n_d)\coloneqq W^{\Lambda^{(j)}}_{n_j}(\Gamma^{(j)})$ where $j$ is the component such that $n_j=n$.
\item If $l(n_1,\dots,n_d)>1$ then 
$$\mathcal{W}(n_1,\dots,n_d)\coloneqq \Bigl(\bigotimes_{\substack{j=1\\n_j\neq 0}}^d W^{\widetilde{\Lambda}^{(j)}}_{n_j}(\Gamma^{(j)})\Bigr)\rtimes C_2^{l(n_1,\dots,n_d)-1}\ .$$
\end{itemize}
\end{corollary}
\begin{proof}
We write $W_n(\Gamma) = \oplus_\beta W_\beta(\Gamma)$ and $W^{\Lambda}_n(\Gamma) = \oplus_\beta W^{\Lambda}_\beta(\Gamma)$, where $\beta$ runs over all the orbits of $I^n$ under the action of $B_n$. We note that if some $n_j(\beta)$ are equal to 0 then, as explained at the beginning of this subsection, we can remove the corresponding components of $\Gamma$ to obtain another quiver $\tilde{\Gamma}$ for which the assumptions of Theorem \ref{theorem:isom_carquois_disjoints_D} are satisfied. Then the proof is a repetition of the proof of Corollary \ref{corollary:isom_dec_klr_n}, using Theorem \ref{theorem:isom_carquois_disjoints_D} for each orbit $\beta$.
\end{proof}

\begin{remark}
As in Remarks~\ref{rem-cancellation_n} and~\ref{remark:cancellation_n_B}, we deduce that we can assume that $\Lambda$ is supported on all the components of $\Gamma$.
\end{remark}

\section{Morita equivalence for cyclotomic quotients of affine Hecke algebras of type~B and D}
\label{section:morita}

In this section, we will combine our previous results Corollaries~\ref{corollary:isom_disjoint_quiver_cyclo_B_n} and~\ref{corollary:disjoint_quiver_isom_D} with~\cite{poulain_dandecy-walker_B,poulain_dandecy-walker_D} to obtain Morita equivalences theorems for cyclotomic quotients of affine Hecke algebras of type~B and D. We emphasize that these Morita equivalences will be deduced from isomorphisms. As they combine the isomorphisms of~\cite{poulain_dandecy-walker_B,poulain_dandecy-walker_D} with those of the previous sections, these isomorphisms can be written down explicitly even though they are rather complicated.

Recall that $K$ is a field with characteristic different from two. Let $p,q\in K\setminus\{0\}$ such that $q^2\neq 1$.
As in Remark~\ref{remark:coherence_PAW}, for any $x \in K\setminus\{0\}$ we define the set 
\[
I_x \coloneqq \{ x^\epsilon q^{2l} : \epsilon \in \{\pm 1\}, l \in \mathbb{Z}\}\ .
\]
Then we take $d\geq 1$ and $x_1, \dots, x_d \in K^\times$ such that the sets $I^{(j)} \coloneqq I_{x_j}$ are pairwise disjoint, and we set
\[I\coloneqq \amalg_{j =1}^d I_{x_j}\ .\]

The quiver $\Gamma$ with involution that we will be considering in this section is the following:
\begin{itemize}
\item The vertex set of $\Gamma$ is $I$  as above.
\item  There is an arrow starting from $v$ and pointing to $q^2 v$ for all $v \in I$. These are all arrows.
\item The involution $\theta$ on $I$ is the scalar inversion $\theta(x)=x^{-1}$ for all $x\in I$. 
\end{itemize}

The partition $I = \amalg_{j = 1}^d I^{(j)}$ induces a decomposition of $\Gamma$ into full subquivers $\Gamma = \amalg_{j = 1}^d \Gamma^{(j)}$ as in Section~\ref{section:isomorphism B}, in particular each $\Gamma^{(j)}$ is stable under the scalar inversion $\theta$.
We also choose a finitely-supported family $\Lambda = (\Lambda_i)_{i \in I}$ of non-negative integers.
Finally, let $L$  be a free $\mathbb{Z}$-module of rank~$n$ with basis $\{\epsilon_i\}_{i=1,\dots,n}$:
\[L\coloneqq\bigoplus_{i=1}^n\mathbb{Z}\epsilon_i\ .\]

\subsection{Morita equivalence for cyclotomic quotients of affine Hecke algebras of type B}
\label{subsection:morita_B}
We set
\[\alpha_0\coloneqq2\epsilon_1\ \ \ \ \ \text{and}\ \ \ \ \ \alpha_i\coloneqq\epsilon_{i+1}-\epsilon_i\,,\ i=1,\dots,n-1\ .\]
For $n \geq 1$, the Weyl group $B_n$ of type B acts on $L$ by
\begin{align*}
r_0(\epsilon_1) &= -\epsilon_1,
\\
r_0(\epsilon_i)&=\epsilon_i \qquad \text{if } i>1,
\\
r_a(\epsilon_i)&=\epsilon_{r_a(i)},
\end{align*}
for $i = 1, \dots, n-1$ and $a = 1, \dots, n-1$.

We denote $q_0\coloneqq p$ and $q_i\coloneqq q$ for $i=1,\dots,n-1$. The affine Hecke algebra $\widehat{H}(B_n)$ is the unitary $K$-algebra generated by elements 
$$g_0,g_1,\dots,g_{n-1}\ \ \ \ \text{and}\ \ \ \ X^x,\ \ x\in L\ .$$
The defining relations are $X^0=1$, $X^xX^{x'}=X^{x+x'}$ for any $x,x'\in L$, and the characteristic equations for the generators $g_i$:
\begin{equation}\label{rel-H1}
g_i^2=(q_i-q_i^{-1})g_i+1 \ \qquad \text{for}\ i\in\{0,\dots,n-1\}\,,
\end{equation}
with the braid relations of type B
\begin{align}
&g_0g_1g_0g_1=g_1g_0g_1g_0 && \label{rel-H2}\\
&g_ig_{i+1}g_i=g_{i+1}g_ig_{i+1} && \text{for } i\in\{1,\dots,n-2\}, \label{rel-H3}
\\
&g_ig_j=g_jg_i && \text{for } i,j\in\{0,\dots,n-1\} \text{ such that } |i-j|>1,
\label{rel-H4}
\end{align}
together with
\[%\label{rel-Lu}
g_iX^x-X^{r_i(x)}g_i=(q_i-q_i^{-1})\frac{X^x-X^{r_i(x)}}{1-X^{-\alpha_i}}\ ,
\]
for any $x\in L$ and $i=0,1,\dots,n-1$. Note that the right-hand side is a well-defined element since there exists $k\in\mathbb{Z}$ such that $r_i(x)=x-k\alpha_i$. Note also that $\widehat{H}(B_0) = K$.

Let $X_i\coloneqq X^{\epsilon_i}$ for $i=1,\dots,n$. An equivalent presentation of the algebra $\widehat{H}(B_n)$ is with generators
\[
%\label{equation:generators_H(Bn)}
g_0,g_1,\dots,g_{n-1},X_1^{\pm1},\dots,X_n^{\pm1}\,,
\]
and defining relations \eqref{rel-H1}--\eqref{rel-H4} together with
\begin{align*}
&X_iX_j=X_jX_i  &&\text{for } i,j\in\{1,\dots,n\},%\label{rel-H5}\\
\\
&g_0X_1^{-1}g_0=X_1,% && \label{rel-H6}\\
\\
&g_iX_ig_i=X_{i+1} && \text{for } i\in\{1,\dots,n-1\},%}\,, \label{rel-H7}\\
\\
&g_iX_j=X_jg_i && \text{for } i\in\{0,\dots,n-1\} \text{ and } j\in\{1,\dots,n\} \text{ such that } j\neq i,i+1.%}\,. \label{rel-H8}
\end{align*}

\begin{definition}
The cyclotomic quotient $H^{\Lambda}(B_n)$ of type B associated with $\Lambda = (\Lambda_i)_{i \in I}$ is the quotient of the algebra $\widehat{H}(B_n)$ over the relation
\[\prod_{i\in I}(X_1-i)^{\Lambda_i}=0\ .\]
\end{definition}

Note that if $\Lambda_i = 0$ for all $i$ then
\[
H^\Lambda(B_n) = \begin{cases}
\{0\}, &\text{if } n \geq 1,
\\
K,&\text{if } n = 0.
\end{cases}
\]

We recall the main result of~\cite{poulain_dandecy-walker_B,poulain_dandecy-walker_D} concerning $H^\Lambda(B_n)$.

\begin{theorem}
\label{theorem:isom_B}
Let $\lambda, \gamma$ be as in Remarks~\ref{remark:coherence_PAW} and~\ref{remark:coherence_PAW_typeD} if $p^2 \neq 1$ and $p^2 = 1$ respectively. 
The algebras $H^\Lambda(B_n)$ and $V_n^\Lambda(\Gamma,\lambda,\gamma)$ are (explicitly) isomorphic.
\end{theorem}

\begin{remark}
Theorem~\ref{theorem:isom_B} is proven for $n \geq 1$, but is also trivially true for $n = 0$.
\end{remark}

We now state the first main application of the results of the preceding sections.

\begin{theorem}
\label{theorem:morita_B}
We have an (explicit) isomorphism of algebras:
\[H^{\Lambda}(B_n) \simeq \bigoplus_{\substack{n_1, \dots, n_d \geq 0 \\ n_1 + \dots + n_d = n}} \mathrm{Mat}_{\binom{n}{n_1, \dots, n_d}} \left(\bigotimes_{j=1}^d H^{\Lambda^{(j)}}(B_{n_j})\right).\]
In particular, $H^{\Lambda}(B_n)$ is Morita equivalent to  $\displaystyle\bigoplus_{\substack{n_1, \dots, n_d \geq 0 \\ n_1 + \dots + n_d = n}}\left(\bigotimes_{j=1}^d H^{\Lambda^{(j)}}(B_{n_j})\right)$.
\end{theorem}

\begin{proof}
Note that the statement is true if $n = 0$, thus we now assume $n \geq 1$.
Let us first assume that $p^2 \neq 1$.
Let $\lambda$ be the indicator function of $\{\pm p\} \cap I$ and $(\gamma_i)_{i \in I}$ be given by $\gamma_i = \begin{cases}
1, &\text{if } \theta(i) = i,
\\
0, &\text{otherwise},
\end{cases}$ as in Remark~\ref{remark:coherence_PAW}. By Theorem~\ref{theorem:isom_B}, we have an isomorphism $H^\Lambda(B_n) \simeq V_n^\Lambda(\Gamma,\lambda,\gamma)$. For any $j \in \{1, \dots, d\}$, the restrictions $\lambda^{(j)}$  and $\gamma^{(j)}$ of $\lambda$ and $\gamma$ respectively to $I^{(j)}$ satisfy, by Corollary~\ref{corollary:isom_disjoint_quiver_cyclo_B_n},
%are still of the above form with respect to the quiver $\Gamma^{(j)}$. Hence, the isomorphism
\[
V_n^\Lambda(\Gamma,\lambda,\gamma) \simeq \bigoplus_{\substack{n_1, \dots, n_d \geq 0 \\ n_1 + \dots + n_d = n}} \mathrm{Mat}_{\binom{n}{n_1, \dots, n_d}} \left(\bigotimes_{j=1}^d V_{n_j}^{\Lambda^{(j)}}(\Gamma^{(j)},\lambda^{(j)}, \gamma^{(j)})\right)\ .
\]
Since $\lambda^{(j)}$ and $\gamma^{(j)}$ are still of the above form with respect to the quiver $\Gamma^{(j)}$, by Theorem~\ref{theorem:isom_B} we have $V_{n_j}^{\Lambda^{(j)}}(\Gamma^{(j)}, \lambda^{(j)}, \gamma^{(j)}) \simeq H^{\Lambda^{(j)}}(B_{n_j})$  for any $n_j$.
We thus deduce the isomorphism of the theorem. We deduce the statement of Morita equivalence since $\mathrm{Mat}_N(A)$ and $A$ are Morita equivalent for any algebra $A$ and $N \in \mathbb{N}^*$. The case $p^2 = 1$ is similar, still by Theorem~\ref{theorem:isom_B}.
\end{proof}

We obtain the following corollary.

\begin{corollary}
\label{corollary:morita_B}
To study an arbitrary cyclotomic quotient of the affine Hecke algebra $\widehat{H}(B_n)$, it is enough to consider cyclotomic quotients given by a relation
\[\prod_{\substack{\epsilon \in \{\pm 1\} \\ l\in \mathbb{Z}}}(X_1-x^\epsilon q^{2l})^{m_{\epsilon, l}}=0\ ,\]
for any finitely-supported family of non-negative integers $(m_{\epsilon, l})_{\epsilon \in \{\pm 1\}, l\in\mathbb{Z}}$, where $x \in K^\times$ satisfies one of the following four cases:
\[\mathbf{(a)}\ x=1\ \qquad \mathbf{(b)}\ x=q\ \qquad \mathbf{(c)}\ x=p\ \qquad \mathbf{(d)}\ x\notin\pm q^{\mathbb{Z}}\cup\pm p^{\pm1}q^{2\mathbb{Z}}\ .\]
\end{corollary}

\begin{proof}
We sketch a proof, in the same spirit as in the introduction of~\cite{poulain_dandecy-walker_B}.
By Theorem~\ref{theorem:morita_B}, it is clear that it suffices to consider cyclotomic quotients given by a relation
\[
\prod_{i \in I_x} (X_1 - i)^{\Lambda_i} = 0,
\]
where $I_x = \{x^\epsilon q^{2l} : \epsilon \in \{\pm 1\}, l \in \mathbb{Z}\}$ with $x \in K^\times$ and $\Lambda = (\Lambda_i)_{i \in I_x}$ is a finitely-supported family of non-negative integers. By Theorem~\ref{theorem:isom_B} and Remark~\ref{remark:coherence_PAW}, this cyclotomic quotient is determined by:
\begin{itemize}
\item the quiver $\Gamma$ with vertex set $I_x$, arrows $v \to q^2 v$ for all $v \in I_x$ and involution $\theta : v \mapsto v^{-1}$ on $I_x$;
\item the set $\{\pm p\} \cap I_x$.
\end{itemize}

A first distinction arises when looking at the number of connected components of $\Gamma$. It has exactly one (respectively two) connected component(s) when $x^2 \in q^{2\mathbb{Z}}$ (resp. $x^2 \notin q^{2\mathbb{Z}}$).

The first case, $x^2 \in q^{2\mathbb{Z}}$, is equivalent to $x \in \pm q^\mathbb{Z}$. We can switch between $x$ and $-x$ by the variable change $X_i \leftarrow -X_i$ for all $i \in \{1, \dots, n\}$, replacing $I_x$ by $-I_x = I_{-x}$ and $\Lambda = (\Lambda_i)_{i \in I_x}$ by $\Lambda' = (\Lambda'_i)_{i \in I_{-x}}$ given by $\Lambda'_i \coloneqq \Lambda_{-i}$ for all $i \in I_{-x}$. Thus, it suffices to consider $x \in q^\mathbb{Z}$, but now a simple shift of $\Lambda$ (that is, setting $\Lambda'_i=\Lambda_{iq^{2N}}$ for appropriate $N$) shows that it suffices to consider the cases $x = 1$ (this is case $\mathbf{(a)}$) or $x = q$ (this is case $\mathbf{(b)}$), according to the parity of the power of $q$.

We now consider the case $x^2 \notin q^{2\mathbb{Z}}$, that is, $x \notin \pm q^\mathbb{Z}$. If $\{\pm p\} \cap I_x = \emptyset$, then $x \notin \pm q^\mathbb{Z} \cup \pm p^{\pm 1} q^{2\mathbb{Z}}$, and all these choices of $x$ lead to isomorphic algebras since moreover $\theta$ has no fixed points (if $x^{\pm 1} q^{2k}$ is fixed by $\theta$ then $x^2 \in q^{4\mathbb{Z}}$ thus $x \in \pm q^{2\mathbb{Z}}$). This is case $\mathbf{(d)}$. Now if $\{\pm p\} \cap I_x \neq \emptyset$, using the variable change $X_i \leftarrow -X_i$ for all $i \in \{1, \dots, n\}$ we can always assume that $p \in I_x$, that is, $x \in p^{\pm 1} q^{2 \mathbb{Z}}$. It suffices in fact to consider $x \in p q^{2\mathbb{Z}}$, since the variable change $g_0 \leftarrow -g_0$ exchanges $p$ and $p^{-1}$. This case reduces to $x = p$ by shifting $\Lambda$ as above, and this is case $\mathbf{(c)}$.
\end{proof}

\begin{remark}\label{rem-final}
We make additional final remarks on the four cases $\mathbf{(a)}$--$\mathbf{(d)}$ to be considered.
\begin{itemize}
\item Cases $\mathbf{(a)}$ and $\mathbf{(b)}$ correspond to a quiver with a single connected component (an infinite oriented line or a finite oriented polygon depending on whether $q$ is a root of unity or not). This quiver is stable by the involution $\theta$, and then Case $\mathbf{(a)}$ corresponds to $\theta$ having a fixed point, while Case $\mathbf{(b)}$ generically corresponds to the situation where there is no fixed point. This latter situation cannot occur if the number of vertices is finite and odd, that is, Case $\mathbf{(b)}$ is not present (or more precisely, is not necessary since it is equivalent to Case $\mathbf{(a)}$) when $q^2$ is an odd root of unity.
\item Cases $\mathbf{(c)}$ and $\mathbf{(d)}$ (generically) correspond to a quiver with two identical connected components (two infinite oriented lines or two finite oriented polygons depending on whether $q$ is a root of unity or not), which are exchanged by the involution $\theta$. Then Case $\mathbf{(c)}$ corresponds to the situation where one of the special values $\pm p^{-1}$ is present, while Case $\mathbf{(d)}$ corresponds to the situation where no such values occur. We see that Case $\mathbf{(c)}$ is not necessary (more precisely, it reduces to one of Cases $\mathbf{(a)}$ or $\mathbf{(b)}$) when $p^2$ is a power of $q^2$.
\item To summarise, there are at least two cases to consider in general: $\mathbf{(a)}$ and $\mathbf{(d)}$, while the additional two cases $\mathbf{(b)}$ and $\mathbf{(c)}$ are to be considered or not depending on $p$ and $q$.
\end{itemize}
\end{remark}

\subsection{Morita equivalence for cyclotomic quotients of affine Hecke algebras of type D}
Let $n\geq 2$. We set
\[\alpha'_0=\epsilon_1+\epsilon_2\ \ \ \ \ \text{and}\ \ \ \ \ \alpha'_i=\epsilon_{i+1}-\epsilon_i\,,\ i=1,\dots,n-1\ .\]
The Weyl group $D_n$ of type D acts on $L$ by
\begin{align*}
s_0(\epsilon_1) &= -\epsilon_2,\\
s_0(\epsilon_2)&=-\epsilon_1,
\\
s_0(\epsilon_i)&=\epsilon_i, \qquad \text{if } i > 2,
\\
s_a(\epsilon_i)&=\epsilon_{r_a(i)},
\end{align*}
for $i = 1,\dots,n-1$ and $a = 1,\dots,n-1$.

The affine Hecke algebra $\widehat{H}(D_n)$ is the unitary $K$-algebra generated by elements 
\[
\{g_i\}_{1 \leq i \leq n-1} \cup \{G_0\}  \cup \{X^x\}_{x\in L}.
\]
The defining relations are $X^0=1$, $X^xX^{x'}=X^{x+x'}$ for any $x,x'\in L$, and
the characteristic equations for the generators $g_i$ and $G_0$:
\begin{equation}\label{rel-HD1}
\begin{aligned}
g_i^2&=(q-q^{-1})g_i+1 & \text{for } i\in\{1,\dots,n-1\},
\\
G_0^2 &= (q - q^{-1})G_0 + 1,
\end{aligned}
\end{equation}
with the braid relations of type D
\begin{align}
\label{rel-HD2}
& G_0g_2G_0=g_2G_0g_2,%\,, &&   \\[0.2em]
\\
\label{rel-HD3}
& G_0g_i=g_{i}G_0 &&  \text{for } i\in\{1,\dots,n-1\}\setminus\{2\},%}\,,\\[0.2em]
\\
\label{rel-HD4}
& g_ig_{i+1}g_i=g_{i+1}g_ig_{i+1} &&  \text{for } i\in\{1,\dots,n-2\},%}\,,\\[0.2em]
\\
\label{rel-HD5}
& g_ig_j=g_jg_i && \text{for } i,j\in\{1,\dots,n-1\} \text{ such that } |i-j|>1,%}\,.
\end{align}
together with
\begin{align*}
g_iX^x-X^{s_i(x)}g_i&=(q-q^{-1})\frac{X^x-X^{s_i(x)}}{1-X^{-\alpha'_i}},
\\
G_0 X^x-X^{s_0(x)}G_0&=(q-q^{-1})\frac{X^x-X^{s_0(x)}}{1-X^{-\alpha'_0}},
\end{align*}%\end{equation}
for any $x\in L$ and $i=1,\dots,n-1$. Note that the right-hand sides are well-defined elements since for any $i \in \{0, \dots, n-1\}$ there exists $k\in\mathbb{Z}$ such that $s_i(x)=x-k\alpha'_i$.

An equivalent presentation of the algebra $\widehat{H}(D_n)$ is with generators (where again $X_i\coloneqq X^{\epsilon_i}$)
\[
\{g_i\}_{1 \leq i \leq n-1} \cup \{G_0\} \cup \{X_i^{\pm1}\}_{1 \leq i \leq n}\,,
\]
and defining relations \eqref{rel-HD1}--\eqref{rel-HD5} together with
\begin{align*}
& X_iX_j=X_jX_i && \text{for } i,j\in\{1,\dots,n\},
\\
& G_0X_1^{-1}G_0=X_2,% && \\[0.2em]
\\
& G_0X_i=X_iG_0 && \text{for } i\in\{3,\dots,n-1\},%}\,,\\[0.2em]
\\
& g_iX_ig_i=X_{i+1} && \text{for } i\in\{1,\dots,n-1\},%}\,,\\[0.2em]
\\
& g_iX_j=X_jg_i && \text{for } i\in\{1,\dots,n-1\} \text{ and } j\in\{1,\dots,n\}  \text{ such that } j\neq i,i+1.
\end{align*}

By convention, we set that $\widehat{H}(D_n)$ coincides with the usual affine Hecke algebra of type $A_n$ if $n \in \{0, 1\}$, that is, we have $\widehat{H}(D_0) = K$ and $\widehat{H}(D_1) = K[X_1^{\pm 1}]$.

\begin{definition}
The cyclotomic quotient $H^{\Lambda}(D_n)$ of type D associated with $\Lambda = (\Lambda_i)_{i \in I}$  is the quotient of the algebra $\widehat{H}(D_n)$ over the relation
\[\prod_{i\in I}(X_1-i)^{\Lambda_i}=0\ .\]
\end{definition}

Note that if $\Lambda_i = 0$ for all $i$ then
\[
H^\Lambda(D_n) = \begin{cases}
\{0\}, &\text{if } n \geq 1,
\\
K, &\text{if } n = 0.
\end{cases}
\]

We recall the main result of~\cite{poulain_dandecy-walker_D} concerning $H^\Lambda(D_n)$. Recall that the quiver $\Gamma$ was defined at the beginning of Section~\ref{section:morita}.

\begin{theorem}
\label{theorem:isom_D}
The algebras $H^\Lambda(D_n)$ and $W_n^\Lambda(\Gamma)$ are (explicitly) isomorphic.
\end{theorem}

\begin{remark}
Theorem~\ref{theorem:isom_D} is proven for $n \geq 2$, but it is immediate with our conventions that it remains true for $n \in \{0, 1\}$.
\end{remark}

\paragraph{Expression as a semi-direct product.} We assume here that $n \geq 1$.
Assuming $p^2 = 1$, we now can see $\widehat{H}(D_n)$ as a subalgebra of $\widehat{H}(B_n)$. Namely, we have an inclusion (see, for instance, \cite[\textsection 2.3]{poulain_dandecy-walker_D}) $\widehat{H}(D_n) \subseteq \widehat{H}(B_n)$, given on the generators by
\begin{align*}
G_0 &\mapsto g_0 g_1 g_0, &  g_i &\mapsto g_i, & X_j^{\pm 1} &\mapsto X_j^{\pm 1},
\end{align*}
for any $i \in \{1, \dots, n-1\}$ and $j \in \{1, \dots, n\}$. Another way to see $\widehat{H}(D_n)$ as a subalgebra of $\widehat{H}(B_n)$ is to write $\widehat{H}(D_n)$ as the subalgebra of fixed points of $\widehat{H}(B_n)$ under the involution $\eta$ given by
\begin{align*}
g_0 &\mapsto -g_0, & g_i &\mapsto g_i, & X_j^{\pm 1} &\mapsto X_j^{\pm 1},
\end{align*}
for each $i \in \{1, \dots, n-1\}$ and $j \in \{1, \dots, n\}$ (note that since $p^2 = 1$ the defining relation for the generator $g_0$ is $g_0^2 = 1$).  In particular, as in~\textsection\ref{subsection:def_W}  we have a vector space decomposition $\widehat{H}(B_n) = \widehat{H}(D_n) \oplus \widehat{H}(D_n) g_0$ and thus an isomorphism of algebras
\[
\widehat{H}(B_n) \simeq \widehat{H}(D_n) \rtimes C_2.
\]
Note that the action of the generator of $C_2$ on the generating set of $\widehat{H}(D_n)$ is given by
\begin{align*}
G_0 &\mapsto g_1, & g_1 &\mapsto G_0, & g_i &\mapsto g_i,
\\
X_1 &\mapsto X_1^{- 1}, & X_1^{-1} &\mapsto X_1, & X_j^{\pm 1} &\mapsto X_j^{\pm 1},
\end{align*}
for all $i \in \{2, \dots, n-1\}$ and $j \in \{2, \dots, n\}$. 

\medskip
The involution~$\eta$ on $\widehat{H}(B_n)$ is compatible with the cyclotomic quotient $H^\Lambda(B_n)$.  Now if $\Lambda$ satisfies the stability condition of Proposition~\ref{prop-fixed-points}\ref{item:fixed_points_VV_cyclo} (which is here $\Lambda_{i^{-1}} = \Lambda_i$ for all $i \in I$), the previous action of $C_2$ on $\widehat{H}(D_n)$ is compatible with the cyclotomic quotient $H^\Lambda(D_n)$ and as above we have
\[
H^\Lambda(B_n) \simeq H^\Lambda(D_n) \rtimes C_2.
\]

\paragraph{Morita equivalence theorem.}

Let $n_1, \dots, n_d \geq 1$. If $\Lambda$ satisfies $\Lambda_{i^{-1}} = \Lambda_i$ for all $i \in I$, the previous action of $C_2$ on $H^\Lambda(D_n)$ extends to a (diagonal) action of $C_2^d$ on $\otimes_{j = 1}^d H^{\Lambda^{(j)}}(D_{n_j})$. As in~\textsection\ref{subsection:disjoint_quiver_isom_W}, we restrict this action to the subgroup $C_2^{d-1}$ of even elements given in~\eqref{subgroup}. Recall also the definition of $\widetilde{\Lambda} = (\widetilde{\Lambda}_i)_{i \in I}$ given in Theorem~\ref{theorem:isom_carquois_disjoints_D}.

We now state the second main application of the paper. As in Corollary~\ref{corollary:disjoint_quiver_isom_D}, for any $(n_1, \dots, n_d)  \in (\mathbb{Z}_{\geq 0})^d$ we denote by $l(n_1, \dots, n_d)$ the number of its non-zero components.

\begin{theorem}
\label{theorem:morita_D}
We have an (explicit) isomorphism of algebras:
\[
H^{\Lambda}(D_n) \simeq \bigoplus_{\substack{n_1, \dots, n_d \geq 0 \\ n_1 + \dots + n_d = n}} \mathrm{Mat}_{\binom{n}{n_1, \dots, n_d}} \Bigl(\mathcal{H}(n_1,\dots,n_d)\Bigr)\ ,
\]
where:
\begin{itemize}
\item If $l(n_1,\dots,n_d)=1$ then $\mathcal{H}(n_1,\dots,n_d)\coloneqq H^{\Lambda^{(j)}}(D_{n_j})$ where $j$ is the component such that $n_j=n$.
\item If $l(n_1,\dots,n_d)>1$ then 
$$\mathcal{H}(n_1,\dots,n_d)\coloneqq \Bigl(\bigotimes_{\substack{j=1\\n_j\neq 0}}^d H^{\widetilde{\Lambda}^{(j)}}(D_{n_j})\Bigr)\rtimes C_2^{l(n_1,\dots,n_d)-1}\ .$$
\end{itemize}
In particular, $H^{\Lambda}(D_n)$ is Morita equivalent to $\displaystyle\bigoplus_{\substack{n_1, \dots, n_d \geq 0 \\ n_1 + \dots + n_d = n}}  \mathcal{H}(n_1,\dots,n_d)$.
\end{theorem}

\begin{proof}
We argue as in  the proof of Theorem~\ref{theorem:morita_B}, using Corollary~\ref{corollary:disjoint_quiver_isom_D} and Theorem~\ref{theorem:isom_D}. Note that the isomorphism of~\cite{poulain_dandecy-walker_D} is compatible with the semi-direct product since the involution $\iota$ (respectively, the element $\psi_0$) of $V^\Lambda_n(\Gamma)$ is sent to the involution $\eta$ (resp., the element $g_0$) of $H^\Lambda(B_n)$ by the isomorphism of \textit{loc. cit.}
\end{proof}

We obtain the following corollary. We note that the situation is a little bit more intricate than for type B because of the presence of semidirect products with products of groups $C_2$. So below, it is implicit that it is enough to consider some special cyclotomic quotients, up to the application of standard Clifford theory to deal with the semidirect products.
\begin{corollary}
\label{corollary:morita_D}
To study an arbitrary cyclotomic quotient of the affine Hecke algebra $\widehat{H}(D_n)$, it is enough to consider cyclotomic quotients given by a relation
\[\prod_{\substack{\epsilon\in\{\pm 1\} \\ l\in \mathbb{Z}}}(X_1-x^\epsilon q^{2l})^{m_{\epsilon, l}}=0\ ,\]
for any finitely-supported family of non-negative integers $(m_{\epsilon,l})_{\epsilon\in\{\pm 1\}, l\in\mathbb{Z}}$, where $x \in K^\times$ satisfies one of the following three cases:
\[\mathbf{(a)}\ x=1\ \qquad \mathbf{(b)}\ x=q\ \qquad \mathbf{(c)}\ x\notin \pm q^{\mathbb{Z}}\ .\]
\end{corollary}

\begin{proof}
We sketch a proof, in the same spirit as in the introduction of~\cite{poulain_dandecy-walker_D}.
We deduce from Theorem~\ref{theorem:morita_D} that it suffices to study the cyclotomic quotients of $\widehat{H}(D_n)$ given by a relation
\[
\prod_{i \in I_x} (X_1 - i)^{\Lambda_i},
\]
where $I_x$ and $\Lambda$ are as in the proof of Corollary~\ref{corollary:morita_B}. By Theorem~\ref{theorem:isom_D}, this cyclotomic quotient is only determined by the quiver $\Gamma$ and its involution $\theta$ as defined in the proof of Corollary~\ref{corollary:morita_B}. In particular, looking at the number of connected components of $\Gamma$ we still have the two cases $x \in \pm q^\mathbb{Z}$ (which give cases $\mathbf{(a)}$ and $\mathbf{(b)}$) and $x \notin \pm q^\mathbb{Z}$ (which is case $\mathbf{(c)}$). In the latter case all the choices of $x$ lead to isomorphic algebras since $\theta$ has no fixed points.
\end{proof}

\begin{remark}\label{rem-final_D}
We make an additional final remark on the three cases $\mathbf{(a)}$--$\mathbf{(c)}$ to be considered, similarly to Remark \ref{rem-final}. Cases $\mathbf{(a)}$ and $\mathbf{(b)}$ correspond to a quiver with a single connected component (an infinite oriented line or a finite oriented polygon depending on whether $q$ is a root of unity or not), while $\mathbf{(c)}$ corresponds to a quiver with two identical connected components exchanged by the involution $\theta$. Case $\mathbf{(a)}$ corresponds to $\theta$ having a fixed point, while Case $\mathbf{(b)}$ generically corresponds to the situation where there is no fixed point. As before, when $q^2$ is an odd root of unity, Case $\mathbf{(b)}$ is not necessary since it is equivalent to Case $\mathbf{(a)}$.
\end{remark}

\appendix

\section{Polynomial realisation}
\label{appendix_section:polynomial_realisation}

We prove here Lemma~\ref{lemma:pbw_action_well_defined}. In this appendix,  for any $f \in K[x, \beta]$ we also systematically write $f$ for the the element of $\mathrm{End}_K(K[x, \beta])$ given by left multiplication and we use concatenation to denote the composition inside  $\mathrm{End}_K(K[x, \beta])$. In particular, for any $w \in B_n$ and $f \in K[x]$ we have $wf = (\actpol{w}{f}) w$ inside $\mathrm{End}_K(K[x, \beta])$.

 We now define some elements of $\mathrm{End}_K(K[x, \beta])$ by
\begin{equation}
\label{equation:polynomial_realisation}
\begin{aligned}
\polrep(e(\tuple{i})) &= \multone_{\tuple{i}},
\\
\polrep(y_a e(\tuple i)) &=  \multx_a \multone_{\tuple i},
\\
\polrep(\psi_b e(\tuple{i})) &= \Bigl(\delta_{i_b, i_{b+1}} (\multx_b - \multx_{b+1})^{-1}(r_b - 1) 
+
P_{i_b, i_{b+1}}(\multx_{b+1}, \multx_b) r_b \Bigr) \multone_{\tuple i},
\\
\polrep(\psi_0 e(\tuple{i})) &= \left( \gamma_{i_1}\multx_1^{-1}(1 - r_0) + \alpha_{i_1}(\multx_1)r_0\right)\multone_{\tuple i},
\end{aligned}
\end{equation}
for any $a \in \{1, \dots, n\}$ and $b \in \{1, \dots, n-1\}$, and extend these formulas to $\polrep(X)$ for $X \in \{y_1, \dots, y_n, \psi_0, \dots, \psi_{n-1}\}$ by $\polrep(X) = \sum_{\tuple i \in \beta} \polrep(X e(\tuple i))$.

We will prove that $\polrep$ extends to an algebra homomorphism $\polrep : V_\beta(\Gamma,\lambda,\gamma) \to \mathrm{End}_K(K[x, \beta])$, which will imply Lemma~\ref{lemma:pbw_action_well_defined}. Indeed, the map $\polrep$ is the homomorphism associated with the action defined in~\textsection\ref{subsection:basis_theorem}. To prove that $\polrep$ extends to an algebra homomorphism, we check the defining relations of $V_\beta(\Gamma, \lambda,\gamma)$.
Recall that $P_{i, j} = 0$ when $i = j$ so that
\[
\polrep(\psi_b e(\tuple i)) = \begin{cases}
(\multx_b - \multx_{b+1})^{-1} (r_b - 1)  \multone_{\tuple i}, &\text{if } i_b = i_{b+1},
\\
P_{i_b, i_{b+1}}(\multx_{b+1}, \multx_b) r_b \multone_{\tuple i}, &\text{otherwise}.
\end{cases}
\]
%Note that $(x_a - x_{a+1})^{-1}(r_a -1) = \partial_a$ is the divided difference operator, hence stabilises $K[x_1, \dots, x_n]$.
Moreover, by~\eqref{equation:condition_gamma} and~\eqref{equation:alpha_0}  we have
\[
\polrep(\psi_0 e(\tuple i)) = \begin{cases}\alpha_{i_1}(\multx_1)r_0  \multone_{\tuple i},
&\text{if } \gamma_{i_1} = 0,
\\
\gamma_{i_1} \multx_1^{-1}(1 - r_0)\multone_{\tuple i}, &\text{otherwise}.
\end{cases}
\]

The relations that do not  involve $\psi_0$ are satisfied since the action is the same as in~\cite[Proposition 3.12]{rouquier_2}.
Relations~\eqref{relation:psi0_e(i)}, \eqref{relation:psi0_psib} and~\eqref{relation:psi0_yj} are immediate. 

To simplify the notation, for any $v \in V_\beta(\Gamma,\lambda,\gamma)$ we  also write $\prep v$ instead of $\polrep (v)$. Note that the composition operation in $\mathrm{End}_K(K[x, \beta])$ is denoted as a simple multiplication. For example, $\prep \psi_0\multx_1$ means composition of the multiplication by $x_1$ with the operator $\phi(\psi_0)$. Concerning~\eqref{relation:psi0_y1}, we have
\begin{align*}
(\prep \psi_0\prep y_1 + \prep y_1 \prep\psi_0)\prep{e(\tuple i)}
&=
\prep \psi_0\multx_1 \mathbf{1}_{\tuple{i}} + \multx_1 \left(\gamma_{i_1}\multx_1^{-1} (1 - r_0) + \alpha_{i_1}(\multx_1)r_0 \right)\multone_{\tuple i}
\\
&=
\Bigl[\left(\gamma_{i_1} \multx_1^{-1}(1 - r_0) \multx_1 +  \alpha_{i_1}(\multx_1) r_0 \multx_1 \right) + \left(\gamma_{i_1}(1 - r_0) + \multx_1 \alpha_{i_1}(\multx_1)  r_0\right)\Bigr]\multone_{\tuple i}
\\
&=
\Bigl[\gamma_{i_1} (1 + r_0) - \multx_1  \alpha_{i_1}(\multx_1) r_0 + \gamma_{i_1}(1 - r_0) + \multx_1 \alpha_{i_1}(\multx_1)  r_0\Bigr]\multone_{\tuple i}
\\
&=
2\gamma_{i_1} \multone_{\tuple i}=\polrep\bigl(2\gamma_{i_1} e(\tuple i)\bigr)\ .
\end{align*}
%Since $\gamma_{i_1} = 0$ if $\tuple i \neq r_0 \cdot \tuple i$ by~\eqref{subequations:conditions_gamma}, we obtain
%\begin{align*}
%(\prep \psi_0\prep y_1 + \prep y_1 \prep\psi_0)\prep{e(\tuple i)}
%&=
%2\gamma_{i_1} \multone_{\tuple i}
%\\
%&=
%\polrep\bigl(2\gamma_{i_1} e(\tuple i)\bigr)
%\\
%&=
%\polrep\bigl(\psi_0 y_1 + y_1 \psi_0)e(\tuple i)\bigr),
%\end{align*}
%thus we conclude.

For~\eqref{relation:psi0square}, if $\gamma_i = 0$ then $\gamma_{\theta(i)} = 0$ by~\eqref{equation:gamma_thetainv} and  we have, noting that $\multone_{\tuple j} r_0 = r_0 \multone_{r_0 \cdot \tuple j}$ inside $\mathrm{End}_K(K[x, \beta])$,
\begin{align*}
{\prep\psi_0}^2 \prep{e(\tuple{i})}
&=
\prep\psi_0 \alpha_{i_1}(\multx_1)r_0\multone_{\tuple i}
\\
&=
\alpha_{\theta(i_1)}(\multx_1)r_0 \alpha_{i_1}(\multx_1)r_0 \multone_{\tuple i}
\\
&=
\alpha_{\theta(i_1)}(\multx_1) \alpha_{i_1}(-\multx_1)\multone_{\tuple i}
\\
&=
(-1)^{\lambda_{\theta(i_1)}} \multx_1^{d(i_1)} \multone_{\tuple i}
\\
&=
\polrep\left((-1)^{\lambda_{\theta(i_1)}} y_1^{d(i_1)} e(\tuple i)\right)
%\\
%&=
%\polrep\bigl(\psi_0^2e(\tuple i)\bigr)
\,,
\end{align*}
by~\eqref{equation:equation_alpha}, and if $\gamma_{i_1} \neq 0$ then $\gamma_{\theta(i_1)}\neq 0$ and we have
\begin{align*}
{\prep\psi_0}^2 \prep{e(\tuple i)}
&=
\prep\psi_0  \gamma_{i_1} \multx_1^{-1}(1 - r_0) \multone_{\tuple i}
\\
&=
\gamma_{\theta(i_1)} \gamma_{i_1} \left(\multx_1^{-1}(1 - r_0)\right)^2 \multone_{\tuple i}
\\
&=
0.
\end{align*}
It remains to check~\eqref{relation:braid4}. 
As in~\eqref{equation:example_short}, we write $i_1 i_2$ and even $12$ instead of $\tuple{i}$, and $\bar{a}$ instead of $\theta(i_a)$. We have, using~\eqref{relation:psi0_e(i)},
\begin{subequations}
\label{subequations:action_braid4}
\begin{align}
(\prep\psi_0 \prep\psi_1)^2 \prep{e(12)}
&=
(\prep\psi_0\multone_{1\bar{2}}) (\prep\psi_1 \multone_{\bar{2}1} ) (\prep\psi_0  \multone_{21} )(\prep\psi_1  \multone_{12}),
\\
(\prep\psi_1 \prep\psi_0)^2 \prep{e(12)} &= (\prep\psi_1\multone_{\bar{2}\bar{1}})(\prep\psi_0\multone_{2\bar{1}})(\prep\psi_1 \multone_{\bar{1}2})(\prep\psi_0  \multone_{12}).
\end{align}
\end{subequations}

\subsection{Case \texorpdfstring{$\gamma_{i_1} = 0 = \gamma_{i_2}$}{gammai1=0=gammai2}}
\label{subsection:gammai1=0=gammai2}

First, recall that by~\eqref{equation:gamma_thetainv} we know that if $\gamma_{i_1} = 0$ and $\theta(i_1) = i_2$ then $\gamma_{i_2} = 0$. Thus, we want to prove that
\begin{equation}
\label{equation:braid4_gamma0=0=gamma1}
\left((\prep\psi_0 \prep\psi_1)^2 - (\prep\psi_1 \prep\psi_0)^2 \right) \multone_{\tuple{i}} = \begin{cases}
 (-1)^{\lambda_{\theta(i_1)}} \frac{(-\prep y_1)^{d(i_1)} - {\prep y_2}^{d(i_1)}}{\prep y_1 + \prep y_2} \prep \psi_1 \multone_{\tuple{i}},
 & \text{if } \theta(i_1) = i_2,
 \\
 0, &\text{otherwise}.
 \end{cases}
\end{equation}

Since $\gamma_{i_1} = \gamma_{\theta(i_1)} = \gamma_{i_2} = \gamma_{\theta(i_2)} = 0$, for any $a, b \in \{1, 2, \bar{1}, \bar{2}\}$ the element $\psi_0$ acts on $\multone_{ab}$ as $\alpha_a(\multx_1) r_0$.

Assume that $\theta(i_1) = i_1$ and $\theta(i_2) = i_2$. By~\eqref{equation:condition_strong} we have $d(i_1) = d(i_2) = 0$, thus~\eqref{equation:braid4_gamma0=0=gamma1} becomes
\begin{equation}
\label{equation:braid4_gamma0=0=gamma1_zerocase}
\left((\prep\psi_0 \prep\psi_1)^2 - (\prep\psi_1 \prep\psi_0)^2 \right) \multone_{\tuple{i}} = 0.
\end{equation}
Since $d(i_1) = d(i_2) =  0$, by~\eqref{equation:equation_alpha} we can assume  $\alpha_{i_1}(y) =\alpha_{i_2}(y) = 1$, thus $\psi_0$ acts on $\multone_{a b}$ as $r_0$ for any $a, b$. Hence, the same calculation as in~\cite[\S3.1]{poulain_dandecy-walker_D} proves that~\eqref{equation:braid4_gamma0=0=gamma1_zerocase} is satisfied.
In the opposite case, if $\theta(i_1) \neq i_1$ and $\theta(i_2) \neq i_2$ we know by the proof of~\cite[Proposition 7.4]{varagnolo-vasserot_canonical} that~\eqref{equation:braid4_gamma0=0=gamma1}  holds. 

Thus, we now assume that $\theta(i_1) = i_1$ and $\theta(i_2) \neq i_2$, in particular $i_1 \neq i_2$ and $\theta(i_1) \neq i_2$. As above, we have $d(i_1) = 0$ thus $\psi_0$ acts on $\multone_{1a}$ as $r_0$. We obtain from~\eqref{subequations:action_braid4}, omitting the idempotents,
\begin{align*}
(\prep\psi_0 \prep\psi_1)^2
&=
r_0 P_{\bar{2}1}(\multx_2, \multx_1)r_1 \alpha_2(\multx_1)r_0 P_{12}(\multx_2, \multx_1)r_1
\\
&= P_{\bar{2}1}(\multx_2, -\multx_1)r_0\alpha_2(\multx_2) r_1 r_0 P_{12}(\multx_2, \multx_1)r_1
\\
&=
P_{\bar{2}1}(\multx_2, -\multx_1)\alpha_2(\multx_2)P_{12}(-\multx_1, -\multx_2)r_0 r_1 r_0r_1,
\end{align*}
and
\begin{align*}
(\prep\psi_1 \prep\psi_0)^2
&=
P_{\bar{2}\bar{1}}(\multx_2, \multx_1)r_1 \alpha_2(\multx_1)r_0 P_{\bar{1}2}(\multx_2, \multx_1)r_1 r_0
\\
&=
P_{\bar{2}\bar{1}}(\multx_2, \multx_1)\alpha_2(\multx_2)r_1 r_0 P_{\bar{1}2}(\multx_2, \multx_1)r_1 r_0
\\
&=
P_{\bar{2}\bar{1}}(\multx_2, \multx_1)\alpha_2(\multx_2)P_{\bar{1}2}(\multx_1, -\multx_2)r_1 r_0 r_1 r_0,
\end{align*}
thus $(\prep\psi_0\prep\psi_1)^2 = (\prep\psi_1 \prep\psi_0)^2$ as desired, where we used $\bar{1} = 1$ and~\eqref{subequations:properties_P}. The case $\theta(i_1) \neq i_1$ and $\theta(i_2) = i_2$ is similar.

\begin{remark}
\label{remark:not_strong_assumption}
(See Remark~\ref{remark:condition_strong}.) Without condition~\eqref{equation:condition_strong}, we have to choose another, more complicated, relation~\eqref{relation:braid4}, if we want it to be compatible with the action on polynomials.
\end{remark}

\subsection{Case \texorpdfstring{$\gamma_{i_1} = 0 \neq \gamma_{i_2}$}{gammai1=0<>gammai2}}
\label{subsection:gammai1=0<>gammai2}

We want to prove that
\[
\left((\prep\psi_0 \prep\psi_1)^2 - (\prep\psi_1 \prep\psi_0)^2\right)  \multone_{\tuple i} = \gamma_{i_2} \frac{Q_{ i_2 i_1}(\prep y_1, -\prep y_2) - Q_{i_2 i_1}(\prep y_1, \prep y_2)}{\prep y_2}\prep\psi_0  \multone_{\tuple{i}},
\]
that is,
\[
\left((\prep\psi_0 \prep\psi_1)^2 - (\prep\psi_1 \prep\psi_0)^2\right) \multone_{\tuple i} =
\gamma_{i_2} \frac{Q_{i_2 i_1}(\multx_1, -\multx_2) - Q_{i_2 i_1}(\multx_1, \multx_2)}{\multx_2}\alpha_{i_1}(\multx_1)r_0 \multone_{\tuple i}.
\]

By~\eqref{equation:condition_gamma} we have $\theta(i_2) = i_2$. Note that $\gamma_{i_1} = 0 \neq \gamma_{i_2}$ implies $i_1 \neq i_2$. By~\eqref{subequations:action_braid4} we have, omitting the idempotents,
\begin{align*}
(\prep\psi_0 \prep\psi_1)^2
&=
\alpha_1(\multx_1) r_0 P_{\bar{2}1}(\multx_2, \multx_1)r_1 \gamma_2 \multx_1^{-1}(1 - r_0) P_{12}(\multx_2, \multx_1)r_1
\\
&=
\alpha_1(\multx_1)P_{\bar{2}1}(\multx_2, -\multx_1)\gamma_2 \multx_2^{-1}r_0 r_1 (1 - r_0)P_{12}(\multx_2, \multx_1)r_1
\\
&=
\alpha_1(\multx_1)P_{\bar{2}1}(\multx_2, -\multx_1)\gamma_2 \multx_2^{-1} \bigl[P_{12}(-\multx_1, \multx_2)r_0 r_1 - P_{12}(-\multx_1, -\multx_2)r_0 r_1 r_0\bigr]r_1
\\
&=
\gamma_2 \multx_2^{-1}\alpha_1(\multx_1)P_{\bar{2} 1}(\multx_2, -\multx_1)\bigl[P_{12}(-\multx_1, \multx_2)r_0r_1 - P_{12}(-\multx_1, -\multx_2)r_0 r_1 r_0\bigr] r_1,
\end{align*}
and
\begin{align*}
(\prep\psi_1 \prep\psi_0)^2
&=
P_{\bar{2}\bar{1}}(\multx_2, \multx_1)r_1 \gamma_2 \multx_1^{-1}(1 - r_0)  P_{\bar{1}2}(\multx_2, \multx_1)r_1 \alpha_1(\multx_1)r_0
\\
&=
P_{\bar{2}\bar{1}}(\multx_2, \multx_1)\gamma_2 \multx_2^{-1}r_1 (1 - r_0)P_{\bar{1}2}(\multx_2, \multx_1)\alpha_1(\multx_2)r_1 r_0
\\
&=
P_{\bar{2}\bar{1}}(\multx_2, \multx_1)\gamma_2 \multx_2^{-1} \bigl[P_{\bar{1}2}(\multx_1, \multx_2)r_1 - P_{\bar{1}2}(\multx_1, -\multx_2)r_1r_0\bigr]\alpha_1(\multx_2)r_1 r_0
\\
&=
P_{\bar{2}\bar{1}}(\multx_2, \multx_1)\gamma_2 \multx_2^{-1} \alpha_1(\multx_1) \bigl[P_{\bar{1}2}(\multx_1, \multx_2)r_1 - P_{\bar{1}2}(\multx_1, -\multx_2) r_1r_0\bigr]r_1 r_0
\\
&=
\gamma_2 \multx_2^{-1} \alpha_1(\multx_1)P_{\bar{2}\bar{1}}(\multx_2, \multx_1)\bigl[P_{\bar 1 2}(\multx_1, \multx_2) - P_{\bar 1 2}(\multx_1, -\multx_2)r_1 r_0r_1 \bigr]r_0.
\end{align*}
Thus, recalling $\bar{2} = 2$ and using the properties~\eqref{equation:Q_ij(uv)=Q_ji(vu)=Q_ij(-v-u)}, \eqref{equation:Qij=Qthetaithetaj}, \eqref{subequations:properties_P}, \eqref{equation:PP=Q} for the families $P$ and $Q$ we have
\begin{align*}
(\prep\psi_0\prep\psi_1)^2 - (\prep\psi_1 \prep\psi_0)^2
&=
\gamma_2 \multx_2^{-1} \alpha_1(\multx_1) \bigl[P_{\bar{2} 1}(\multx_2, -\multx_1) P_{12}(-\multx_1, \multx_2) - P_{\bar{2}\bar{1}}(\multx_2, \multx_1)P_{\bar 1 2}(\multx_1, \multx_2)\bigr]r_0
\\
&=
\gamma_2 \multx_2^{-1} \bigl[Q_{2 1}(\multx_2, -\multx_1) - Q_{2\bar 1}(\multx_2, \multx_1)\bigr]\alpha_1(\multx_1)r_0
\\
&=
\gamma_2 \multx_2^{-1} \bigl[Q_{21}(\multx_1, -\multx_2) - Q_{21}(\multx_1, \multx_2)\bigr]\alpha_1(\multx_1)r_0,
\end{align*}
as desired.

\subsection{Case \texorpdfstring{$\gamma_{i_1} \neq 0 = \gamma_{i_2}$}{gammai1<>0=gammai2}}

We want to prove that
\[
\left((\prep\psi_0 \prep\psi_1)^2 - (\prep\psi_1 \prep\psi_0)^2\right) \multone_{\tuple i} = 0.
\]

Similarly to~\textsection\ref{subsection:gammai1=0<>gammai2} we have $\theta(i_1) = i_1 \neq i_2$. By~\eqref{subequations:action_braid4} we have, omitting the idempotents,
\begin{align*}
(\prep\psi_0\prep\psi_1)^2
&=
\gamma_1 \multx_1^{-1}(1 - r_0) P_{\bar{2}1}(\multx_2,\multx_1)r_1 \alpha_2(\multx_1) r_0 P_{12}(\multx_2,\multx_1)r_1
\\
&=
\gamma_1 \multx_1^{-1}\left[P_{\bar 2 1}(\multx_2,\multx_1) - P_{\bar 2 1}(\multx_2,-\multx_1)r_0\right] \alpha_2(\multx_2)P_{12}(\multx_1,-\multx_2)r_1 r_0 r_1
\\
&=
\gamma_1 \multx_1^{-1}\alpha_2(\multx_2)\left[ P_{\bar 2 1}(\multx_2,\multx_1)P_{12}(\multx_1,-\multx_2) - P_{\bar 2 1}(\multx_2,-\multx_1)P_{12}(-\multx_1,-\multx_2)r_0\right]r_1 r_0 r_1
\\
&=
\gamma_1 \multx_1^{-1}\alpha_2(\multx_2)P_{\bar 2 1}(\multx_2,\multx_1)P_{12}(\multx_1,-\multx_2)(1 - r_0)r_1 r_0 r_1
\end{align*}
by~\eqref{subequations:properties_P}, and
\begin{align*}
(\prep\psi_1\prep\psi_0)^2
&=
P_{\bar{2} 1}(\multx_2,\multx_1)r_1 \alpha_2(\multx_1)r_0 P_{ 1 2}(\multx_2,\multx_1)r_1  \gamma_1 \multx_1^{-1}(1 - r_0)
\\
&=
\gamma_1 \multx_1^{-1} \alpha_2(\multx_2) P_{\bar 2 1}(\multx_2,\multx_1)P_{ 1 2}(\multx_1,-\multx_2) r_1 r_0 r_1(1 - r_0),
\end{align*}
Thus $(\prep\psi_0\prep\psi_1)^2 = (\prep\psi_1\prep\psi_0)^2$ as desired.

\subsection{Case \texorpdfstring{$\gamma_{i_1} \neq 0 \neq \gamma_{i_2}$}{gammai1<>0<>gammai2}}

We want to prove that (recalling from~\eqref{equation:def_Q} that $Q_{ii} = 0$)
\[
\left((\prep\psi_0 \prep\psi_1)^2 - (\prep\psi_1 \prep\psi_0)^2\right) \multone_{\tuple i} = \begin{cases}
\gamma_{i_2}  \frac{Q_{i_2 i_1}(\prep y_1, -\prep y_2) - Q_{i_2 i_1}(\prep y_1, \prep y_2)}{\prep y_1 \prep y_2} \left(\prep y_1 \prep\psi_0 - \gamma_{i_1}\right) \multone_{\tuple i},
&
\text{if } i_1 \neq i_2,
\\
0,
&\text{otherwise},
\end{cases}\]
that is, since $\psi_0$ acts on $\multone_{\tuple i}$ as $\gamma_{i_1}\multx_1^{-1} (1 - r_0)$ (recalling that $\theta(i_1) = i_1$ by~\eqref{equation:condition_gamma}),
\[
\left((\prep\psi_0 \prep\psi_1)^2 - (\prep\psi_1 \prep\psi_0)^2\right) \multone_{\tuple i} = \begin{cases}
\gamma_{i_1}\gamma_{i_2}  \frac{Q_{i_2 i_1}(\multx_1, \multx_2) - Q_{i_2 i_1}(\multx_1, -\multx_2)}{\multx_1 \multx_2} r_0 \multone_{\tuple i},
&
\text{if }  i_1 \neq i_2,
\\
0,
&\text{otherwise}.
\end{cases}\]

The next result is an easy calculation.

\begin{lemma}
\label{lemma:commutation_rule_partial0}
Let $P$ be a polynomial in $\multx_1, \multx_2$ and let $w \in \langle r_0, r_1\rangle$. Then
\[
\multx_1^{-1}(1 - r_0) P w-  Pw \multx_1^{-1}(1 - r_0) = \left(\multx_1^{-1} - \actpol{w}{\multx}_1^{-1}\right)Pw + \actpol{w}{\multx}_1^{-1} P wr_0 - \multx_1^{-1} \prescript{r_0}{}{P} r_0w\,,
%\\= \multx_1^{-1}(Pw - \prescript{r_0}{}{P} r_0 w) - P \actpol{w}{\multx_1}^{-1}(w - wr_0)
\]
inside $\mathrm{End}_K(K[x,\beta])$.
\end{lemma}

By~\eqref{equation:condition_gamma} we have $\theta(i_2) = i_2$. If $i_1 \neq i_2$ we obtain from~\eqref{subequations:action_braid4}
\begin{align*}
\prep\psi_1 \prep\psi_0 \prep\psi_1  \multone_{12}
&=
P_{2 1}(\multx_2, \multx_1)r_1 \gamma_2 \multx_1^{-1}(1 - r_0) P_{12}(\multx_2, \multx_1)r_1 \multone_{12}
\\
&=
 \gamma_2  \multx_2^{-1}  P_{2 1}(\multx_2, \multx_1)r_1 (1 - r_0) P_{12}(\multx_2, \multx_1)r_1 \multone_{12}
\\
&=
 \gamma_2  \multx_2^{-1}  P_{2 1}(\multx_2, \multx_1)r_1 \bigl[P_{12}(\multx_2, \multx_1) -  P_{12}(\multx_2, -\multx_1)r_0\bigr]r_1 \multone_{12}
 \\
&=
 \gamma_2  \multx_2^{-1}  P_{2 1}(\multx_2, \multx_1)\bigl[P_{12}(\multx_1, \multx_2)r_1 -  P_{12}(\multx_1, -\multx_2)r_1 r_0\bigr]r_1 \multone_{12}
  \\
&=
 \gamma_2  \multx_2^{-1}  P_{2 1}(\multx_2, \multx_1)\bigl[P_{12}(\multx_1, \multx_2) -  P_{12}(\multx_1, -\multx_2)r_1 r_0r_1\bigr] \multone_{12}.
\end{align*}
%minuscules modifs ici
Since $\prep\psi_0  \multone_{12} = \gamma_1 \multx_1^{-1}(1 - r_0) \multone_{12}$, we can apply Lemma~\ref{lemma:commutation_rule_partial0} for the two above summands. We obtain that second summand will vanish in $\left((\prep\psi_0 \prep\psi_1)^2 - (\prep\psi_1 \prep\psi_0)^2\right) \multone_{12}$ since $x_1^{-1} \in K(x_1)$ is invariant under $r_1 r_0 r_1$ and $P_{21}(x_2, x_1) P_{12}(x_1, -x_2) \in K[x_1,x_2]$ is invariant under $r_0$ by~\eqref{subequations:properties_P}. Thus, we only consider the first summand, which is equal to $\gamma_2 \multx_2^{-1} Q_{21}(\multx_2, \multx_1)$, and we obtain, omitting the idempotents and using~\eqref{equation:Q_ij(uv)=Q_ji(vu)=Q_ij(-v-u)} and~\eqref{equation:Qij=Qthetaithetaj},
\begin{align*}
(\prep\psi_0 \prep\psi_1)^2 - (\prep\psi_1 \prep\psi_0)^2
&=
\gamma_1\gamma_2 \multx_1^{-1}\multx_2^{-1} \bigl[Q_{21}(\multx_2, \multx_1) - Q_{21}(\multx_2, -\multx_1)\bigr]r_0
\\
&=
\gamma_1\gamma_2 \multx_1^{-1}\multx_2^{-1} \bigl[Q_{21}(\multx_1, \multx_2) - Q_{21}(\multx_1, -\multx_2)\bigr]r_0,
\end{align*}
as desired.

Finally, assume that $i_1 = i_2$. We have
\begin{multline*}
(\prep\psi_0\prep\psi_1)^2 - (\prep\psi_0 \prep\psi_1)^2 = \gamma_1^2\left[\multx_1^{-1}(1 - r_0) (\multx_1 - \multx_2)^{-1}(r_1 - 1)\multx_1^{-1}(1 - r_0)(\multx_1 - \multx_2)^{-1}\right.
\\
\left.- (\multx_1 - \multx_2)^{-1}(r_1 - 1)\multx_1^{-1}(1 - r_0)(\multx_1 - \multx_2)^{-1} \multx_1^{-1}(1 - r_0)\right] = 0,
\end{multline*}
since this is just the braid relation for the divided difference operators $\partial_0 \coloneqq \multx_1^{-1}(1 - r_0)$ and $\partial_1 \coloneqq (\multx_1 - \multx_2)^{-1}(r_1 - 1)$ (see~\cite{bgg,demazure}).

\end{document}